\documentclass[preprint,12pt]{elsarticle}
\usepackage{amsfonts,mathtools}
\usepackage{hyperref}
\usepackage{adjustbox}
\usepackage{ulem}
\usepackage{cancel}
\usepackage{todonotes}
\usepackage{bm}
\usepackage{longtable} 
\usepackage{commath}
\usepackage{multirow}
\usepackage{stmaryrd}
\usepackage{amsmath}
\DeclareMathOperator*{\argmin}{arg\,min}

\newdefinition{rmk}{Remark}
\newproof{proof}{Proof}
\newcommand{\jump}[1]{\llbracket #1 \rrbracket}

\newtheorem{theorem}{Theorem}[section]
\newtheorem{remark}[theorem]{Remark}

\usepackage[percent]{overpic}
\usepackage[font=small,labelfont=bf]{caption}
\newcommand{\bx}{\bm{x}}

\newcommand{\balpha}{\bm{\alpha}}

\newcommand{\TN}{\text{TN}}
\newcommand{\NN}{\text{NN}}
\newcommand{\exact}{\text{e}}

\newcommand{\WTN}{\text{WTN}}
\newcommand{\stf}{\text{SF}}
\newcommand{\F}{\text{F}}
\newcommand{\DRM}{\text{DRM}}
\newcommand{\pou}{\text{PoU}}

\newcommand{\hatu}{\hat{u}}
\newcommand{\R}{\text{R}}
\newcommand{\M}{\text{M}}
\usepackage{algorithm}
\usepackage{algpseudocode}

\journal{Journal of Computational Physics}
\makeatletter
\def\ps@pprintTitle{%
 \let\@oddhead\@empty
 \let\@evenhead\@empty
 \let\@oddfoot\@empty
 \let\@evenfoot\@empty
}
\makeatother 

\begin{document}

\begin{frontmatter}

\title{Weak TransNet: A Petrov-Galerkin based neural network method for solving elliptic PDEs}
\author[UHaddress]{Zhihang Xu}
\ead{zxu29@central.uh.edu}
\address[UHaddress]{Department of Mathematics, University of Houston, TX 77204, USA}
\author[UHaddress]{Min Wang\corref{mycorrespondingauthor}}
\cortext[mycorrespondingauthor]{Corresponding author}
\ead{mwang55@central.uh.edu}
\author[SCaddress]{Zhu Wang}
\address[SCaddress]{Department of Mathematics, University of South Carolina, Columbia, SC 29208, USA}
\ead{wangzhu@math.sc.edu} 

\begin{abstract}
While deep learning has achieved remarkable success in solving partial differential equations (PDEs), it still faces significant challenges, particularly when the PDE solutions have low regularity or singularities. To address these issues, we propose the Weak TransNet (WTN) method, based on a Petrov-Galerkin formulation, for solving elliptic PDEs in this work, though its framework may extend to other classes of equations. 
Specifically, the neural feature space defined by TransNet~\cite{Zhang2023TransNetTN} is used as the trial space, while the test space is composed of radial basis functions. 
Since the solution is expressed as a linear combination of trial functions, the coefficients can be determined by minimizing the weak PDE residual  via least squares.  
Thus, this approach could help mitigate the challenges of non-convexity and ill-conditioning that often arise in neural network training. 
Furthermore, the WTN method is extended to handle problems whose solutions exhibit multiscale features or possess sharp variations.
Several numerical experiments are presented to demonstrate the robustness and efficiency of the proposed methods. 
\end{abstract}



\begin{keyword}
Weak formulation\sep Petrov-Galerkin formulation \sep TransNet \sep  physics-informed learning 
\end{keyword}
\end{frontmatter}

\section{Introduction}
Many scientific research and engineering applications require numerical solutions to PDEs, for which traditional numerical methods, such as the finite element method~(FEM), have been established. For example, in groundwater hydrology and oil and gas reservoir simulations,  classical Darcy flow~\cite{whitaker1986flow} serves as a fundamental model in subsurface fluid dynamics that needs to be effectively simulated. When permeability varies across multiple spatial scales in heterogeneous porous media, the generalized multiscale finite element method (GMsFEM)~\cite{efendiev2013generalized} is designed to efficiently capture the multiscale structure characteristic in numerical simulations of Darcy flow. Although these traditional numerical methods have been successful in scientific computing, they are mesh-based methods and thus suffer from the \textit{curse of dimensionality}. That is, the number of required grid points grows exponentially with respect to the dimensionality $d$ of the problem, leading to prohibitive memory consumption and computational costs. Therefore, meshless approaches, including neural networks \cite{yang2016data,zhu2018bayesian,zhu2019physics,raissi2019physics} and kernel methods \cite{chen2021solving}, are emerging to solve PDEs, especially for high-dimensional problems. 

Among these approaches, deep learning-based methods approximate the solution of a PDE, $u(\bx)$ (or $u(\bx,t)$), using a neural network (NN) model,  $u_{\NN}(\bx;\Theta)$ (or $u_{\NN}(\bx,t;\Theta)$), with its parameters $\Theta$ learned by minimizing a loss function during training \footnote{Hereafter, we omit $\Theta$ in $u_{\NN}(\bx;\Theta)$ for brevity unless explicitly needed.}. These methods can be categorized into two classes: \textit{data-driven}, where training is supervised and the objective function measures the discrepancy between $u_{\NN}$ and $u$ on a set of labeled data without incorporating any physics~\cite{mo2019deep}; 
and \textit{physics-driven}, 
where the loss enforces how well $u_{\NN}$ satisfies the governing equations. 
Within the later category, 
by embedding the governing equations into the learning process, the physics-informed neural network (PINN)~\cite{raissi2019physics,cai2021physics} has
shown promising results in solving various PDEs and its associated inverse problems.

Despite the remarkable success of deep learning in solving PDEs, several critical limitations remain. 
First, the accuracy of NN-based methods depends on the regularity of the solution, as indicated by the NN approximation theory~\cite{devore2021neural}. The regularity of the PDE solution is strongly influenced by the smoothness of the source term, boundary conditions, and domain geometry.
When the domain has geometric singularities, such as corners or fractures, the solution may exhibit low regularity, losing smoothness in certain local regions.
Meanwhile, due to its inherent spectral bias or F-principle, the NN model favors approximating smooth, low-frequency functions~\cite{rahaman2019spectral,xu2019frequency}. Therefore, finding the NN approximation for problems involving complex or high-frequency solutions becomes more challenging. 

To overcome this limitation, one line of research focuses on minimizing the PDE's residual in different norms within the NN-based framework, which is equivalent to solving the PDEs weakly~\cite{sirignano2018dgm, yu2018deep, zhu2019physics, zang2020weak, kharazmi2021hp, sheng2021pfnn,  shang2023randomized, liu2024subspace}. 
In particular, the deep Ritz method (DRM)~\cite{yu2018deep} is based on the energy functional that corresponds to the weak form of the PDEs. 
The deep Nitsche method~\cite{liao2019deep} employs Nitsche's variational formulation to enforce the essential boundary condition. 
The loss function of weak adversarial networks (WAN)~\cite{zang2020weak} is defined by the operator norm induced by the weak form of the PDEs. 
The penalty-free neural network (PFNN)~\cite{sheng2021pfnn} also adopts the weak form of the original problem to avoid higher-order derivative evaluations.  
Besides, there is a growing body of work~\cite{wang2020deep,xu2021weak} that combines both physics-driven and data-driven methods to leverage their complementary advantages.  
Another line of research focuses on problems containing high-frequency components, where vanilla NN-based methods often struggle to capture fine-scale structures due to the aforementioned F-principle. To address this limitation, several approaches have been proposed to enhance the learning of high-frequency information. Among them, the multi-scale deep neural network (MscaleDNN)~\cite{liu2020multi,wang2020multi}, inspired by the idea of radial scaling in the frequency domain, introduces multiple scaling factors to improve the approximation of high-frequency features. 

Another limitation lies in the computational cost of solving the optimization problem.
Gradient-based methods such as gradient descent method and quasi-Newton algorithms are often computationally expensive~\cite{raissi2017physics, markidis2021old}, primarily due to the nonconvex nature of the objective function.
Although stochastic gradient descent (SGD) and its variants can improve efficiency, achieving highly accurate neural network approximations remains a significant challenge. 

To address this limitation, recent studies~\cite{dong2021local, chen2022bridging,Zhang2023TransNetTN, shang2023randomized,wang2024extreme, lu2025multiple} have proposed the use of randomized neural networks, extreme learning method, or TransNet-based methods. 
In these approaches, the hidden-layer parameters of the NN model are frozen, defining a predetermined feature space. Training reduces to solving for the parameters in the output layer, which linearly combine these features to approximate the solution. 
For linear PDEs, they can be efficiently determined using a least-squares algorithm. For nonlinear PDEs, the solution requires a nonlinear iterative solver combined with least-squares minimization.
By doing this, one can avoid nonconvex optimization and significantly reduce training difficulties. 
Specifically, 
Chen et al. proposed random feature method (RFM)~\cite{chen2022bridging}, which uses random feature functions to approximate PDE solutions. The loss function is taken to be the strong form of the PDE residual at collocation points, with the boundary condition incorporated as a penalty term. 
Local extreme learning machines~(LocELM)~\cite{wang2024extreme} combined the ideas of extreme learning machine (ELM)~\cite{ding2015extreme} and domain decomposition.  
In this paper, we propose a weak TransNet (WTN) method to solve PDEs in their weak form. 
In particular, 
the trial basis functions are defined by those generated from TransNet, which form a neural feature space~\cite{Zhang2023TransNetTN}.
To ensure accuracy while minimizing computational costs, radial basis functions (RBFs) are selected as the test functions due to their inherent local support.
For problems with multi-scale features, we combine the Fourier feature mapping and WTN, proposing a Fourier-Weak TransNet (F-WTN) method. 
In addition, the current framework can be seamlessly integrated with the partition of unity (PoU) method to enhance representation capability.
As a first step toward addressing the multiscale challenges associated in the Darcy flow problems, 
we consider several representative cases in our numerical examples. Through them, we demonstrate the effectiveness and robustness of the proposed methods.

Our contributions in this work are summarized as follows:
\begin{itemize}
    \item We propose a novel WTN method to solve PDEs based on their weak formulation. The method leverages the TransNet architecture to construct trial basis functions within a neural feature space, enabling an efficient and flexible training-free approximation framework.
    \item Building upon WTN, we develop several extensions to enhance its capability and applicability. Specifically, we incorporate Fourier feature mapping to capture multiscale behaviors, and apply the PoU method to improve the representation power and scalability of the model. 
\end{itemize}

The rest of the paper is organized as follows. In Sec.~\ref{sec:method}, we first provide a brief review of different formulations employed by deep learning methods for solving PDEs and TransNet. 
In Sec.~\ref{sec:weak_form} and Sec.~\ref{sec:wtn-extensions}, we detail the WTN approach and its extensions, especially, including F-WTN and PoU-WTN.
In Sec.~\ref{sec:num}, we provide several numerical examples to demonstrate the effectiveness of our proposed method. 
Finally, Sec.~\ref{sec:conclusion} presents the conclusion.

\section{Problem setting and TransNet}
\label{sec:method}
Consider the following equation over a bounded and connected polygon domain $\Omega\subset \mathbb{R}^d$ with boundary $\partial \Omega$:
\begin{subequations}
\label{eq:problem}
\begin{align}
&\mathcal{L} [u(\bx)] = f(\bx)\,,\quad \bx\in \Omega
\,,
\label{eq:pde}
\\
&\mathcal{B}[u(\bx)] = g(\bx)\,,\quad \bx\in \partial \Omega\,,
\label{eq:bc}
\end{align}
\end{subequations}
\noindent
where $\bx= [x_1,\ldots, x_d]^{\top}$ are independent variables \footnote{We focus on elliptic equations in this work, which are independent of time. If the problem is time dependent, both spatial and time variables can be included in $\bx$.}, and $u(\bx): \Omega \to \mathbb{R}$ is the state variable. In (\ref{eq:pde}), $\mathcal{L}$ is a differential operator that characterizes the equation and $f(\bx)$ is the source term; and in (\ref{eq:bc}), $\mathcal{B}$ is a boundary operator and $g(\bx)$ specifies the given boundary condition.  
For instance, when considering Darcy flow with Dirichlet boundary condition, $\mathcal{L}[u]\coloneqq -\nabla \cdot (\kappa(\bx) \nabla u )$ and $\mathcal{B}[u] \coloneqq u$, where $\kappa(\bx)$ denotes the spatially varying permeability field. 
The \textit{strong-form} PDE residual of this problem is defined as 
\[
\mathcal{R}[u(\bx)]\coloneqq \mathcal{L}[u(\bx)] - f(\bx). 
\]

To seek an approximation of $\hat{u}(\bx)$ of $u(\bx)$, one often uses the ansatz that $\hat{u}(\bx;\balpha): = \sum_{j=0}^M \alpha_j \phi_j(\bx)$,
where $\{\phi_j\}_{j=0}^M$ is a set of \textit{trial} basis functions that span the approximation space $\mathcal{U}$, and $\balpha = [\alpha_0, \ldots, \alpha_M]^\top$ is the coefficient vector to be determined. 
To find $\balpha$, the weighted residual method can be applied. It enforces that $\mathcal{R}$ vanishes in a weighted sense: 
\begin{equation}
\int_{\Omega} \psi_i(\bx) \mathcal{R}[\hatu(\bx; \balpha)] \dif \bx=0\,,\quad i=1,\ldots,N\,,
\label{eq:weighted_residual}
\end{equation}
in which the weights $\{\psi_i(\bx)\}_{i=1}^{N}$ are referred to as the \textit{test} functions. Different choices of test functions lead to different methods: When $\psi_i = \phi_i$, it is known as the Galerkin method.  
When they are different, it is the Petrov-Galerkin method. Specifically, if $\psi_i = \frac{\partial \mathcal{R}[\hat{u}]}{\partial \alpha_i}$, the approach is known as the least-squares method~\cite{bochev2009least}. 
When $\psi_i = \delta(\bx - \bx_i)$, (\ref{eq:weighted_residual}) becomes $\mathcal{R}[\hat{u}(\bx_i)] = 0$ for $i= 1,\ldots,N$, which is referred as the collocation method.

Recently developed neural network-based PDE solvers can also fit within this framework, in which $\hat{u}(\bx)$ is modeled as a neural network and the associated trial space $\mathcal{U}$ is spanned by functions corresponding to the output of neurons in the last hidden layer of the neural network. Various types of test functions have been explored, including but not limited to: Dirac delta functions $\psi_i = \delta(\bx - \bx_i)$ in PINN~\cite{raissi2019physics}; pre-trained neural networks in a Galerkin framework~\cite{liu2024subspace}; finite element basis functions in~\cite{shang2023randomized};  piecewise polynomials in the variational PINNs~(VPINNs) method~\cite{kharazmi2019variational} and its $hp$-VPINN extension~\cite{kharazmi2021hp}, and also in the randomized neural networks with Petrov-Galerkin method (RNN-PG)~\cite{shang2024randomized}; and Lagrange polynomials in~\cite{roop2024randomized}.

\subsection{TransNet}
\label{sec:transnet}
TransNet \cite{Zhang2023TransNetTN} uses a single-hidden-layer fully connected neural network comprising $M$ neurons to approximate $\hat{u}(\bx)$. The approximation, denoted by $u_{\TN}: \Omega \to \mathbb{R}$, is defined as
\[
u_{\TN}(\bx) : =  \sum_{j=1}^M \alpha_j \sigma(\bm{w}_j^\top \bx + b_j) + \alpha_0\,,
\label{eq:TransNet0}
\]
where $\bm{w}_j\in \mathbb{R}^d$ and $b_j\in\mathbb{R}$ are the weights and bias parameters of the $j$-th neuron, respectively, $\sigma(\cdot)$ is a nonlinear activation function, and $\{\alpha_j\}_{j=0}^M$ are undetermined coefficients. 
Each hidden neuron output, $\sigma(\bm{w}_j^\top \bx+ b_j)$, can be regarded as a neural function, denoted by $\phi_j\coloneqq \sigma(\bm{w}_j^\top \bx+ b_j)$ for $j=1,\ldots,M$. We further let $\phi_0 =1$, then $u_{\TN}$ can be rewritten as 
\begin{equation}
u_{\TN}(\bx) = \sum_{j=0}^M \alpha_j \phi_j(\bx)\,.
\label{eq:TransNet}
\end{equation}
However, with common choices of activation functions such as ReLU or Tanh, $\{\phi_j\}_{j=0}^M$ are typically \textit{non-orthogonal} and \textit{globally-supported}. Nevertheless, we name them \textit{neural basis functions} and define the associated \textit{neural feature space} as  
\[
\mathcal{U}_{\TN} \coloneqq \text{span}\{\phi_0, \ldots, \phi_M\}\,.  
\] 

TransNet introduces a novel way to ensure neural basis functions are nearly uniformly distributed in $\Omega$. 
To achieve this, $\bm{w}_j^\top\bx + b_j$ is first re-parameterized as 
$\gamma_j(\bm{a}_j^\top\bx+r_j)\,,$
where $\bm{a}_j$ is a unit hyperparameter vector, the shape parameter $\gamma_j > 0$ and $r_j$ are two additional scalar hyperparameters. Geometrically, $(\bm{a}_j, r_j)$ determines the location of the hyperplane $\bm{w}_j^\top\bx + b_j=0$ and $\gamma_j$ describes its steepness. According to \cite[Theorem 1]{Zhang2023TransNetTN}, if $\{\bm{a}_j\}_{j=1}^M$ are 
i.i.d. random vectors that are uniformly distributed on the $d$-dimensional unit sphere,
and $\{r_j\}_{j=1}^M$ are i.i.d. random variables obeying a uniform distribution in $[0,1]$, then the neural basis functions are \textit{uniformly distributed} in the unit ball. 
It has been generalized in \cite[Theorem 2]{lu2025multiple} to construct uniformly distributed neural basis functions on a $d$-dimensional ball centered at a point $\bx_c\in \mathbb{R}^d$ of a radius $R >0$.
Meanwhile, shape parameters $\{\gamma_j\}_{j=1}^M$ are critical in controlling the expressivity of the neural feature space $\mathcal{U}_{\TN}$. 
The higher the value of $\gamma_j$, the \textit{steeper} the basis $\phi_j$ becomes. 
In \cite{Zhang2023TransNetTN}, all the shape parameters are set to the same value $\gamma$, which is tuned by minimizing the projection error of a set of auxiliary functions, generated by Gaussian random fields, onto the neural feature space. 
In~\cite{lu2025multiple}, an empirical formula, 
$
\gamma \approx C {M^{1/d}}{R^{-1}},
$
is suggested, where $C$ is a constant related to the settings of the problem and $R$ is the radius of the ball. 

A key difference between TransNet and other methods, such as RFM~\cite{chen2022bridging} and locELM~\cite{wang2024extreme}, lies in the interpretability of its hyperparameters $\{\bm{a_j}, r_j, \gamma_j\}$ and, consequently, $\{\bm{w}_j, b_j\}$.
Once selected, the neural basis functions become fixed. Hence, determining $u_{\TN}(\bx)$ reduces to solving for the coefficient vector $\balpha=[\alpha_0, \ldots, \alpha_M]^\top$. 
When solving linear PDEs, enforcing $u_{\TN}(\bx)$ to satisfy the PDE yields a linear system for $\balpha$, 
since $u_{\TN}(\bx)$ linearly depends on $\balpha$.
This is elaborated in~\ref{sec:strong_loss} and \ref{sec:drm_loss}, where the loss functions are set as the strong-form residual and the Ritz energy, respectively. 
This contrasts with PINNs, in which all learnable parameters $\{\alpha_j, \bm{w}_j, b_j\}_{j=0}^M$ are optimized via gradient-based algorithms, which has been shown in~\cite{krishnapriyan2021characterizing} to suffer from ill-conditioning and convergence difficulties. On the other hand, TransNet reduces the reliance of neural network-based PDE solvers on extensive training and mitigates their sensitivity to the initialization of learnable parameters. As a result, it improves the stability and reproducibility of the numerical results. 
Furthermore, by adjusting shape parameters, the resulting neural feature space can represent solutions with high frequency or sharp gradients. This property can potentially alleviate the spectral bias and consequently improve the accuracy of the approximate solution. 
Therefore, we adopt the neural basis functions generated by TransNet as trial functions in this work, but, unlike the collocation method used in \cite{Zhang2023TransNetTN,lu2025multiple}, we employ the Petrov-Galerkin approach to effectively handle problems with weak regularization. 

\section{Weak TransNet method}
\label{sec:weak_form}
Setting the trial space $\mathcal{U}_{\TN}= \text{span}\{\phi_0, \phi_1,\ldots, \phi_{M} \}$ and the test space $\mathcal{V} \coloneqq \text{span}\{\psi_1,\ldots, \psi_{N}\}$, we can recast the equation (\ref{eq:weighted_residual}) into the following form: 
to find $\hat{u}(\bx) \in \mathcal{U}_{\TN}$, satisfying the Petrov-Galerkin formulation
\begin{equation}
\label{eq:pg}
a(\hat{u}, \psi_i) = l(\psi_i)\,,\quad \forall \psi_i \in \mathcal{V}\,,
\end{equation}
where $l(\psi_i) = \int_{\Omega}f(\bx)\psi_i(\bx)\dif \bx$ and the bilinear form $a(\hat{u},\psi_i)$ is derived from $\int_\Omega \psi_i(\bx) \mathcal{L}[u(\bx)]\dif \bx$ after certain integration by parts, along with the enforcement of boundary conditions.

Obviously, evaluations of integrals over the domain $\Omega$ are required in this formulation. Since the trial basis functions are globally defined, we choose test functions with only local support or behave like locally supported functions in order to reduce the computational cost. To this end, we select the \textit{RBFs}~\cite{buhmann2000radial} with a diagonal covariance matrix:
\begin{equation}
    \psi(\bx;\bm{\mu},\bm{\sigma}) = \frac{1}{(2\pi)^{d/2}\prod_{k=1}^d \sigma_k} \exp\left( -\frac{1}{2}\sum_{k=1}^d\frac{(x_k - \mu_k)^2}{\sigma_k^2} \right),
\label{eq:test_function}
\end{equation}
where $\bm{\mu} = [\mu_1, \ldots, \mu_d]^\top$ is the mean vector, and $\bm{\sigma} = [\sigma_1,\ldots, \sigma_d]^{\top}$ is the standard deviation vector. 
Consider evaluating an integral involving the test function: $\int_{\Omega} \widetilde{h}(\bx) \psi(\bx)\dif \bx$, where $\widetilde{h}(\bx)$ denotes an arbitrary function.
Due to the bell-shaped profile of the test function, the integration region can be reduced to $\Omega_{\psi} \coloneqq \left(\prod_{k=1}^d[\mu_k-N_l \sigma_k, \mu_k + N_l\sigma_k] \right)\cap \Omega$ where $N_l$ is a small integer (typically not greater than 10).
Since the volume of this region is $\prod_{k=1}^d (2N_l\sigma_k)$, the ratio of computational saving for evaluating integrations is about  
$1 - [{\prod_{k=1}^d(2N_l \sigma_k)}]/{|\Omega|}$ for achieving a desired numerical accuracy. For instance, if we integrate over $\Omega= (0,1)^2$, when $N_l=10$ and $\sigma_k =0.03$, for $k=1,\ldots, d$, are taken, the actual computational saving in time is about 64\%. 

To generate $N$ test functions, we vary $\bm{\mu}$ while keeping the same $\bm{\sigma}$, and obtain
\[\psi_i(\bx) = \psi(\bx;\bm{\mu}^{(i)},\bm{\sigma})\,, \text{ for } i = 1, \ldots, N.
\]
Since system (\ref{eq:pg}) may lack a unique solution, we find a least-squares solution, $u_{\WTN}(\bx)$, that minimizes the following \textit{weak residual} loss: 
\begin{equation}
\mathfrak{L}_{\text{weak}}\coloneqq \sum_{i=1}^N \|a(u_{\WTN}, \psi_i) - l(\psi_i)\|^2\,,
\label{eq:weak_residual}
\end{equation}
subject to the given boundary condition. 
The constraint on the boundary can be imposed either softly (e.g., via penalty terms in the loss function) or hardly (e.g., by modifying the neural network's architecture to inherently satisfy constraints). Correspondingly, we propose the weak TransNet method (WTN): to find $u_{\WTN}(\bx)\in \mathcal{U}_{\TN}$ that minimizes the loss function 

 \begin{equation}
\mathfrak{L} =\mathfrak{L}_{\text{weak}}
 + \underbrace{\beta\int_{\partial \Omega}\left\| \mathcal{B}[u_{\WTN}(\bx)]-g(\bx) \right\|^2\dif s}_{\text{boundary loss}}\,,
 \label{eq:weak_loss}
\end{equation} 
where $\beta> 0$ for soft constraints, and $\beta=0$ for hard constraints. 
In the former case, we uniformly generate boundary samples $S_{\partial \Omega} =\{ \bx_{\partial \Omega}^{(m)}\}_{m=1}^{N_{\partial \Omega}} \subset \partial \Omega$, and then 
approximate the boundary loss in \eqref{eq:weak_loss} by the Monte Carlo method, leading to the following empirical loss:
\begin{equation}
\mathfrak{L}_{\text{E}} =\mathfrak{L}_{\text{weak}}
 + \underbrace{\frac{\beta|\partial \Omega|}{ N_{\partial \Omega}} \sum_{m=1}^{N_{\partial \Omega}}\left(\mathcal{B}[u_{\WTN}(\bx_{\partial \Omega}^{(m)})]-g(\bx_{\partial \Omega}^{(m)})\right)^2 }_{\text{boundary loss}}\,.
 \label{eq:weak_loss_empical}
\end{equation} 
\paragraph{Linear case}
Suppose the operators $\mathcal{L}[u]$ and $\mathcal{B}[u]$ are both linear, using the ansatz $u_{\WTN} = \sum_{j=0}^{M}\alpha_j \phi_j$ in~\eqref{eq:weak_loss_empical}, 
we solve for $\balpha$ from the following minimization problem:  
\[
\min_{\balpha} \| \bm{L} \balpha - \bm{r}\|_2^2\,,
\]
where $
\bm{L} = \begin{bmatrix}
    \bm{A} \\
    \widetilde{\beta}\bm{B}
\end{bmatrix}$,  $
\bm{r} = \begin{bmatrix}
    \bm{f} \\
    \widetilde{\beta} \bm{g}
\end{bmatrix}$,
and $\widetilde{\beta} = \sqrt{\frac{\beta |\partial \Omega|}{N_{\partial \Omega}}}$ is the adjusted weight. 
Here $\bm{A}\in \mathbb{R}^{N\times{(M+1)}}$ with 
$\bm{A}_{ij}= a(\phi_j, \psi_i)$, and 
$\bm{f}\in \mathbb{R}^N$ with $\bm{f}_i = l(\psi_i)$, $\bm{B}\in \mathbb{R}^{N_{\partial\Omega}\times (M+1)}$ with 
$\bm{B}_{mj} = \mathcal{B}[\phi_j(\bx^{(m)})]$, 
and 
$\bm{g}\in \mathbb{R}^{N_{\partial \Omega}}$ with $\bm{g}_m = g(\bx^{(m)})$.  
A detailed description of WTN in this case is summarized in Alg.~\ref{alg:wtn}.


\begin{algorithm}[!htb]
\caption{A Weak TransNet (WTN) method for PDEs}\label{alg:wtn}
\begin{algorithmic}[1]
\Require{Number of trial basis $M$, number of test basis $N$, PDE form \eqref{eq:problem}, collocation points on the boundary $S_{\partial \Omega}\coloneqq\{\bx_{\partial \Omega}^{(m)}\}_{m=1}^{N_{\partial\Omega}}$, shape parameter $\gamma$, adjusted boundary weight $\widetilde{\beta}$}
\Ensure{$u_{\WTN}(\bx; \balpha^{\star})$}
\State{Construct the trial basis $\{\phi_j(\bx)\}_{j=0}^M$}
\State{Construct the test basis $\{\psi_i(\bx; \bm{\mu}, \bm{\sigma}\}_{i=1}^{N}$}
\For{$i=1,\ldots, N$}
\State{Identify the local support $\Omega_{\psi_i}$of $\psi_i$}
\State{Compute the the right-hand-side vector entry locally $\bm{f}_i = l(\psi_i)_{\Omega_{\psi_i}\cap \Omega}$}
\For{$j=0,\ldots,M$}
\State{Compute the stiffness matrix entry locally $\bm{A}_{ij} = a(\phi_j,\psi_i)_{\Omega_{\psi_i}\cap \Omega}$}
\EndFor
\EndFor
\State{Compute the boundary matrix $\bm{B}\in \mathbb{R}^{N_{\partial \Omega}\times (M+1)}$ with entry $\bm{B}_{mj} = [\mathcal{B}[\phi_j(\bx^{(m)})]]$ and the boundary condition vector $\bm{g}\in \mathbb{R}^{N_{\partial \Omega}}$ with entry $\bm{g}_m = g(\bx^{(m)})$}
\State{Solve the LS system:
\[
\balpha^{\star} = \argmin_{\balpha} \left\| \begin{bmatrix}
    \bm{A} \\
    \widetilde{\beta}\bm{B} 
\end{bmatrix} \balpha - \begin{bmatrix}
    \bm{f} \\
    \widetilde{\beta}\bm{g}
\end{bmatrix} \right\|_2^2
\]
}
\State{Obtain the approximate solution $u_{\WTN}(\bx;\balpha^{\star})$.}
\end{algorithmic}
\end{algorithm} 

\paragraph{Nonlinear case} 
If the operator $\mathcal{L}[u]$ and/or the boundary $\mathcal{B}[u]$ are nonlinear, one can employ existing nonlinear iterative solvers such as Picard's iteration to solve it.
Within each iteration, the nonlinear problem is linearized around the current approximate solution.
Consequently,
a linear system for $\balpha$ can be formulated and solved using the same procedure as in the linear case.


\subsection{Comparison of numerical quadrature}

Because WTN involves inner products, the choice of numerical quadrature affects the accuracy of numerical approximations. Take the Poisson equation with the homoegeneous Dirichlet boundary condition for example, that is, $\mathcal{L}[u] = -\Delta u$, $g=0$ in (\ref{eq:weighted_residual}), and the computation domain be $\Omega = (0,1)^2$. 
The exact solution is prescribed as 
$u(x,y) = \sin(\pi x)\sin(\pi y)$. 
To impose the boundary condition, we adopt the hard constraint approach by multiplying each trial basis function $\phi_j$ by the function $h(x,y) = x(1-x)y(1-y)$~\cite{liu2024subspace}. 
Although this problem can be effectively solved using the strong form loss~\eqref{eq:stf_loss}, we use the WTN method here, but focusing on comparing two numerical quadratures: Monte Carlo (MC) and composite Simpson’s rule (SIMP).

\begin{figure}[!htb]
    \centering
    \includegraphics[width=0.5\linewidth]{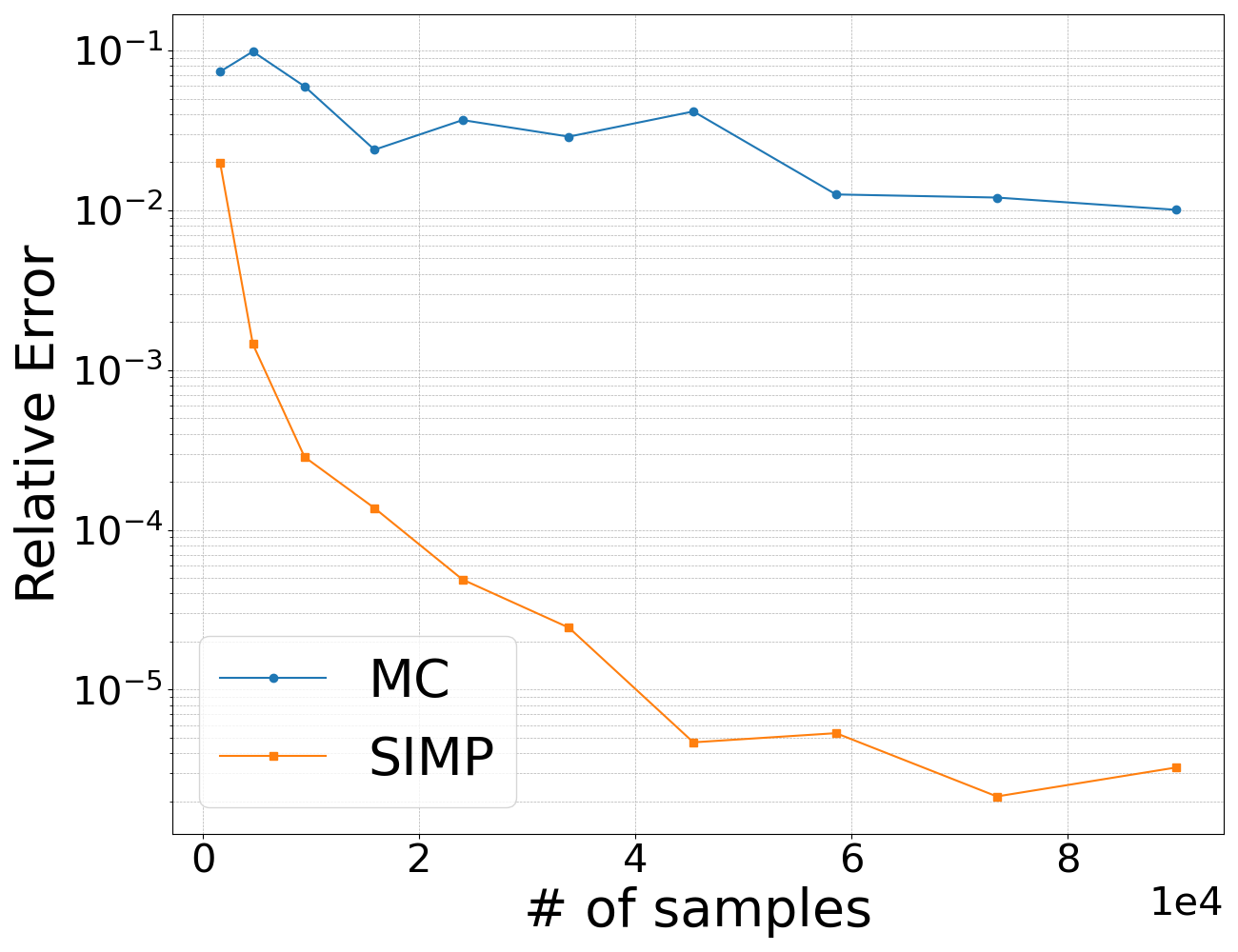}
    \caption{Relative errors of $u_{\text{WTN}}$ using Monte Carlo compared to the composite Simpson’s rule for numerical quadratures. Note that $\gamma = 1$ is used in the neural basis functions, with $\bm{\sigma} = 0.03$ and $N_l = 10$ for all test functions. }
    \label{fig:int_compare}
\end{figure}

The number of trial and test basis functions are set to 200. 
For MC, we generate $N_{\text{MC}}$ samples from the normal distribution with probability density $\psi_i$. For Simpson's rule, we generate $N_{S}$ uniformly spaced grid points. In both cases, boundary integrals are computed via the Simpson's rule. 
We then compare the relative errors of the approximate solutions $u_{\WTN}$ obtained using these two quadrature methods, evaluated over the same testing set of $129\times 129$ uniformly spaced grid points. The results are shown in Fig.~\ref{fig:int_compare}. It is seen that, as expected, Simpson's rule leads to more accurate approximation in this case than MC. The error decreases as $N_{\text{SIMP}}$ increases.  
It is worth noting that low integration accuracy often emerges as a common bottleneck of neural network-based methods when solving PDEs. 

Although finding the weak solution can handle problems with low regularity, many practical problems involve solutions that exhibit complex behavior. To address these challenges, we next propose several extensions of the WTN method. 

\section{Extensions of Weak TransNet method}\label{sec:wtn-extensions}
We improve the WTN method in two aspects: (i) incorporating Fourier features into WTN to solve multiscale problems; (ii) synthesizing partition of unity (PoU) functions with WTN to handle PDEs with sharply varying solutions.  

\subsection{A Fourier-WTN method}
\label{sec:f-wtn}
To enhance the capability of neural networks for approximating functions of multiple scales,  Fourier features are incorporated into the neural network models~\cite{tancik2020fourier,wang2021eigenvector,li2023deep}. 
To solve multiscale problems, we can combine WTN with the Fourier features. 

In general, the Fourier feature mapping is defined as $\mathfrak{F}(\cdot):\mathbb{R}^{d}\to \mathbb{R}^{2P}$
\begin{equation*}
\label{eq:fourier-map}
\mathfrak{F}(\bx) =\begin{bmatrix}\cos(2\pi \mathfrak{B}\bx)\\ \sin(2\pi \mathfrak{B}\bx)
\end{bmatrix},
\end{equation*}
where $\mathfrak{B}\in \mathbb{R}^{P\times d}$ is a random Gaussian matrix, meaning that each element of $\mathfrak{B}$ is sampled from a Gaussian distribution $\mathcal{N}(0,\sigma_B^2)$.
To incorporate it into the TransNet, the input $\bx$ is mapped to a high dimensional feature space before passing them through the network.
The resulting neural basis functions are modified from $\phi_j = \sigma( \gamma_j (\bm{a}_j^\top \bx + b_j)$ to the Fourier neural basis functions $\widetilde{\phi}_j = \sigma( \gamma_j (\bm{a}_j^\top \mathfrak{F}(\bx) + b_j))$. 

Defining the mapped space as $\Omega_{\F}\coloneqq\{ \mathfrak{F}(\bx) \in \mathbb{R}^{2P}|\bx\in \Omega\}$.
For any $\bx\in \Omega$, we have $\|\mathfrak{F}(\bx)\|_2 = \sqrt{P}$ regardless of the choice of $\Omega$, 
which implies that $\Omega_\F$ lies on the sphere of radius $\sqrt{P}$ in the $2P$-dimensional space.
Considering a ball centered at the origin with a slightly larger radius $R = \sqrt{P} + \epsilon_{\F}$, where $\epsilon_{\F} >0$, i.e., $B_{R}(\bm{0}) = \{ \bx \in \mathbb{R}^{2P} | \|\bx\|_2 \le R \}$. Then we have $\Omega_\F \subset B_{R}(\bm{0})$.
Following \cite[Theorem 2]{lu2025multiple}, if $\{a_j\}_{j=1}^M$ are i.i.d. and uniformly distributed on the unit sphere in $\mathbb{R}^{2P}$, and $\{r_j\}_{j=1}^M$ are i.i.d. and uniformly distributed in $[0, R]$. Then, the resulting hyperplanes are uniformly distributed within $B_{R}(\bm{0})$.

Consequently, we use 
\[
u_{\text{F-WTN}}(\bx) : = \sum_{j=0}^M \alpha_j \widetilde{\phi}_j(\bx)\,
\]
to approximate the weak solution in the framework of WTN, and name the resulting approach the Fourier-weak TransNet (F-WTN) method. 

\subsection{A PoU-WTN method}
\label{sec:pou_weaktransnet}
For problems exhibiting sharp gradients or singularities, the \textit{partition of unity} (PoU) method provides a flexible framework for blending locally adapted basis functions to accurately capture the solution’s behavior in critical regions. 
It can be seamlessly integrated with WTN to augment its ability to solve such challenging problems. 

To this end, we first consider a non-overlapping domain decomposition of $\Omega$ such that $\overline{\Omega}\coloneqq \cup_{\ell=1}^{L} \overline{\Omega}^{(\ell)}$, where each subdomain $\Omega^{(\ell)}$ is an open set whose boundary is denoted by $\partial\Omega^{(\ell)}$ and closure denoted by $\overline{\Omega}^{(\ell)}$. For any $\bx \in \partial\Omega^{(\ell)}$, we define the set of neighboring subdomains as  
\[
\Lambda^{(\ell)} (\bx) \coloneqq \{q \in \{1,2,\ldots, L\}|q\ne \ell, \bx \in \overline{\Omega}^{(q)}\}\,,
\] 
and the number of neighboring subdomains at $\bx$ is the cardinality $|\Lambda^{(\ell)} (\bx)|$.
We then define the PoU function $ \chi^{(\ell)}: \mathbb{R}^d \to \mathbb{R}$ as 
\begin{equation}
\chi^{(\ell)}(\bx) = \begin{cases}
    1 \,,\quad \bx \in \Omega^{(\ell)}\,,\\ 
    \frac{1}{|\Lambda^{(\ell)}(\bx)| +1} \,,\quad \bx \in \partial \Omega^{(\ell)}\,,\\ 
    0\,,\quad \text{otherwise}\,, 
\end{cases}
\text{ for } \ell=1, \ldots, L.
\label{eq:chi}
\end{equation}
Obviously, the property $\sum_{\ell=1}^{L} \chi^{(\ell)}(\bx)=1$ holds for any $ \bx \in \Omega$. 
Note that different types of PoU functions can also be considered~\cite{melenk1996partition}.

In the $\ell$-th subdomain, we generate $M^{(\ell)}$ neural basis functions and approximate a local solution in this subdomain by
\[u_{\TN}^{(\ell)}(\bx) = \sum_{j=0}^{M^{(\ell)}}\alpha_{j\ell}\phi_j^{(\ell)}(\bx)\,,\quad \bx \in \overline{\Omega}^{(\ell)}\,.
\]
Note that the neural basis functions are globally supported, to better represent local features, we use the PoU functions and, consequently, define the solution over the entire domain as  
\begin{equation}
u_{\pou-\TN}(\bx) = \sum_{\ell=1}^L \chi^{(\ell)}(\bx)u_{\TN}^{(\ell)}(\bx)\,,\quad \bx \in \overline{\Omega}\,.
\label{eq:pou_wtn}
\end{equation}
Introducing a new notation $\zeta_\nu(\bx)$ to denote a localized basis $\chi^{(\ell)}(\bx)\phi_j^{(\ell)}(\bx)$, where the index $\nu$ iterates over all possible pairs $(\ell, j)$,  
we represent the entire set of localized basis functions as $\{\zeta_\nu(\bx)\}_{\nu=1}^{\overline{M}}$. Here $\overline{M}=\sum_{\ell=1}^L (M^{(\ell)} + 1)$ is the total number of basis functions.   
Then, \eqref{eq:pou_wtn} can be rewritten as   
\[
u_{\pou-\TN}(\bx) = \sum_{\nu=1}^{\overline{M}} \alpha_\nu \zeta_\nu(\bx)
\,.
\]
To determine it, we can use the framework of WTN that minimizes the loss function \eqref{eq:weak_loss}. Thanks to the use of localized basis, $u_{\pou-\TN}$ is able to better capture local features of the PDE solution than $u_{\TN}$. However, certain conditions across interfaces of subdomains must be imposed for it to reach the same regularity as the weak solution. 
To proceed, we consider two arbitrary adjacent subdomains $\Omega^{(\ell)}$ and $\Omega^{(q)}$, denote their shared interface by $\Gamma^{(\ell,q)} \coloneqq \partial \Omega^{(\ell)}\cap \partial \Omega^{(q)}$, 
and elaborate the procedure. 

The first fundamental interface condition is the continuity of the solution across the interface, which requires 
$\jump{u_{\pou-\TN}(\bx)}^{(\ell,q)} = 0$ for any $\bx \in \Gamma^{(\ell,q)}$, where the notation $\jump{\cdot}^{(\ell,q)}$ denotes the jump of a quantity across the interface $\Gamma^{(\ell,q)}$, for example,  
\begin{align*}
\jump{u_{\pou-\TN}(\bx)}^{(\ell,q)}
& = 
\lim_{\bx\in \Omega^{(\ell)}, \bx \to \Gamma^{(\ell,q)}} u_{\pou-\TN}(\bx) -
\lim_{\bx\in \Omega^{(q)}, \bx \to \Gamma^{(\ell,q)}} u_{\pou-\TN}(\bx)  \\
& = 
u^{(\ell)}_{\TN}(\bx) - u^{(q)}_{\TN}(\bx)\,.
\end{align*} 
To impose this constraint, we introduce a set of \textit{interface samples} $S_{\Gamma^{(\ell,q)}} =\{\bx^{(m)}\}_{m=1}^{N_{\Gamma^{(\ell,q)}}}$ on $\Gamma^{(\ell,q)}$ and require
\begin{equation}
u^{(\ell)}_{\TN}(\bx^{(m)}) = u^{(q)}_{\TN}(\bx^{(m)})\,, \quad\forall\, \bx^{(m)}\in S_{\Gamma^{(\ell,q)}}\,.
\label{eq:interface}
\end{equation}
Define $\Phi\in \mathbb{R}^{N_{\Gamma^{(\ell,q)}}\times \overline{M}}$ with the entry $\Phi_{mk} = \zeta_k(\bx^{(m)})$,
and a mask matrix of the same size  
\begin{equation}
\mathfrak{M}^{(\ell,q)}\coloneqq [\eta^{(1)}\bm{E}^{(1)}|\eta^{(2)}\bm{E}^{(2)}|\cdots|\eta^{(L)}\bm{E}^{(L)}]\,
\label{eq:mask_matrix}
\end{equation}
with $\bm{E}^{(p)}\in \mathbb{R}^{N_{\Gamma^{(\ell,q)}}\times (M^{(p)}+1)}$ an all-one matrix ($p=1,\ldots,L)$, $\eta^{(\ell)}=1$,  $\eta^{(q)} =-1$ and $\eta^{(p)}=0$ for $p=1,\ldots, L$ and $p\ne \ell,q$. Then~\eqref{eq:interface} can be equivalently written as 
$
[\mathfrak{M}^{(\ell,q)}\circ \Phi]\balpha= \bm{0}\,,
$
where $\circ$ denotes the Hadamard product. 
Denoting the block as $\mathcal{M}^{0,(\ell,q)} \coloneqq [\mathfrak{M}^{(\ell,q)}\circ \Phi]$, where the superscript $\cdot^{0}$ indicates the zeroth-order derivative continuity across the interface.
Then, by iterating over all the interfaces, we concatenate the row blocks $\mathcal{M}^{0,(\ell,q)}$ into one global interface matrix $\mathcal{M}^{0}$. 

For elliptic PDEs, such as the Poisson equation and Darcy flow, enforcing the continuity of the solution across subdomain interfaces is not sufficient.
To ensure global conservation and maintain physical consistency, we impose the additional flux continuity condition:
\begin{equation}
\jump{\kappa(\bx)\nabla u_{\pou-\TN}\cdot \vec{\bm{n}}^{(\ell,q)}}^{(\ell,q)} = 0\,,\quad \forall \bx \in \Gamma^{(\ell,q)}\,,
\label{eq:flux_condition}
\end{equation} 
where $\vec{\bm{n}}^{(\ell,q)} = (n_1, n_2,\ldots,n_d)$ denotes the unit normal to the interface $\Gamma^{(\ell,q)}$, 
$\kappa(\bx)$ is a scalar function representing thermal conductivity or diffusion coefficient in the diffusion equation and Poisson, and representing the permeability field in Darcy flow.  

When $\kappa(\bx)\equiv \kappa$ is a constant, the flux continuity condition~\eqref{eq:flux_condition} reduces to the requirement 
$\nabla u_{\TN}^{(\ell)}(\bx) \cdot \vec{\bm{n}}^{(\ell,q)} - \nabla u_{\TN}^{(q)}(\bx) \cdot \vec{\bm{n}}^{(\ell,q)} =0$, for all $\bx \in \Gamma^{(\ell,q)}$. 
To enforce this condition,
 we define $\Phi_{x_k}\in \mathbb{R}^{N_{\Gamma^{(\ell,q)}}\times \overline{M}}$ $(k=1,\ldots,d)$, where each entry corresponds to the spatial derivative of basis functions at the interface nodes $[\Phi_{x_k}]_{m\nu} = \frac{\partial \zeta_\nu}{\partial x_k}(\bx^{(m)})$. 
Using the same mask matrix defined in~\eqref{eq:mask_matrix}, the flux continuity can be enforced via $[\mathfrak{M}^{(\ell,q)} \circ (\Phi_{x_1} n_1 + \cdots + \Phi_{x_d}n_d)] \balpha = \bm{0}$.
We denote the matrix $\big[\mathfrak{M}^{(\ell,q)} \circ (\Phi_{x_1} n_1 + \cdots + \Phi_{x_d} n_d)\big]$ as $\mathcal{M}^{1,(\ell,q)}$, where the superscript $\cdot^{1}$ indicates the first order derivative continuity across the interface. 

When $\kappa(\bx)$ is discontinuous across the interface, the definition of the mask matrix needs to be modified. For simplicity, we illustrate it using the case where $\kappa(\bx)$ is a piecewise constant. Specifically, let $\lim_{\bx\in \Omega^{(\ell)}, \bx \to \Gamma^{(\ell,q)} }\kappa(\bx) = \kappa^{(\ell)}$ and  $\lim_{\bx\in \Omega^{(q)}, \bx \to \Gamma^{(\ell,q)} } \kappa(\bx)  = \kappa^{(q)}$, where $\kappa^{(\ell)} $ and $\kappa^{(q)}$ are constants, 
then the corresponding mask matrix is defined as 
\[
\mathfrak{M}^{(\ell,q)} = [\bm{F}^{(1)}|\cdots |\bm{F}^{(L)}]\,.
\]
Each $\bm{F}^{(p)}\in \mathbb{R}^{N_{\Gamma^{(\ell,q)}}\times (M^{(p)}+1)}$ is an all-zero matrix for $p=1,\ldots,L$ and $p\ne \ell,q$, while 
$\bm{F}^{(\ell)} = \kappa^{(\ell)} \bm{E}^{(\ell)}$ and $\bm{F}^{(q)} = - \kappa^{(q)}\bm{E}^{(q)}$ with $\bm{E}^{(\ell)}$ and $\bm{E}^{(q)}$ defined in \eqref{eq:mask_matrix}. 
Consequently, $\mathcal{M}^{1,(\ell,q)}$ is constructed as described in the preceding discussion.

For each interface $\Gamma^{(\ell,q)}$, the matrices $\mathcal{M}^{r,(\ell,q)}$ (with $r = 0, 1$) are vertically stacked and form the global interface matrix $\mathcal{M}$. 
To enforce these interface conditions, an extra term weighted by $\lambda$ is introduced to the loss function \eqref{eq:weak_loss}:
\begin{equation}
\lambda\|\mathcal{M} \balpha  \|_2^2\,.
\label{eq:interface1}
\end{equation} 
It simultaneously penalizes discontinuities in the solution and its derivatives. Obviously, the same procedure can be extended to accommodate more general or alternative types of interface conditions.

This method is referred to as the partity of unity-weak TransNet (PoU-WTN) method. As in Section~\ref{sec:weak_form}, we present a detailed description of it for solving the linear PDE in Alg.~\ref{alg:pou-wtn}. 

\begin{algorithm}[!htp]
\caption{A PoU-Weak TransNet ($\pou-$WTN) method for PDEs}\label{alg:pou-wtn}
\begin{algorithmic}[1]
\Require{PoU strategy ($\{\Omega^{(\ell)}\}_{\ell=1}^{L}$),
number of local trial basis $\{M^{(\ell)}\}_{\ell=1}^{L}$, number of test basis $N$, PDE form \eqref{eq:problem}, shape parameters $\{\gamma_j^{(\ell)}\}_{j=1,\ldots,M^{(\ell)}, \ell =1,\ldots,L}$, boundary samples $S_{\partial \Omega}\coloneqq\{\bx_{\partial \Omega}^{(m)}\}_{m=1}^{N_{\partial \Omega}}$, 
interface samples $S_{\Gamma^{(\ell,q)}} = \{\bx^{(m)}\}_{m=1}^{N_{\Gamma^{(\ell,q)}}}$ for all interfaces, interface weight $\lambda$, adjusted boundary weight $\widetilde{\beta}$}
\Ensure{$u_{\text{PoU}-\WTN}(\bx; \balpha^{\star})$}
\State{Construct the trial basis $\{\phi_j^{(\ell)}(\bx)\}_{j=0,\ldots,M^{(\ell)}, \ell=1,\ldots,L}$}
\State{Construct the test basis $\{\psi_i(\bx; \bm{\mu}, \bm{\sigma})\}_{i=1}^N$}
\For{$\ell = 1,\ldots,L$}
\For{$i=1,\ldots,N$}
\State{Identify the local support $\Omega_{\psi_i}$of $\psi_i$} 
\State{Compute the the right-hand-side vector entry locally $\bm{f}_i^{(\ell)} = l(\psi_i)_{\Omega^{(\ell)} \cap \Omega_{\psi_i}} $
}
\For{$j = 0,\ldots, M^{(\ell)}$}
\State{Compute the local stiffness matrix entry locally $\bm{A}_{ij}^{(\ell)} = a(\phi_j^{(\ell)}, \psi_i)_{\Omega^{(\ell)}\cap \Omega_{\psi_i}}$} 
\EndFor
\EndFor
\EndFor
\State{Assemble the global stiffness matrix 
$\bm{A} = [\bm{A}^{(1)}|\cdots |\bm{A}^{(L)}]$ and the global RHS vector $\bm{f} = \sum_{\ell=1}^L \bm{f}^{(\ell)}$}
\State{Compute the boundary matrix $\bm{B}\in \mathbb{R}^{N_{\partial \Omega}\times (M+1)}$ with entry $\bm{B}_{mj} = [\mathcal{B}[\phi_j(\bx^{(m)})]]$ and the boundary condition vector $\bm{g}\in \mathbb{R}^{N_{\partial \Omega}}$ with entry $\bm{g}_m = g(\bx^{(m)})$}
\For{any interface $(\ell,q)$}
\State{Compute the corresponding interface matrix $\mathcal{M}^{r, (\ell,q)}$ for $r=0,1$}
\EndFor
\State{Grouping all interface matrices as $\mathcal{M}$}
\State{Solve the LS system:
\[
\balpha^{\star} = \argmin_{\balpha} \left\| \begin{bmatrix}
    \bm{A} \\
    \widetilde{\beta} \bm{B} \\
     \gamma\mathcal{M}
\end{bmatrix} \balpha - \begin{bmatrix}
    \bm{f} \\
    \widetilde{\beta}  \bm{g}\\
    \bm{0}
\end{bmatrix} \right\|_2^2
\]
}
\State{Obtain the approximate solution $u_{\text{PoU}-\WTN}(\bx;\balpha^{\star})$.}
\end{algorithmic}
\end{algorithm} 

\begin{remark}
Due to the specific choice of PoU functions in this work, the PoU-WTN method naturally aligns with the non-overlapping domain decomposition method. However, one may instead choose smooth PoU functions to avoid implementation of interface conditions, which we plan to explore in future work.  
\end{remark}

\section{Numerical examples}
\label{sec:num} 
Several numerical experiments are conducted in this section. The first experiment in Sec.~\ref{sec:empirical_shape} investigates the empirical behavior of the shape parameter $\gamma$, while the remaining examples are used to validate the WTN method (presented in Alg.~\ref{alg:wtn}) and its extensions, including F-WTN and PoU-WTN (presented in Alg.~\ref{alg:pou-wtn}). 
More specifically, 
Sec.~\ref{sec:onlyweak_v2} studies a Darcy flow problem lacking a strong solution,  
Sec.~\ref{sec:high-fre} considers a Darcy flow with multiscale features, and 
Sec.~\ref{sec:channel} investigates a Darcy flow problem featuring  channelized permeability.
Sec.~\ref{sec:vpinn} considers the Poisson equation with a sharp-gradient solution,  
while Sec.~\ref{sec:l_singularity} discusses an L-shape domain problem with a solution singularity. 

We compare the proposed methods with TransNet using different loss functions, including the strong form (referred to as SF, presented in Sec.~\ref{sec:strong_loss}) and Ritz energy (referred to as DRM, presented in Sec.~\ref{sec:drm_loss}). 
Specifically, we use the Monte Carlo method to compute the integrals in SF and DRM. 
On a test set, distinct from the training data, 
we evaluate the accuracy of the obtained solution $u_{\star}$ ($\star$ = method)  using the following relative $L_2$ error:
\begin{equation}
e_{\star} = \frac{\|u^{\exact}-u_{\star}\|_2}{\|u^{\exact}\|_2}
\label{eq:error}
\end{equation}
where $u^{\exact}$ is the exact solution or a benchmark solution obtained from a finite element simulation.
In all numerical examples, for accurate numerical integration in the weak formulation, we employ the compSimpson's rule \footnote{Implemented using the \textsf{scipy.integrate.simpson} function from the \textsf{SciPy} library.} to compute the stiffness matrix.
Also, we employ the \textsf{numpy.linalg.lstsq} function to solve all the least-squares problems. 
Unless otherwise specified, the centers $\{\mu_i\}_{i=1}^N$ of the test functions are randomly sampled in $\Omega$ according to a uniform distribution, and $N_l$ is set to 10 by default.

All methods are implemented in \textsf{python} and the numerical results are conducted on Intel Xeon Gold 6252 CPUs. 

\subsection{An empirical study on the shape parameter}
\label{sec:empirical_shape}
Intuitively, a larger value of the shape parameter $\gamma$ is required for the neural basis to accurately approximate a target function with a steeper gradient. However, an optimal choice of $\gamma$ remains an open question for numerical solutions of PDEs. Therefore, we first perform an empirical study on the choice of $\gamma$. 
In this example, all neurons are assumed to have the same shape parameter, i.e., $\gamma_1=\cdots = \gamma_M = \gamma$.

To this end, we take the following two-dimensional Gaussian function as the target function 
\[
T(x,y) = \frac{1}{2\pi \sigma_f^2} \exp\left( \frac{x^2-y^2}{2\sigma_f^2} \right)\,,
\]
and check how well it can be approximated in the neural feature space spanned by basis functions of different shape parameter values. 
To measure it, we compute the projection error evaluated on a set of testing samples 
$\{(x^{(i)},y^{(i)})\}_{i=1}^{N_{\text{test}}}$, as follows: 
\[
e_{\text{proj}}\coloneqq \frac{\| \bm{A}\balpha^{\star}   - \bm{t}\|_2}{\|\bm{t}\|_2}
\]
where $\bm{t}= [T(x^{(1)}, y^{(1)}), T(x^{(2)}, y^{(2)}), \ldots, T(x^{(N_{\text{test}})},y^{(N_{\text{test}} )})]^{\top}$, $\bm{A}$ is a $N_{\text{test}}\times (M+1)$ matrix whose entry $\bm{A}_{ij}= \phi_j(x^{{(i)}}, y^{(i)})$, 
and 
$\balpha^{\star} = \argmin \| \bm{A}\balpha - \bm{T} \|_2^2.$

\begin{figure}[!htb]
\centerline{
    \includegraphics[width=1.2\linewidth]{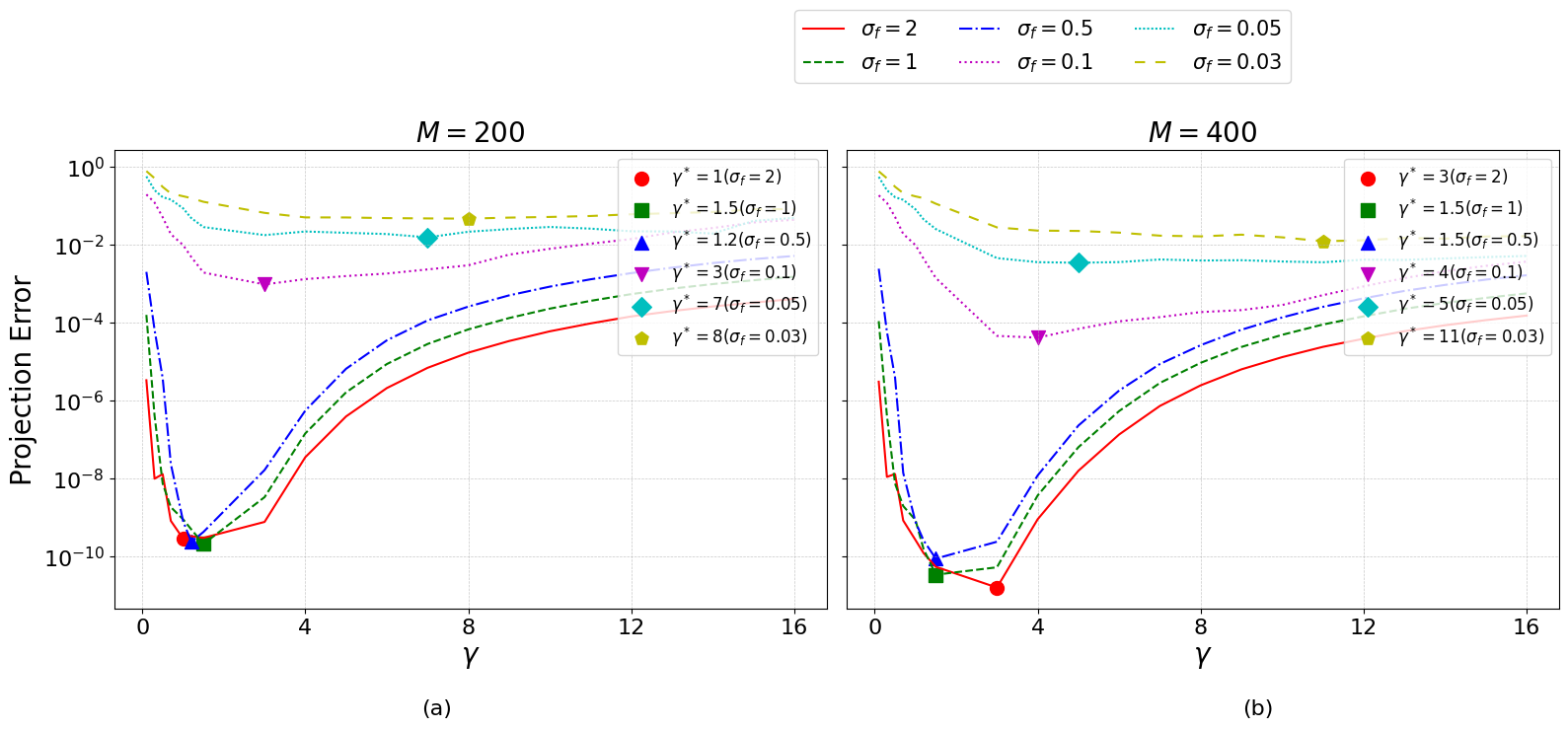}
}
    \caption{Projection error versus the shape parameter $\gamma$ for different values of $\sigma_f$: (a) $M=200$ and (b) $M=400$. 
    Colored curves correspond to different $\sigma_f$ values, with matching markers indicating the optimal $\gamma$ associated to the least projection error for each case.}
    \label{fig:gamma}
\end{figure}

While keeping the number of trial basis functions fixed, we vary $\sigma_f\in [0.03, 0.05, 0.1, 0.5, 1, 2]$ and the shape parameter $\gamma\in [0.1,16]$, then compute the associated projection errors. 
Note that $\sigma_f$ directly determines the sharpness of $f$.
For $M = 200$ and $M=400$, Tab.~\ref{tab:re_gamma} lists the corresponding optimal values of $\gamma$; and Fig.~\ref{fig:gamma} displays the projection error as a function of $\gamma$ for selected values of $\sigma_f$.  
It is observed that for slowly varying target functions ($\sigma_f\ge 0.5$), optimal projection accuracy is achieved at $\gamma\in [1,3]$. For moderately sharp functions ($0.05 \le \sigma_f \le 0.5$), the optimal $\gamma$ ranges from $[3, 7]$. While for sharp functions ($\sigma_f=0.03$), the optimal $\gamma$ increases further to $[8, 11]$. 

Based on these findings, in the subsequent experiments, we set $\gamma_1 = \cdots = \gamma_M =1$ for problems without sharp gradients (see Sec.~\ref{sec:onlyweak_v2}, Sec.~\ref{sec:high-fre} and Sec.~\ref{sec:channel}), $\gamma_1 = \cdots = \gamma_M= 5$ for problems with sharp gradients (see Sec.~\ref{sec:vpinn}), and discuss a special treatment for singular problems where the shape parameters are not identical for all neurons (see Sec.~\ref{sec:l_singularity}).

\begin{table}[!htp]
	\caption{Optimal shape parameters for different pairs of $(\sigma_f, M)$}
    \vspace{0.1cm}
	\centering
    \adjustbox{max width=\textwidth}{
    \begin{tabular}{ccccccc}
\hline
 & $\sigma_f = 2$ & $\sigma_f = 1$ & $\sigma_f = 0.5$ & $\sigma_f = 0.1$ & $\sigma_f = 0.05$ & $\sigma_f = 0.03$ \\    
\hline
 $M=200$& 1 & 1.5 & 1.2 & 3 & 7& 8 \\ 
 $M=400$& 3 & 1.5 & 1.5 & 4 & 5 &  11\\
\hline
\end{tabular}
}
	\label{tab:re_gamma}
\end{table}


\subsection{Darcy flow without a strong solution}
\label{sec:onlyweak_v2}

\begin{figure}
\centerline{
    \begin{tabular}{cc} 
    \includegraphics[width=0.5\textwidth]{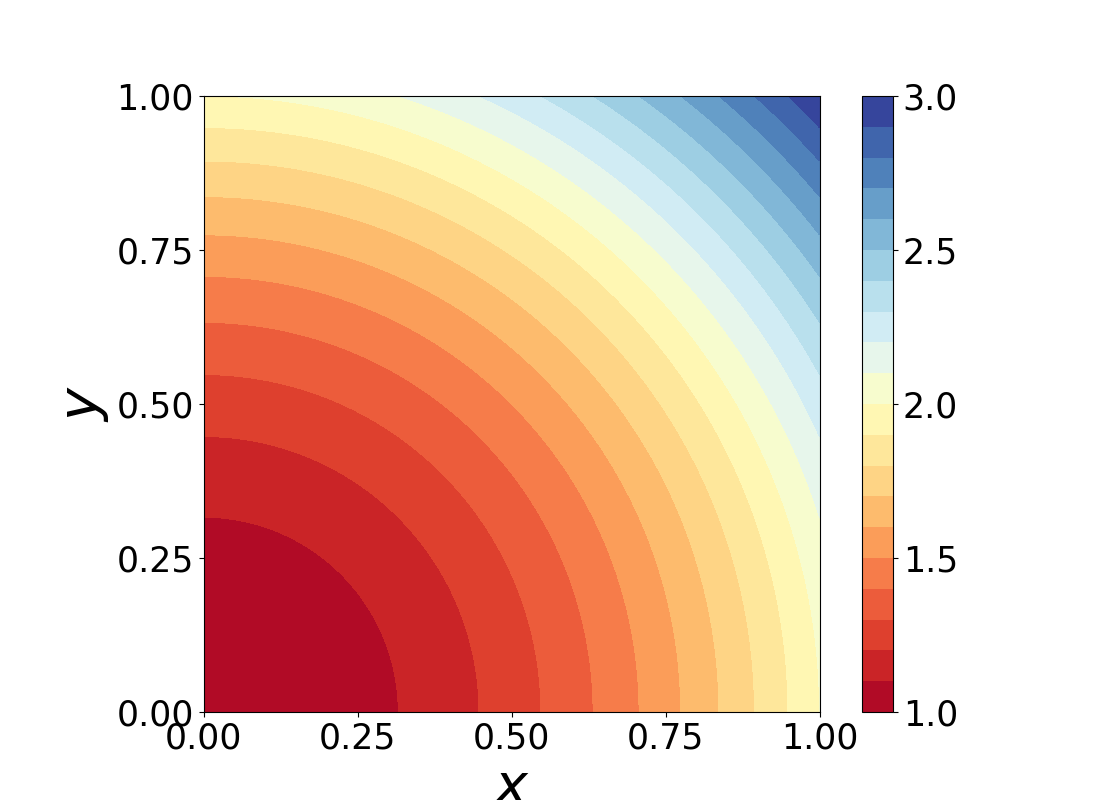} & 
    \includegraphics[width=0.5\textwidth]{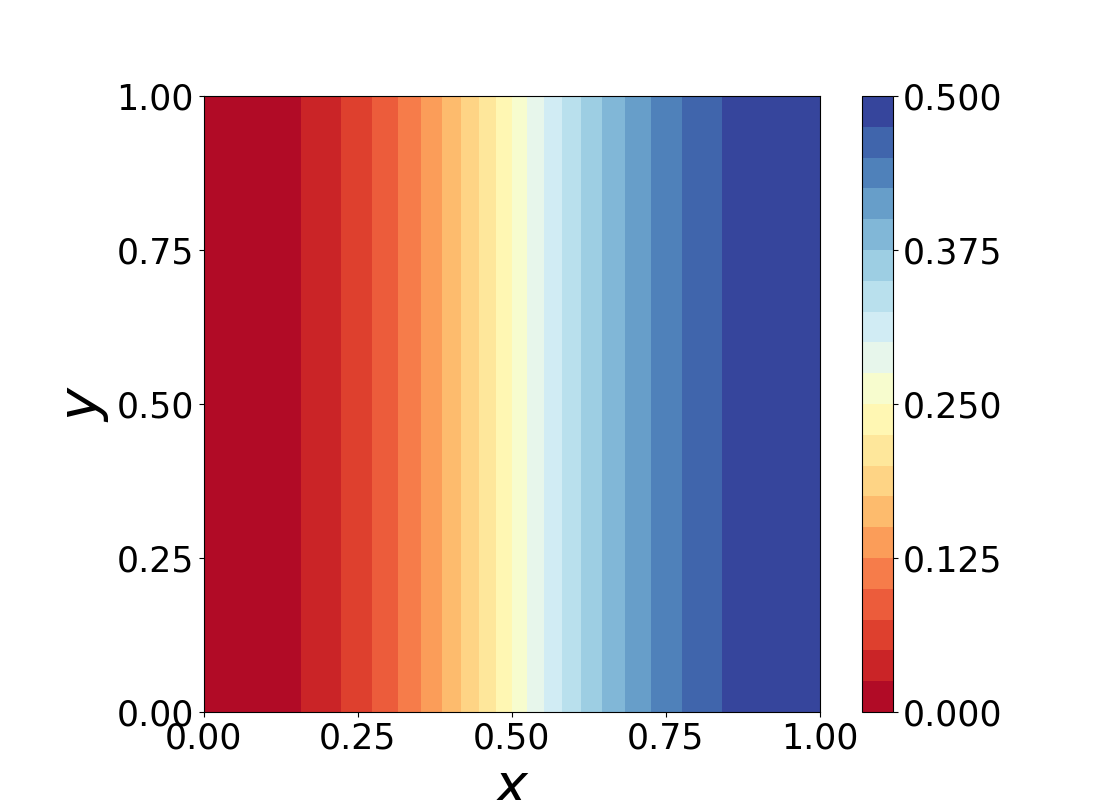} 
    \\
    (a) $\kappa(\bx)$ & (b)  $u^{\exact}$\\
    \end{tabular} \\
    }
\caption{Darcy flow in Sec.~\ref{sec:onlyweak_v2}: (a) permeability $\kappa(\bx)$; and (b) reference solution.}
\label{fig:setup_onlyweak_v2}
\end{figure}

Next, we consider the two dimensional Darcy flow with Dirichlet boundary condition, that is $\mathcal{L}[u]  = -\nabla \cdot (\kappa(\bx) \nabla u)$ and $\mathcal{B}[u]=u$ in~\eqref{eq:problem} with 
$\bx\in\Omega = [0,1]^2$, $\kappa(\bx) = 1 +|\bx|^2 = 1 + x^2 + y^2$ as shown in Fig.~\ref{fig:setup_onlyweak_v2}(a), the source term   
\[
f = \begin{cases}
    -2 -6x^2 -2y^2\,,\quad 0\le x \le \frac{1}{2} \,,\\
    6x^2 + 2y^2 - 4x + 2\,,\quad \frac{1}{2}< x \le 1\,,
\end{cases}
\]
and the boundary condition $g(x,0) = g(x,1) = x^2$ for $0\le x \le \frac{1}{2}$, $g(x,0) = g(x,1) = -x^2 + 2x - 0.5$ for $\frac{1}{2}< x \le 1$ on $\partial \Omega$. Due to the discontinuity of the source term, 
this problem does not admit a strong solution, but only a unique weak solution $u^{\exact} = x^2$ for $0\le x \le \frac{1}{2}$, and $u^{\exact} =-x^2 + 2x - 0.5 $ for $\frac{1}{2}< x \le 1$, which is shown in Fig.~\ref{fig:setup_onlyweak_v2}(b).

 In this case, the matrix $\bm{A}$ in WTN has the entry as 
\[
\bm{A}_{ij} = a(\phi_j,\psi_i) = \int_{\Omega} \kappa \nabla \phi_j \cdot \nabla \psi_i \dif \bx - \int_{\partial \Omega} \psi_i\frac{\partial \phi_j}{\partial \vec{n}} \dif s\,. 
\]
When SF is used, $\kappa$ has to be differentiated, and the corresponding loss is   
\begin{align*}
\mathfrak{L}_{\stf} &= \int_{\Omega}\|- \kappa(\bx) \Delta u_{\TN}(\bx)  - \nabla \kappa(\bx) \cdot \nabla u_{\TN}(\bx) - f(\bx) \|_2^2 \dif \bx \\
&+ \beta_{\stf}\int_{\partial \Omega}\|\mathcal{B}[u_{\TN}(\bx)] - g(\bx) \|_2^2 \dif s \,,
\end{align*}
where $\beta_{\stf}$ is a weight that has been optimized.
If DRM method is instead used, the corresponding loss is 
\begin{align*}
\mathfrak{L}_{\DRM}
 &= \int_{\Omega}\left(\frac{1}{2}\kappa(\bx)|\nabla u_{\TN}(\bx)|^2
 - f(\bx)u_{\TN}(\bx)\right)\dif \bx \\
 &+ \beta_{\DRM}\int_{\partial \Omega}\|\mathcal{B}[u_{\TN}(\bx)] - g(\bx) \|_2^2 \dif s \,,
\end{align*}
where $\beta_{\DRM}$ is a weight that has been optimized.

To obtain $u_{\stf}$ and $u_{\DRM}$, we use randomly generated $10^3$ interior samples
and $800$ boundary samples (200 collocation samples on each side of $\partial \Omega$) according to a uniform distribution. 
The same set of boundary samples is also used for finding $u_{\WTN}.$ 
The weight $\beta$ for the boundary loss is set to as 1. 
The shape parameter of all neurons in the TransNet are set to 1 in this example. 
The standard deviation $\sigma_i$ in all test functions \eqref{eq:test_function} is set to $0.05$.

To measure the accuracy of these methods, a test set is used that consists of $129\times 129$ uniformly spaced samples generated over $\Omega$. 
A comparison of the errors across different numbers of trial basis and test basis functions is presented in Tab.~\ref{tab:rel_err_onlyweak_v2}.   

\begin{table}[]
\caption{Darcy flow in Sec. \ref{sec:onlyweak_v2}: relative errors of $u_{\WTN}$, $u_{\stf}$, $u_{\DRM}$ when different number of trial and test basis is used.}
\vspace{0.1cm}
\centering
\adjustbox{max width=\textwidth}{ 
\begin{tabular}{c|ccc|ccc}
\hline
$M$ & \multicolumn{3}{c|}{100}                                       & \multicolumn{3}{c}{200}                                        \\ \hline
$e_{\stf}$ & \multicolumn{3}{c|}{$8.28\times 10^{-2}$}& \multicolumn{3}{c}{$8.30\times 10^{-2}$}\\ \hline
$e_{\DRM}$& \multicolumn{3}{c|}{$2.14\times 10^{-2}$}& \multicolumn{3}{c}{$2.60\times 10^{-2}$}\\ \hline
\multirow{2}{*}{$e_{\WTN}$ } & \multicolumn{1}{c}{$N$=50} & \multicolumn{1}{c}{$N$=100} & $N$=200 & \multicolumn{1}{c}{$N$=100} & \multicolumn{1}{c}{$N$=200} & $N$=300 \\ \cline{2-7}  & \multicolumn{1}{c}{$9.59\times 10^{-2}$}     & \multicolumn{1}{c}{$4.35\times 10^{-3}$} &   $5.60\times 10^{-3}$    & \multicolumn{1}{c}{$1.62\times 10^{-2}$} & \multicolumn{1}{c}{$3.26\times 10^{-3}$}  & $1.81\times 10^{-3}$ \\\hline
\end{tabular}
}
\label{tab:rel_err_onlyweak_v2}
\end{table}

\begin{figure}[!htb]
\centerline{
    \begin{tabular}{ccc}

    \includegraphics[width=0.4\textwidth]{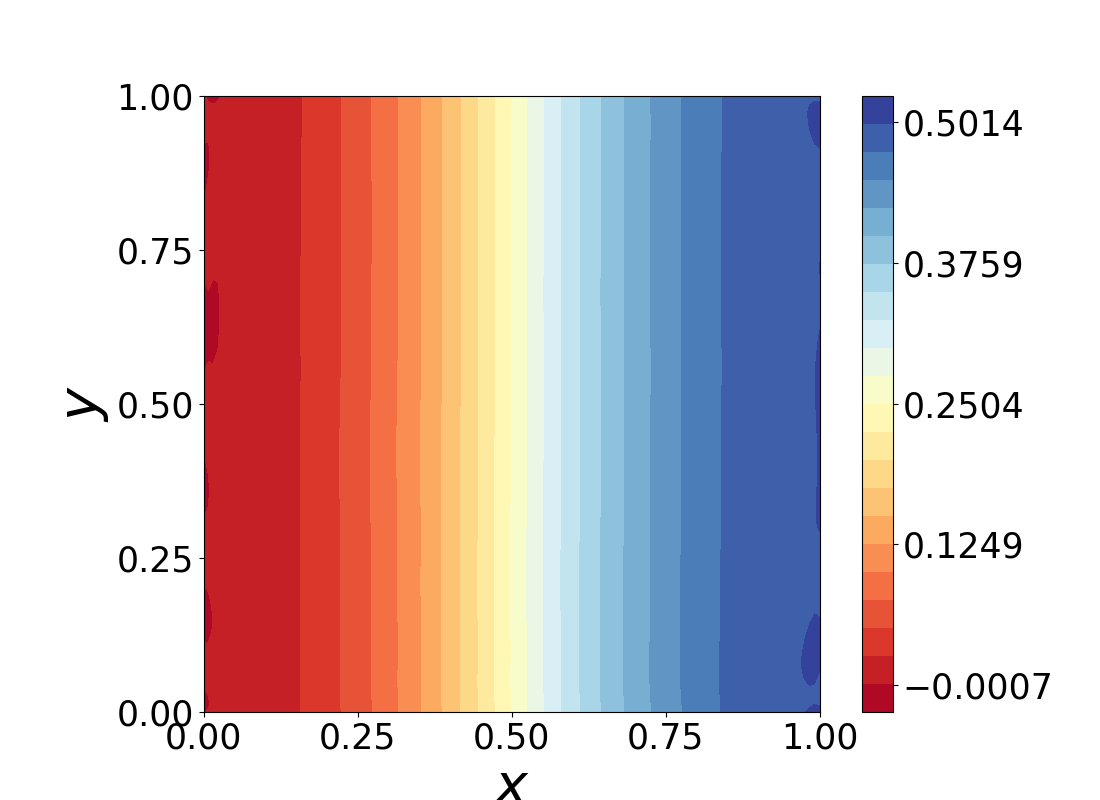} & 
    \includegraphics[width=0.4\textwidth]{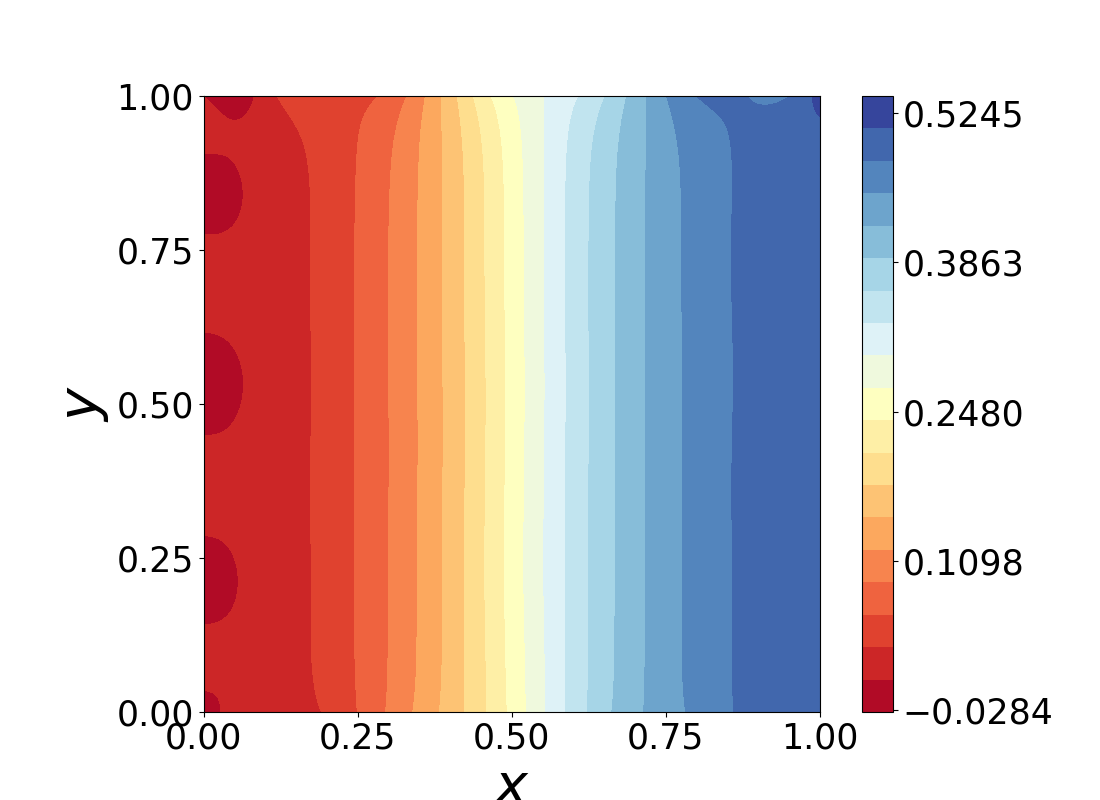}  & 
    \includegraphics[width=0.4\textwidth]{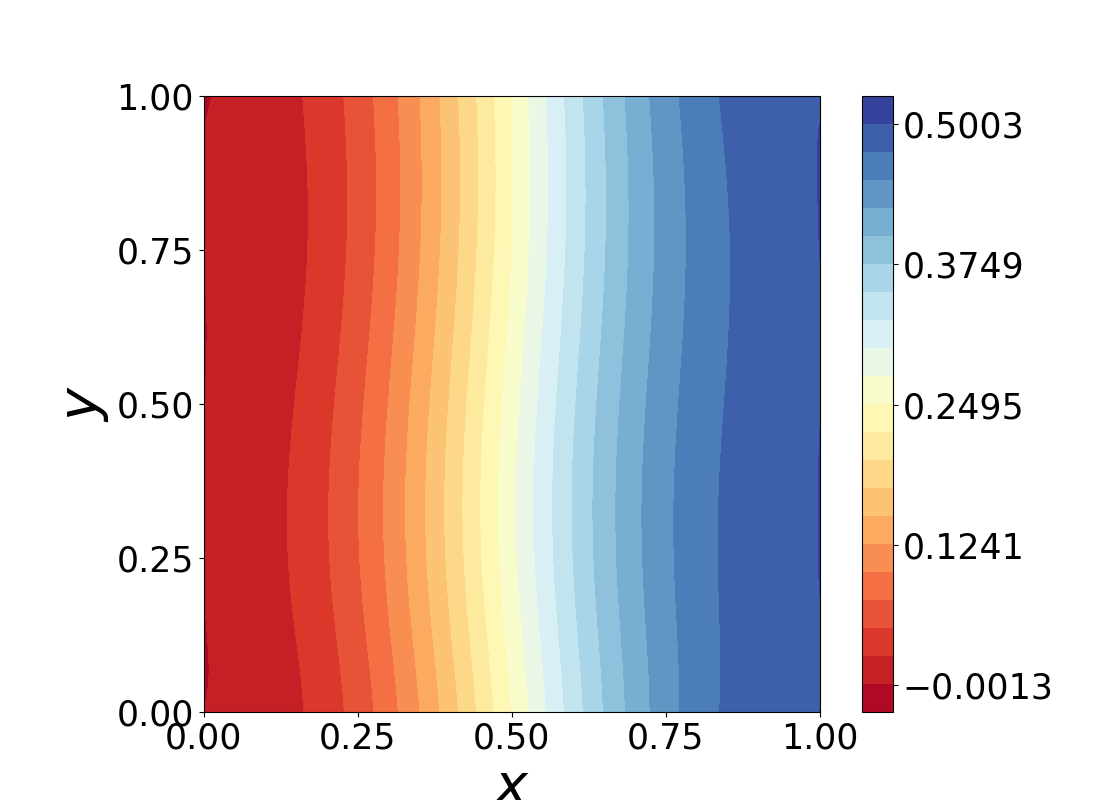}
    \\
    (a) $u_{\WTN}$ & (c)  $u_{\stf}$ &  (d) $u_{\DRM}$ \\
    \includegraphics[width=0.4\textwidth]{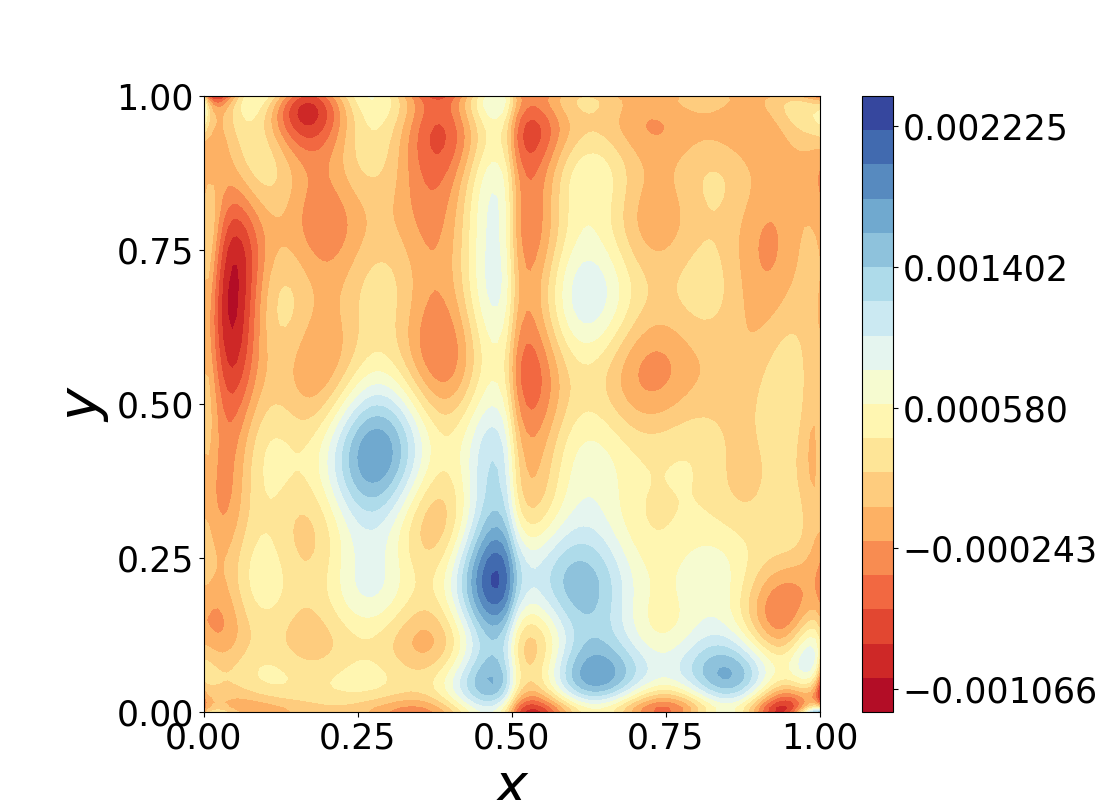} & 
    \includegraphics[width=0.4\textwidth]{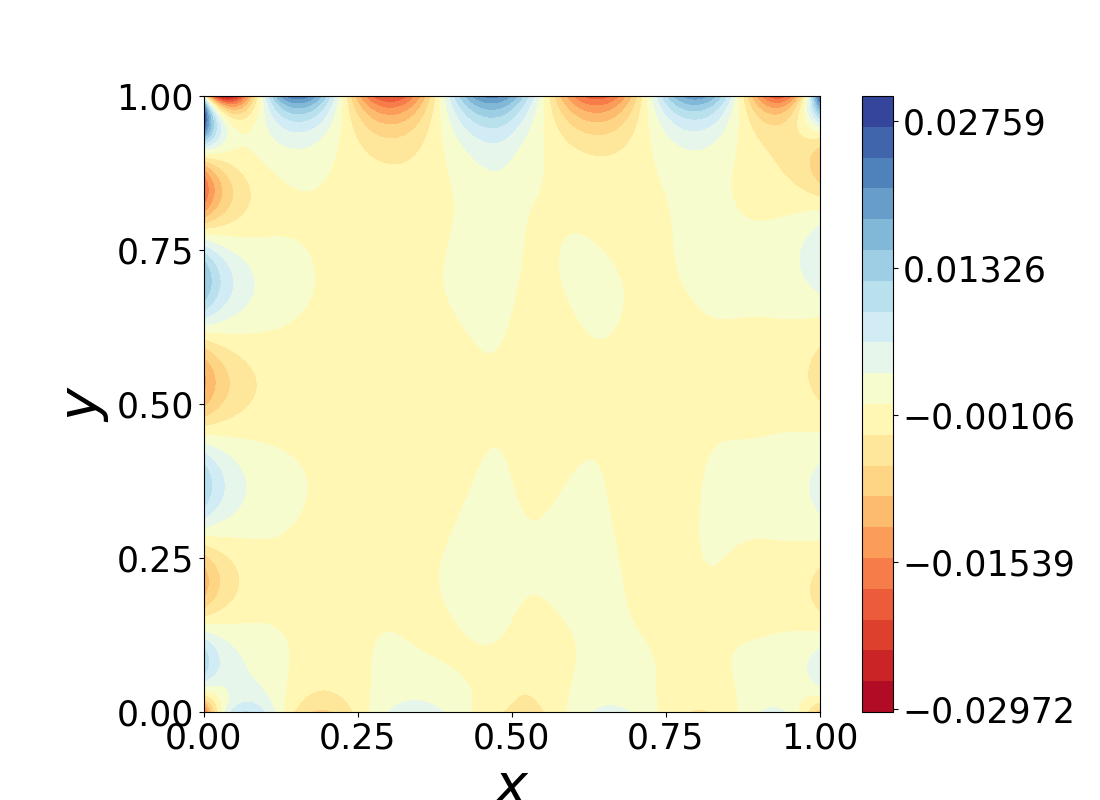}
    & 
    \includegraphics[width=0.4\textwidth]{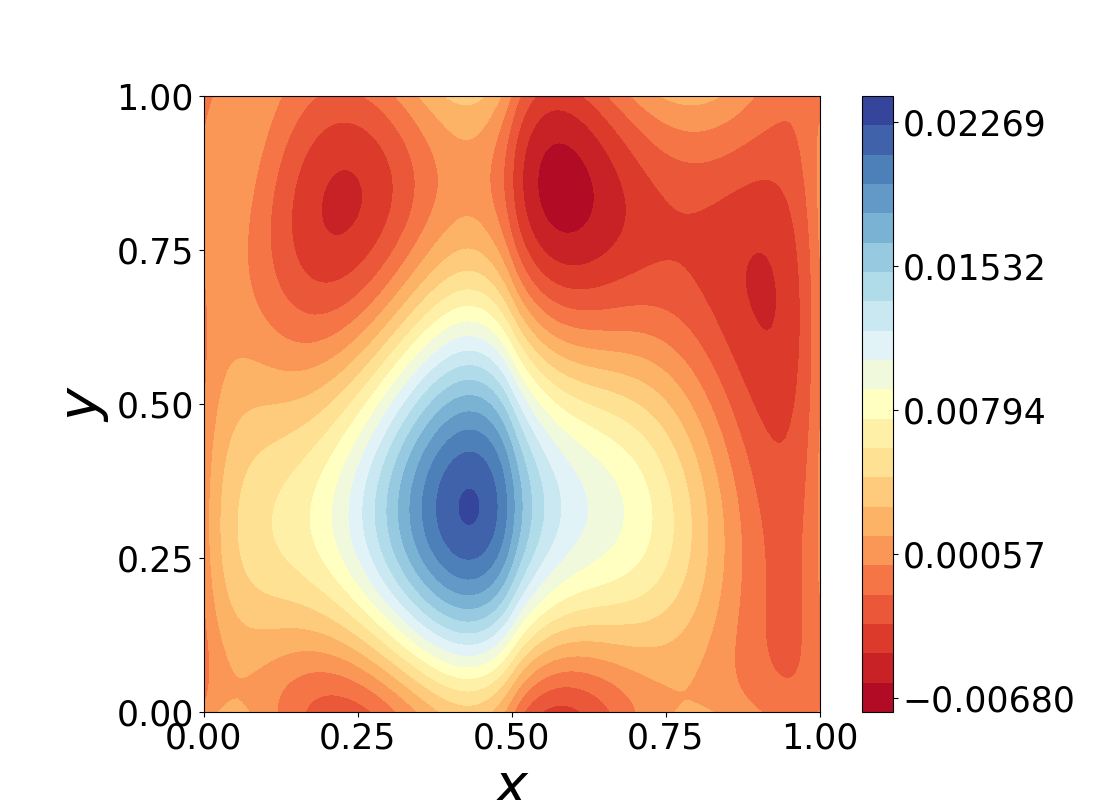}\\
     (e) $u_{\WTN} - u^{\exact}$ &(f)  $u_{\stf} - u^{\exact}$& (g)  $u_{\DRM} - u^{\exact}$ \\
    \end{tabular}
}
\caption{Darcy flow in Sec.~\ref{sec:onlyweak_v2}: (a) $u_{\WTN}$ obtained by Alg.\ref{alg:wtn}; (b) $u_{\stf}$ obtained by minimizing the strong form loss; (c) $u_{\DRM}$ obtained by minimizing the DRM loss; (d), (e), (f) show the pointwise errors $u_{\star} - u^{\exact}$ for $\star = \WTN, \stf$, and $\DRM$, respectively.}
\label{fig:results_onlyweak_v2}
\end{figure}

It is worth noting that in cases where the number of test basis functions is insufficient (e.g., $(M,N)=(100,50)$ and $(M,N)=(200,100)$), the results obtained from $u_{\WTN}$ are unsatisfactory. 
Moreover, for a fixed number of trial basis functions, the relationship between $e_{\WTN}$ and $N$ does not exhibit the expected linear scaling. 
These numerical experiments suggest that choosing $N$ equal to or slightly larger than $M$ provides a reasonable balance between accuracy and computational efficiency.
We observe that, when $N$ is sufficiently large in WTN, $u_{\WTN}$ is more accurate thant $u_{\DRM}$, and both of them perform better than $u_{\stf}$. This matches our expectation: since a strong solution does not exist for this problem, methods based on weak formulations, such as $\WTN$ and $\DRM$, would perform better than $\stf$. On the other hand, a more accurate quadrature is used in WTN, reducing the numerical integration errors in computing $u_{\WTN}$. 
The numerical solutions and corresponding errors of all three methods, using $M=200$ (and $N=300$ in WTN), are presented in Fig.~\ref{fig:results_onlyweak_v2}.

\subsection{Darcy flow with multiscale features}
\label{sec:high-fre}
\begin{figure}[!htb]
\centerline{
    \begin{tabular}{cc} 
    \includegraphics[width=0.5\textwidth]{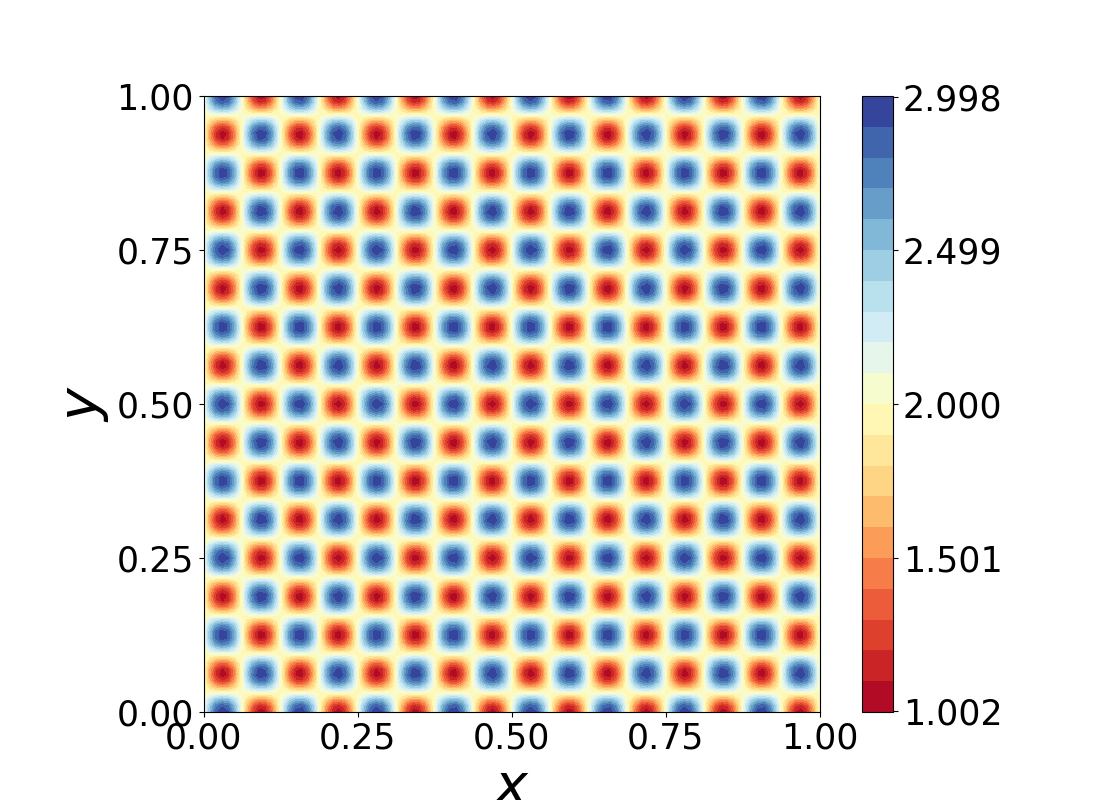} & 
    \includegraphics[width=0.5\textwidth]{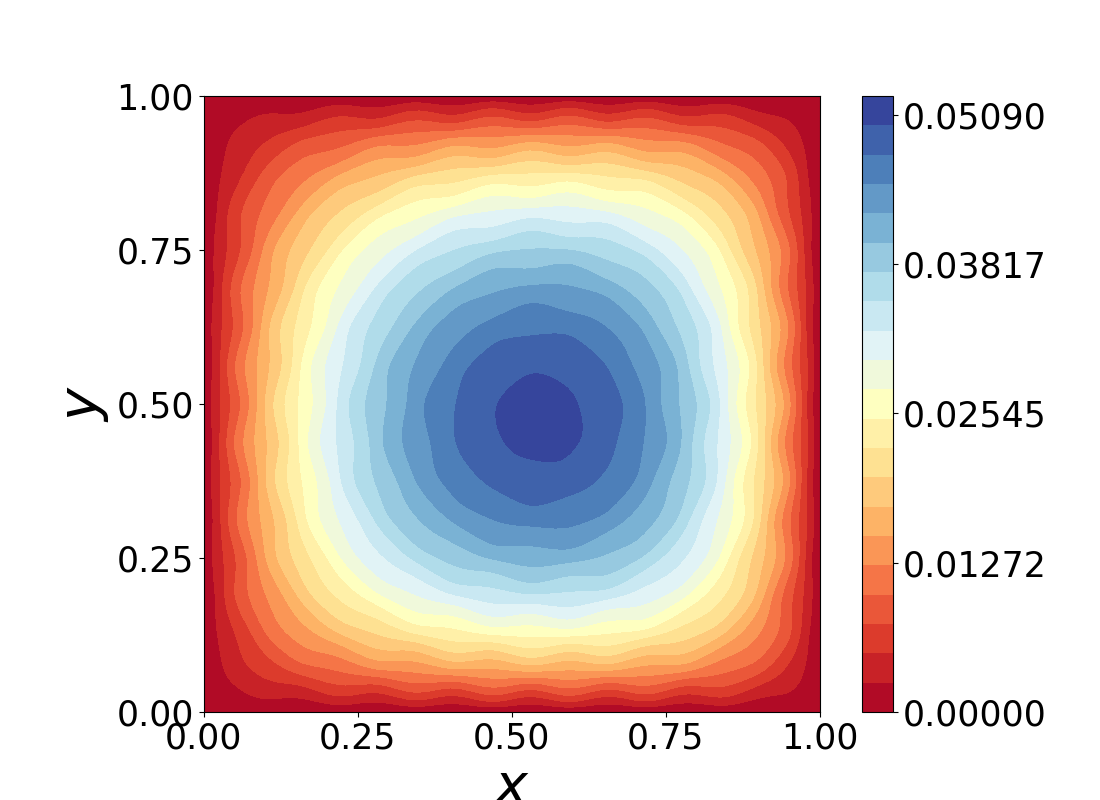} 
    \\
    (a) $\kappa(\bx)$ & (b)  $u^{\exact}$\\
    \end{tabular} \\
    }
\caption{Darcy flow in Sec.~\ref{sec:high-fre}: (a) permeability $\kappa(\bx)$; and (b) reference solution.}
\label{fig:results_high_fre_setup}
\end{figure}

Next, we consider again the two dimensional Darcy flow with homogeneous Dirichlet boundary condition, that is $\mathcal{L}[u]  = -\nabla \cdot (\kappa(\bx) \nabla u)$ and $\mathcal{B}[u]=u$ in~\eqref{eq:problem} with 
$\bx\in\Omega = [0,1]^2$, but change the permeability function $\kappa(\bx) = 2 + \sin(2\pi x/\epsilon) \cos(2\pi y/\epsilon)$ with $\epsilon = \frac{1}{8}$ and $f(\bx) = \sin(x) + \cos(y)$. 
As shown in Fig.~\ref{fig:results_high_fre_setup}(a), $\kappa(\bx)$ is multiscale due to the presence of high-frequency oscillations superimposed on a constant mean field, which brings the multiscale feature into the solution. 
The reference solution, as shown in Fig.~\ref{fig:results_high_fre_setup}(b), is computed by the finite element method on a $101\times 101$ mesh \footnote{FEM is implemented by FEniCS \cite{alnaes2015fenics}. The domain is discretized using a structured mesh consisting of rectangles, and the solution space is approximated using continuous Lagrange finite elements of degree one.}. 

As discussed in \cite{leung2022nh}, the PINN method fails to learn this solution in a finite training budget. 
From the perspective of training, this is because it takes much more efforts for a neural network to learn the high-frequency information due to it spectral bias \cite{rahaman2019spectral}.

\begin{figure}
\centerline{
    \begin{tabular}{ccc} 
    \includegraphics[width=0.4\textwidth]{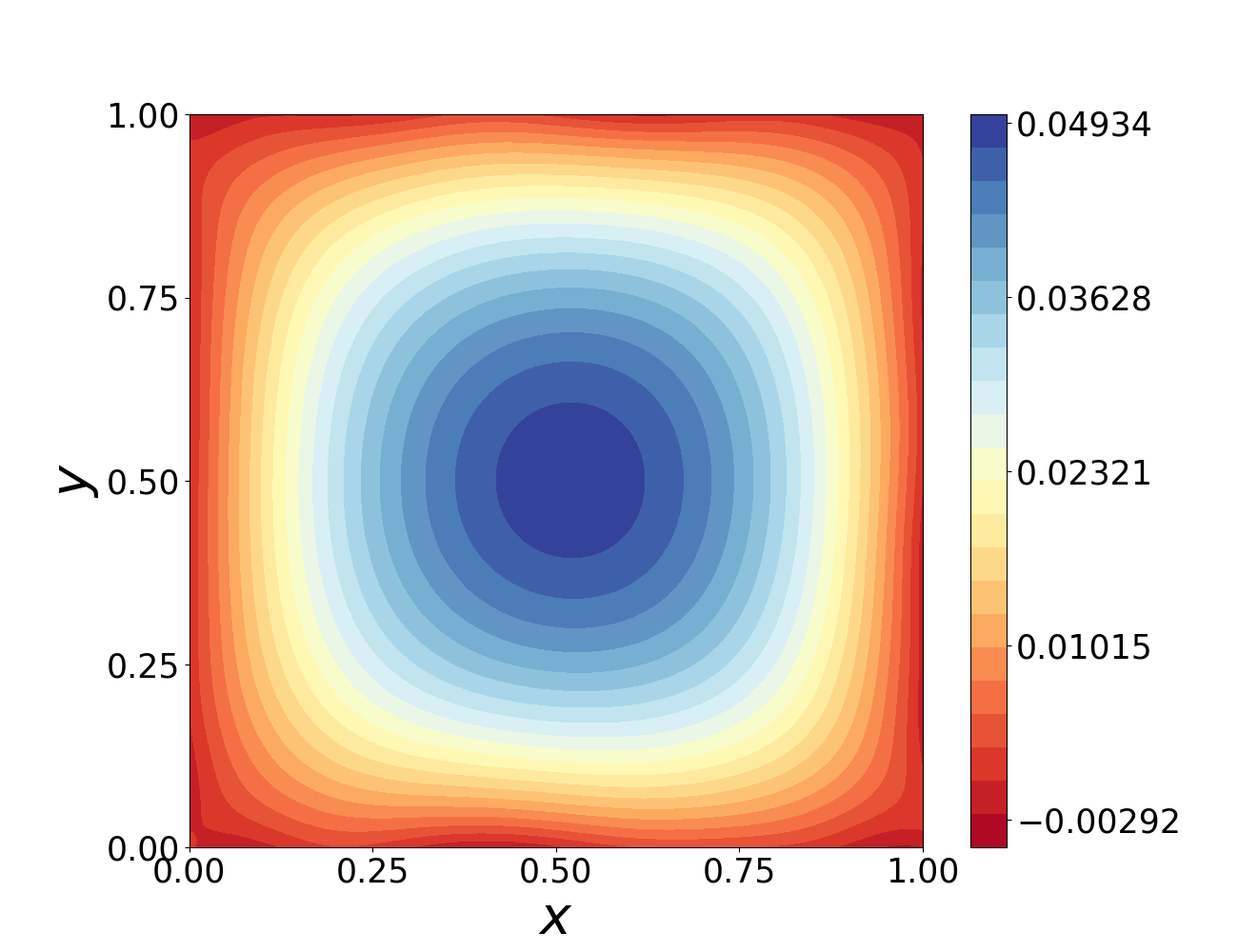} & 
    \includegraphics[width=0.4\textwidth]{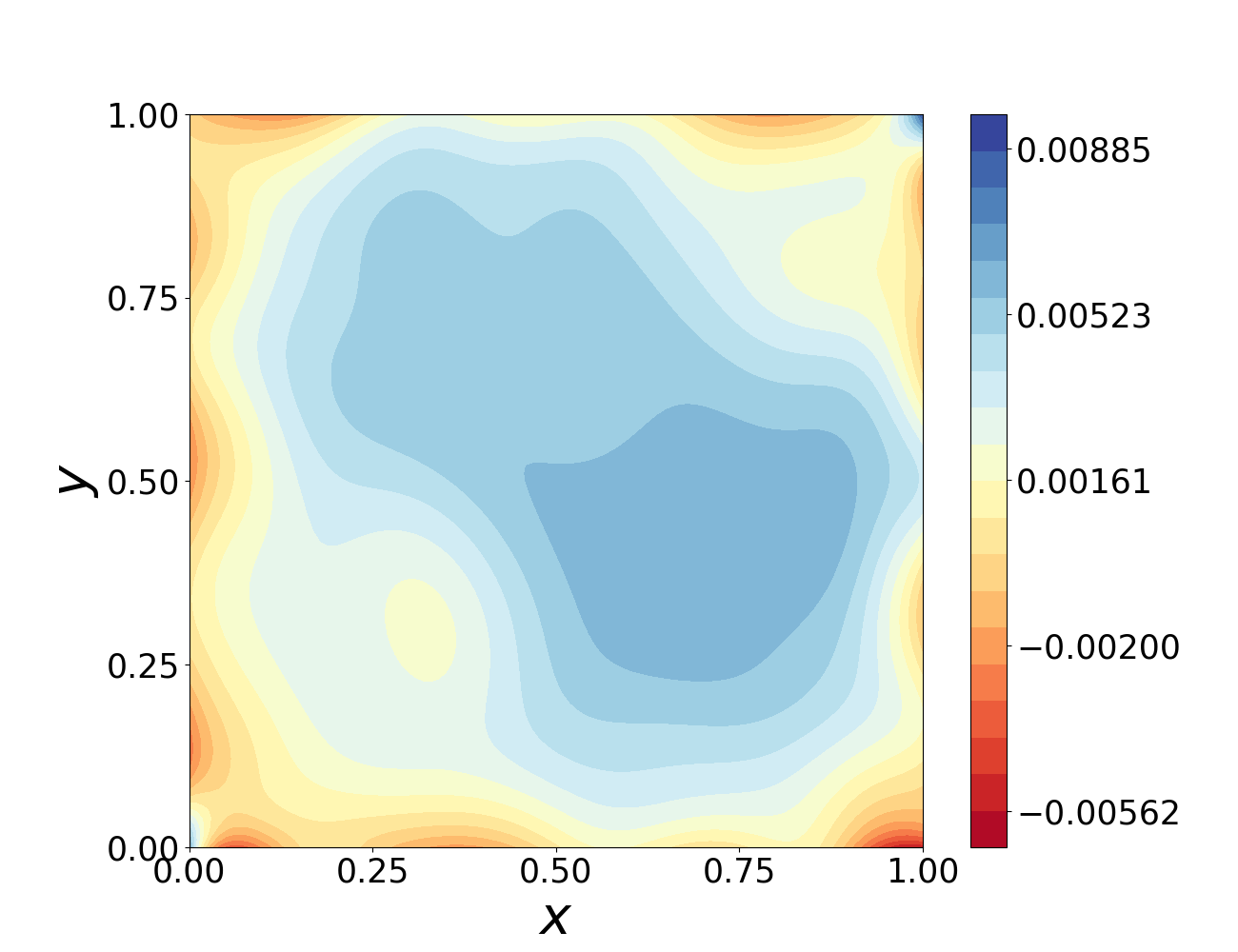} & 
    \includegraphics[width=0.4\textwidth]{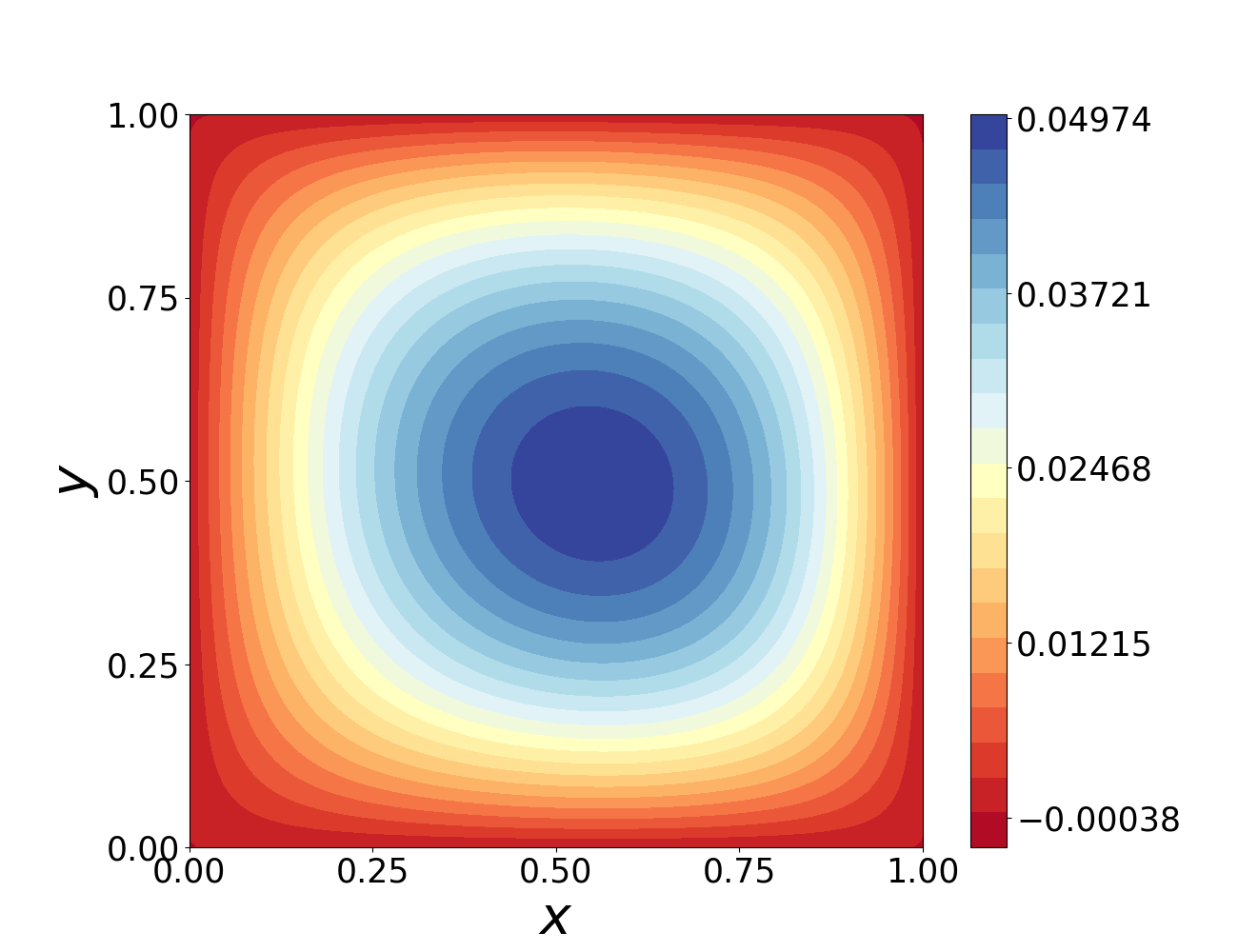}
    \\
    (a) $u_{\WTN}$ & (b)  $u_{\stf}$& (c) $u_{\DRM}$ \\
    \includegraphics[width=0.4\textwidth]{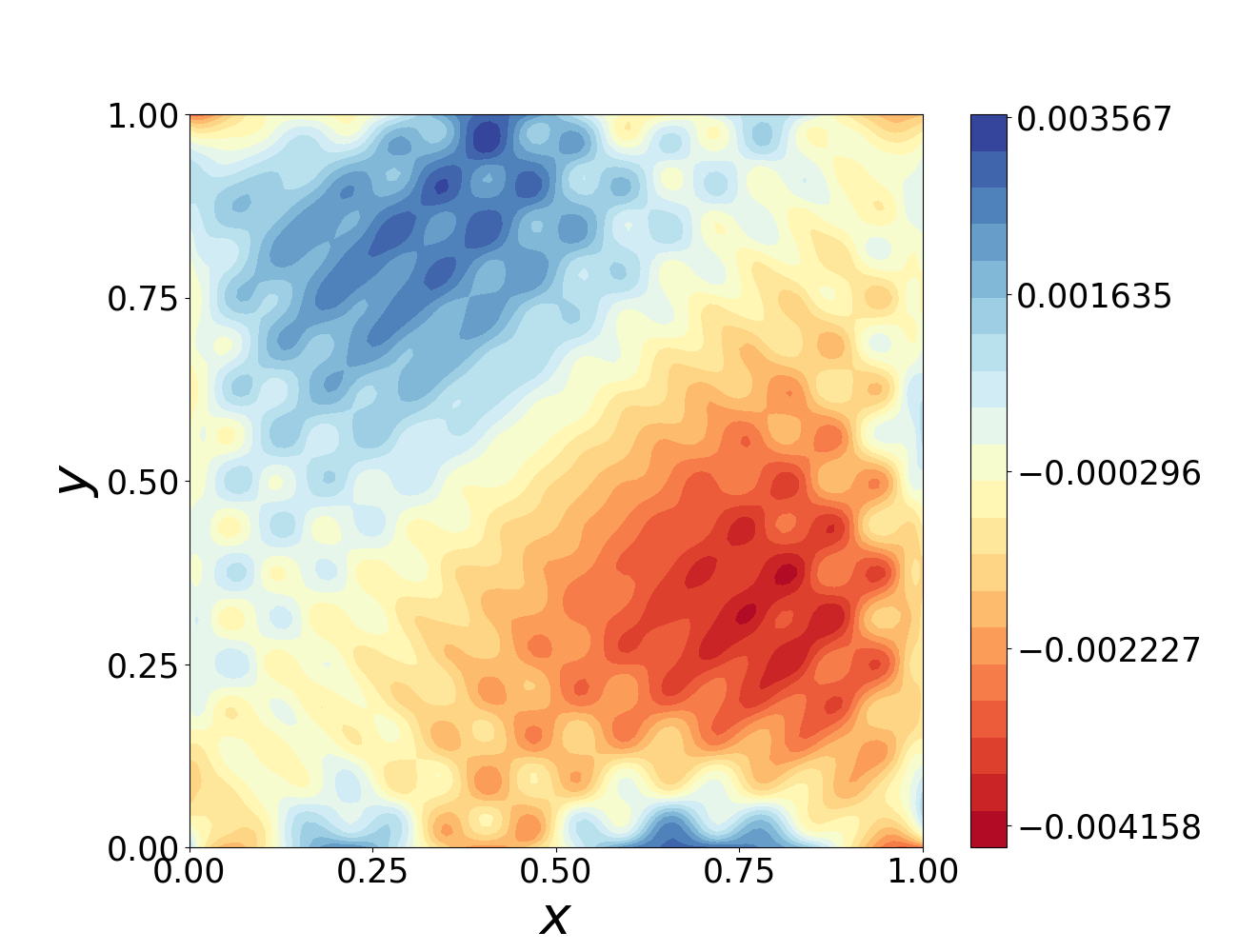} & 
    \includegraphics[width=0.4\textwidth]{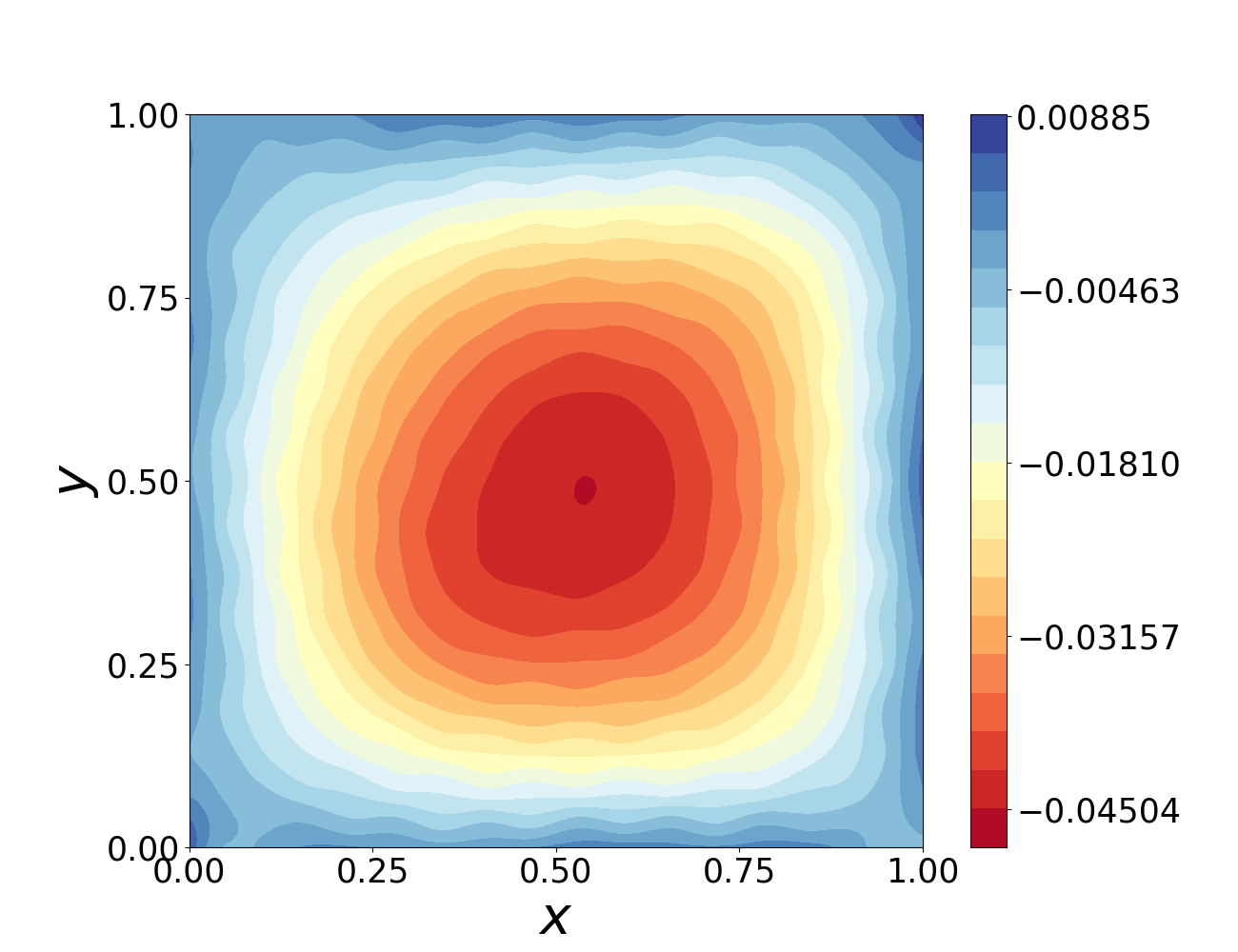} &
    \includegraphics[width=0.4\textwidth]{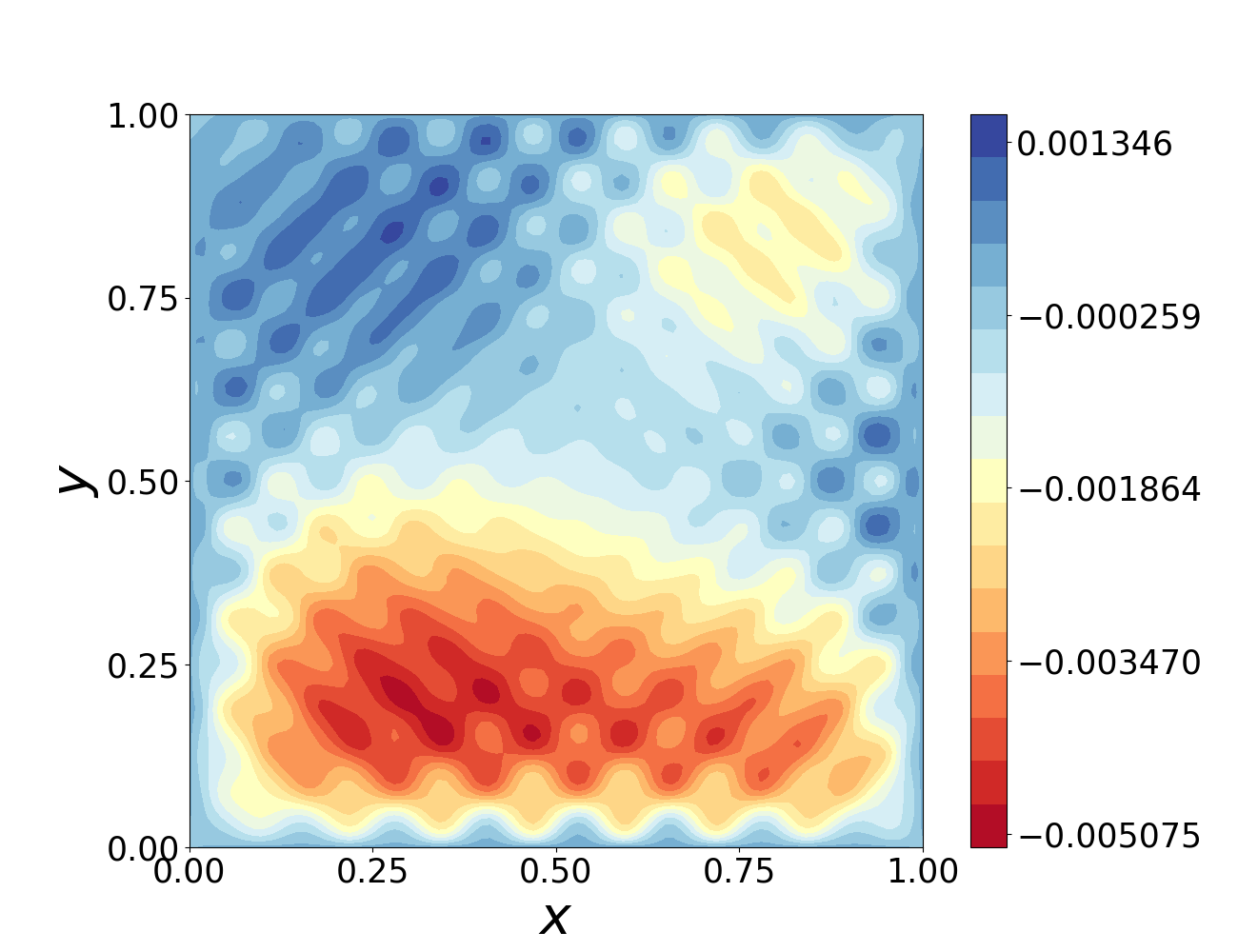}  
    \\
      (d) Pointwise error $u_{\WTN} - u^{\exact}$ &(e) Pointwise error $u_{\stf} - u^{\exact}$& (f)Pointwise error $u_{\DRM} - u^{\exact}$\\
    \end{tabular} \\
    }
\caption{Darcy flow in Sec.~\ref{sec:high-fre}: (a) $u_{\WTN}$ obtained by Alg. \ref{alg:wtn}; (b) $u_{\stf}$ obtained by strong form loss; (c) $u_{\DRM}$ obtained by DRM loss; (e),(f), (g) show pointwise errors $u_{\star} - u^{\exact}$ for $\star =\WTN, \stf$, and $\DRM$, respectively.}
\label{fig:results_high_fre}
\end{figure}

We then employ TransNet and compare the results obtained by minimizing the loss functions of WTN, SF, and DRM. For this comparison, we set $M = 2000$, shape parameter $\gamma = 1$ in all neural basis functions across all three methods, we additionally take $N=2000$ and $\sigma_i =0.05$ in WTN. The numerical results and corresponding pointwise errors are presented in Fig.~\ref{fig:results_high_fre}. 
It is seen that both $u_{\WTN}$ and $u_{\DRM}$ achieve higher accuracy than $u_{\stf}$. Specifically, the relative $L_2$ errors for $u_{\WTN}$ and $u_{\DRM}$ are $6.11\%$ and $7.12\%$, respectively, indicating that WTN marginally outperforms DRM in this setting. 

\begin{figure}[!htb]
\centerline{
    \begin{tabular}{cc} 
    \includegraphics[width=0.5\textwidth]{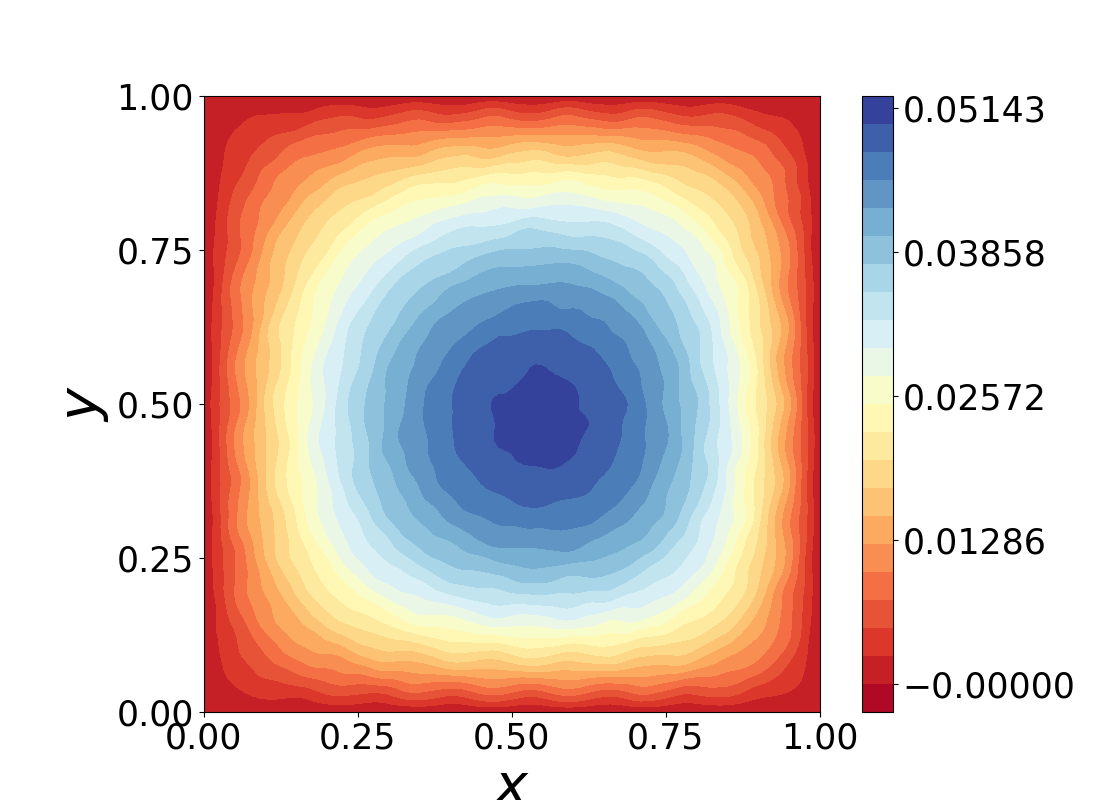} & 
    \includegraphics[width=0.5\textwidth]{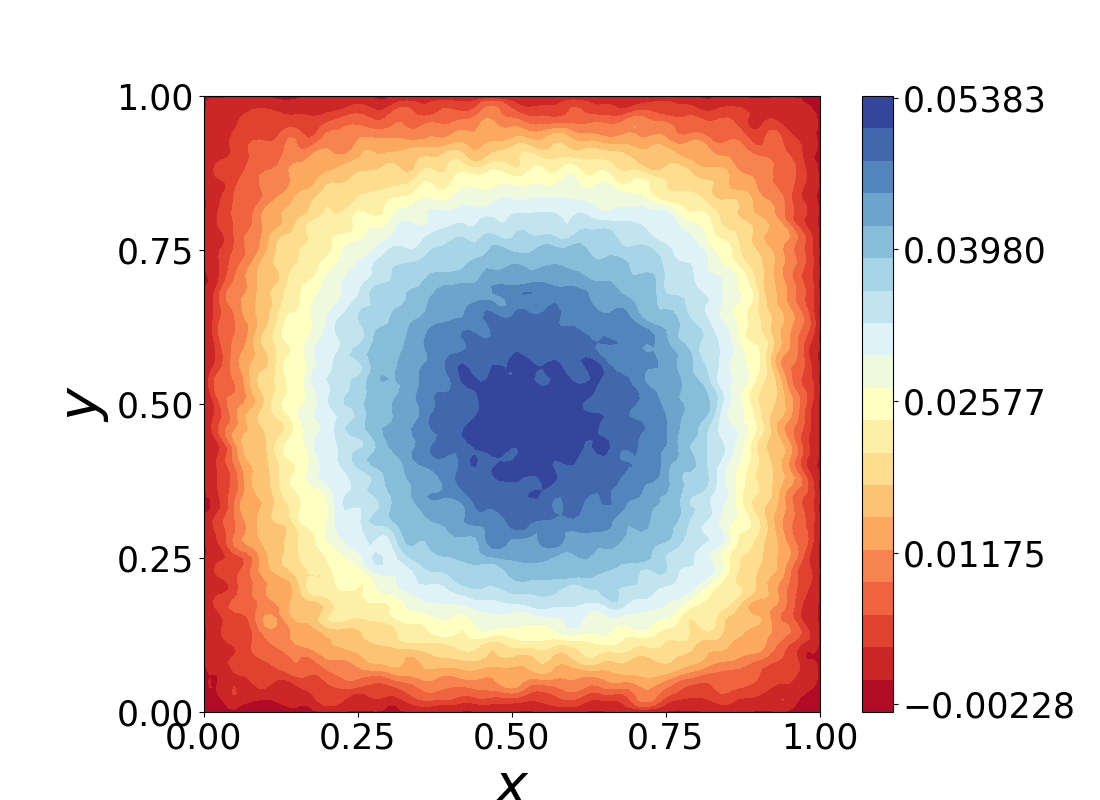}
    \\
    (a) $u_{\F-\WTN}$ & (b) $u_{\F-\DRM}$ \\
    \includegraphics[width=0.5\textwidth]{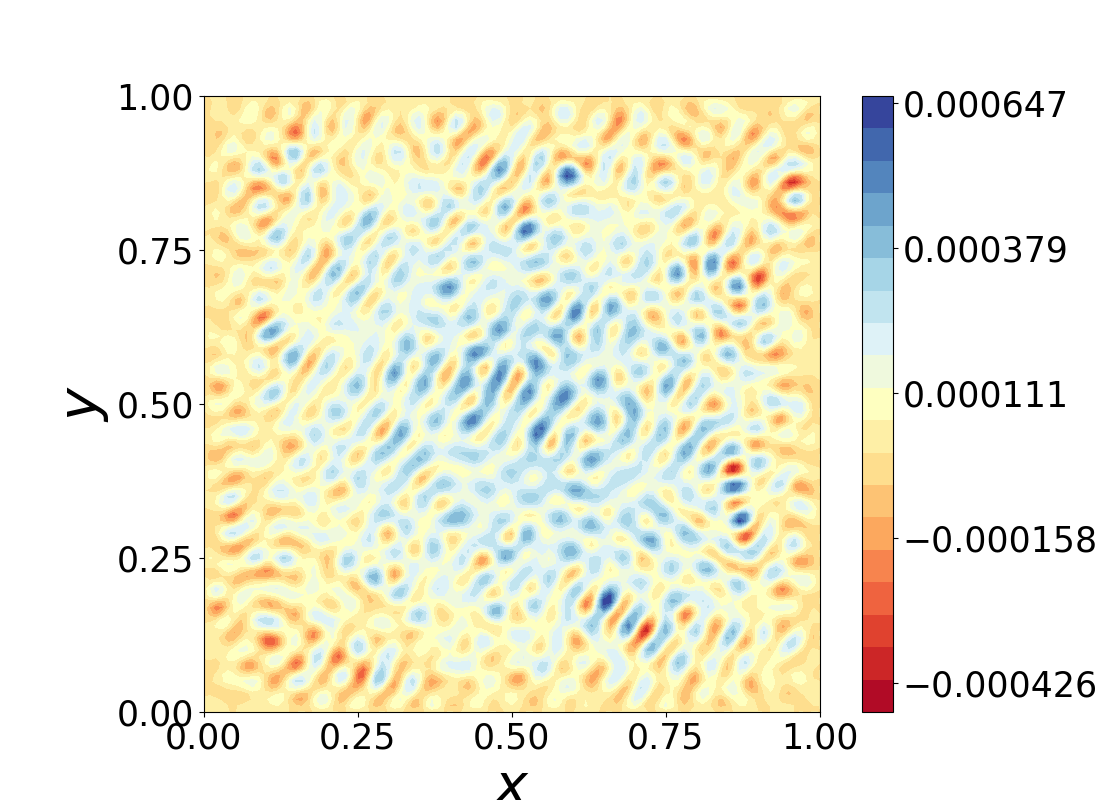} & 
    \includegraphics[width=0.5\textwidth]{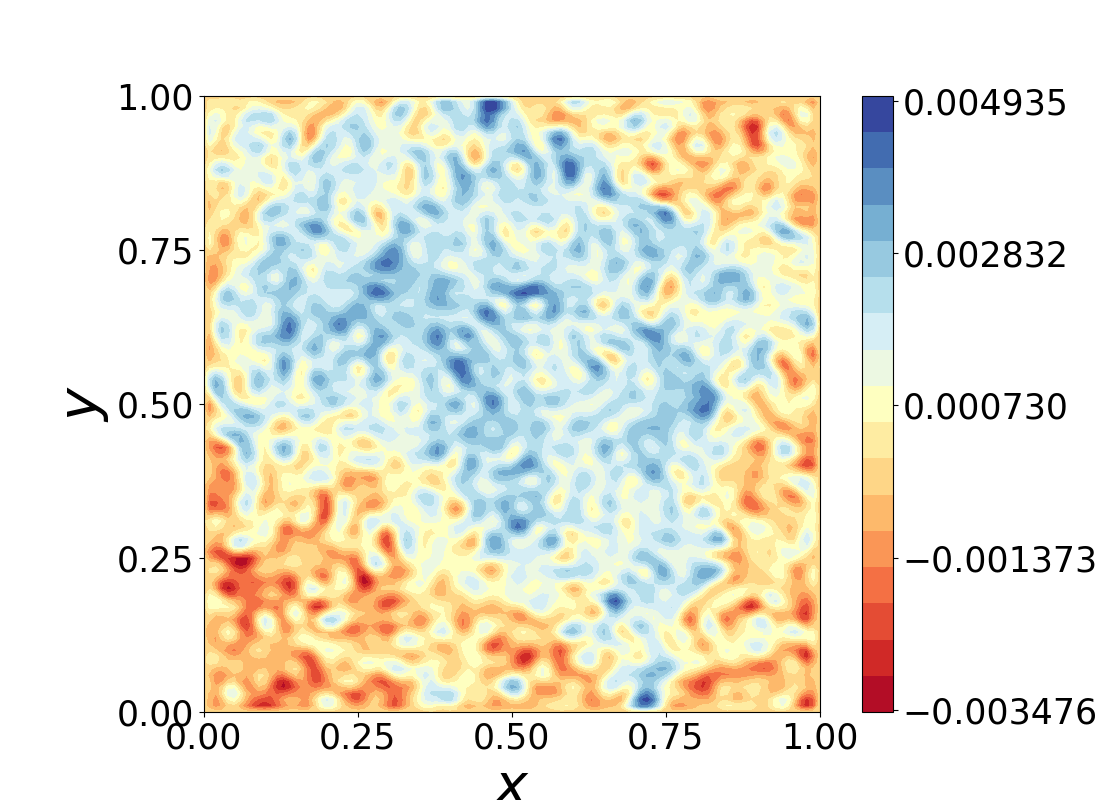}  
    \\
      (c)  $u_{\F-\WTN} - u^{\exact}$ & (d)$u_{\F-\DRM} - u^{\exact}$\\
    \end{tabular} \\
    }
\caption{Darcy flow in Sec.~\ref{sec:high-fre}: (a) $u_{\F-\WTN}$; (b) $u_{\F-\DRM}$; (c), (d) Pointwise error $u_{\star} - u^{\exact}$ for $\star =\F-\WTN$ and $\F-\DRM$, respectively.}
\label{fig:results_high_fre_fourier}
\end{figure}

Given the multiscale nature of $\kappa(\bx)$, we further use F-WTN (see Sec.~\ref{sec:f-wtn}) to solve this problem. For the Fourier feature mapping, the dimension of $\mathfrak{B}$ is taken to be $P = 64$ and $d=2$.
The hyper-parameter $\sigma_B$ controls the frequency range encoded by the Fourier features. In the absence of prior knowledge, we use a mixture of $\sigma_B \in \{1, 3\}$ to capture both low and higher frequencies. 
We take $\{r_j\}_{j=1}^M$ to be i.i.d. and uniformly distributed on $[0,9]$.
By minimizing the WTN and DRM loss functions, we obtain numerical solutions denoted by $u_{\F-\WTN}$ and $u_{\F-\DRM}$, respectively. The results and corresponding pointwise errors are shown in Fig.~\ref{fig:results_high_fre_fourier}. 
Both Fourier-enhanced approaches demonstrate superior performance to their standard TransNet counterparts. In particular, $u_{\F-\DRM}$ achieves a relative $L_2$ error of $5.74\%$, improving upon $u_{\DRM}$ of $7.12\%$. 
Meanwhile, $u_{\F-\WTN}$ attains exceptional accuracy with a relative $L_2$ error of just $0.58\%$, outperforming both $u_{\WTN}$ (error of $6.11\%$) and $u_{\DRM}$ (error or $7.12 \%$ ) by an order of magnitude. 
For ease of comparison, we summarize all the relative errors in Tab.~\ref{tab:re_highfre}.

\begin{table}[!htp]
	\caption{Darcy flow in Sec. \ref{sec:high-fre}: relative $L_2$ errors of $u_{\WTN}$, $u_{\stf}$, $u_\DRM$,  $u_{\F-\WTN}$,  $u_{\F-\DRM}$.}
    \vspace{0.1cm}
	\centering
    \adjustbox{max width=\textwidth}{
    \begin{tabular}{ccccc}
\hline
$e_{\WTN}$ & $e_{\DRM}$ & $e_{\stf}$ & $e_{\F-\WTN}$ & $e_{\F-\DRM}$   \\
$6.11\times 10^{-2} $ & $7.12\times 10^{-2}$ & $8.67\times 10^{-1}$ & $5.78\times 10^{-3}$&  $5.74\times 10^{-2}$ \\ 
\hline
\end{tabular}
}
	\label{tab:re_highfre}
\end{table} 

\subsection{Darcy flow with channelized permeability fields}
\label{sec:channel}

\begin{figure}[!htb]
\centerline{
    \begin{tabular}{cc} 
    \includegraphics[width=0.44\textwidth]{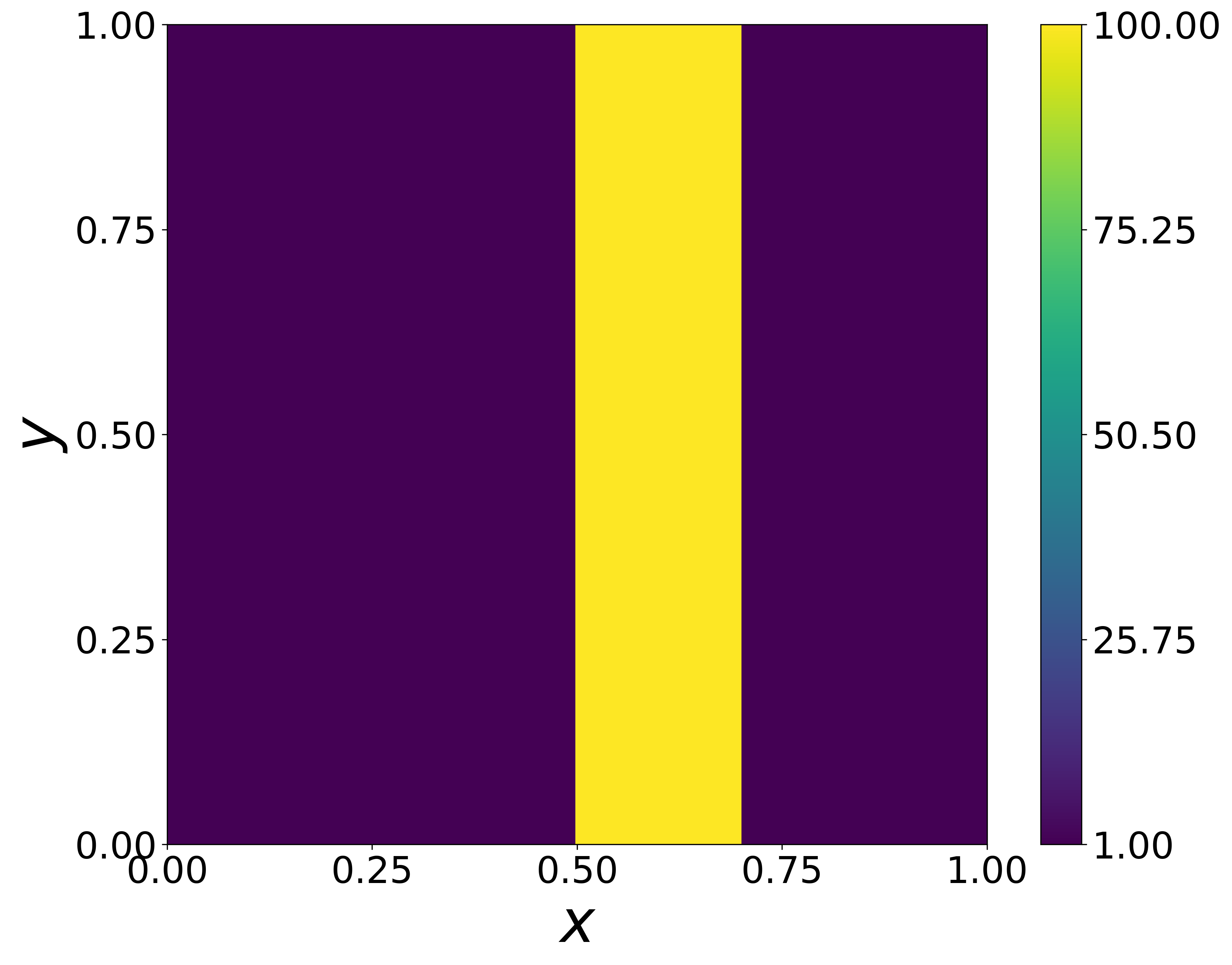} 
    & 
    \includegraphics[width=0.5\textwidth]{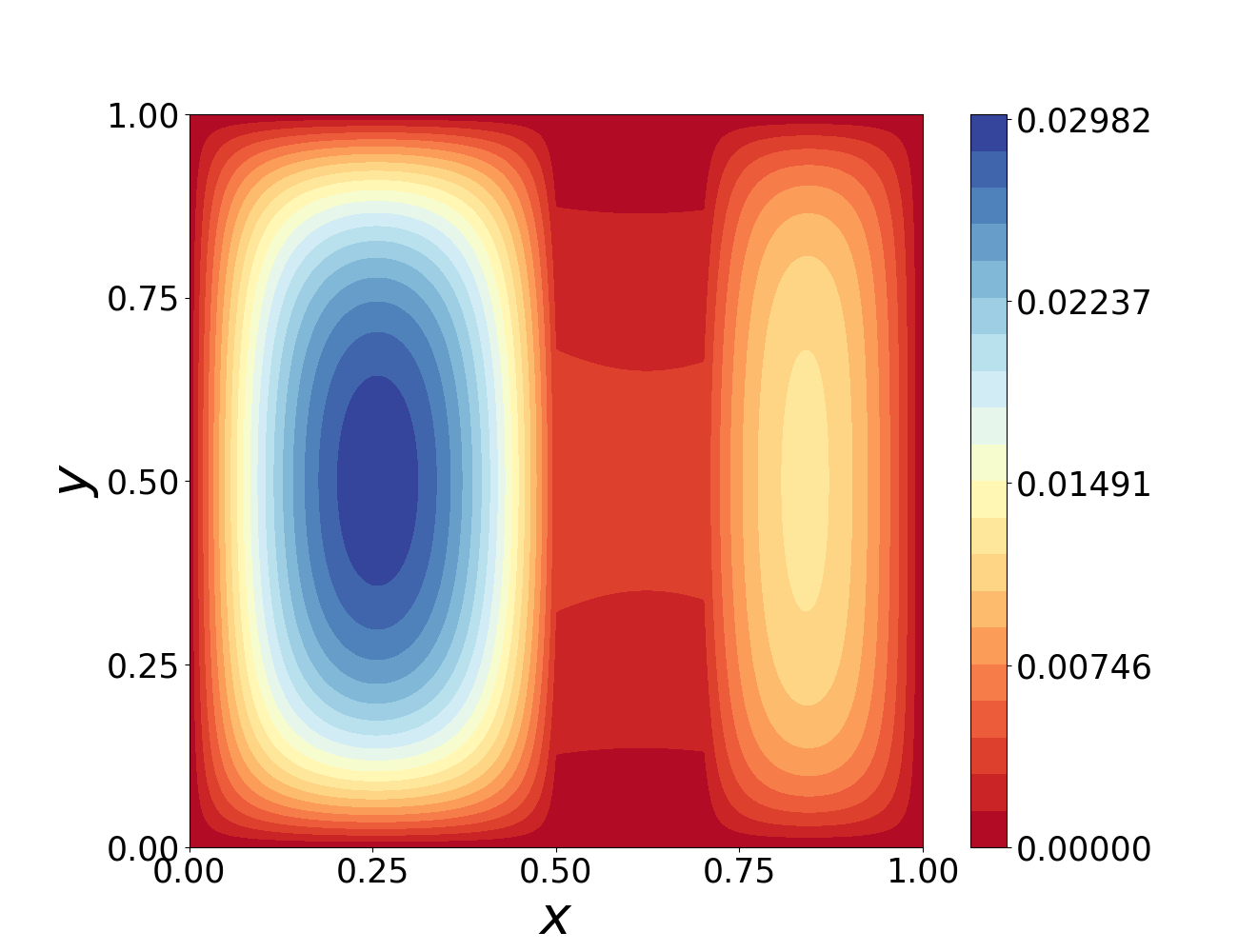}
    \\
    (a) $\kappa(\bx)$ & (b) $u^{\exact}(\bx)$
    \end{tabular} \\
    }
\caption{Darcy flow in Sec.~\ref{sec:channel}: (a) permeability $\kappa(\bx)$; and (b) reference solution $u^{\exact}$.}
\label{fig:channel}
\end{figure}

Next, we consider the Darcy flow problem with a discontinuous permeability function (see Fig.~\ref{fig:channel}(a)) 
\[
\kappa(\bx) = \begin{cases}
    1 \,, &\bx \in (0,0.5) \times (0,1) \cup  (0.7,1) \times (0,1)\,, \\
    100\,, &\bx \in [0.5,0.7] \times (0,1) \,.
\end{cases}
\]
The reference solution (see Fig.\ref{fig:channel}(b)) is obtained by the finite element method on a $257\times 257$ mesh\footnote{FEM is implemented by FEniCS \cite{alnaes2015fenics}. The domain is discretized using a structured mesh consisting of rectangles, and the solution space is approximated using continuous Lagrange finite elements of degree 1 as implemented in FEniCS.}.

\begin{figure}[!htb]
\centerline{
    \begin{tabular}{cc} 
    \includegraphics[width=0.44\textwidth]{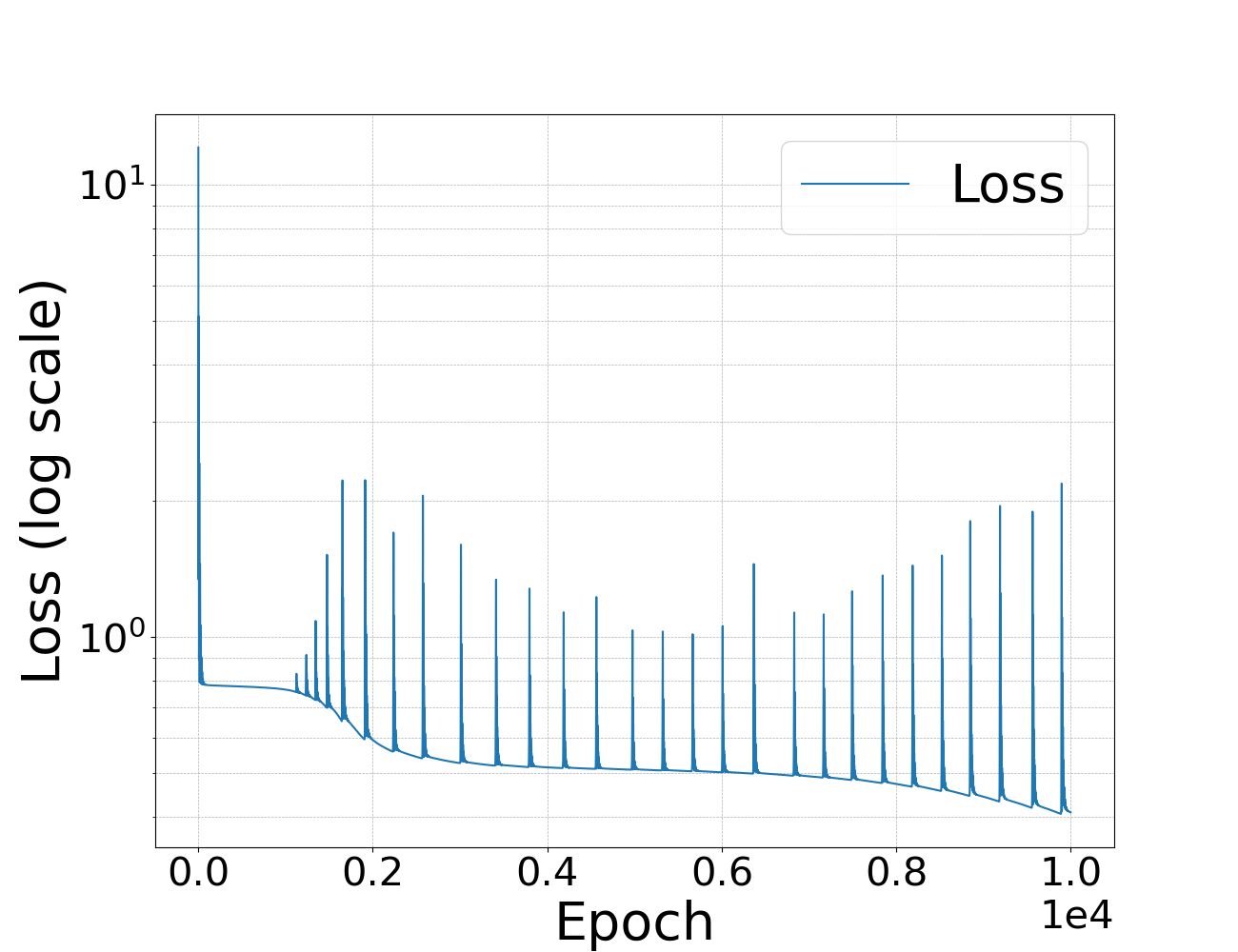} 
    & 
    \includegraphics[width=0.5\textwidth]{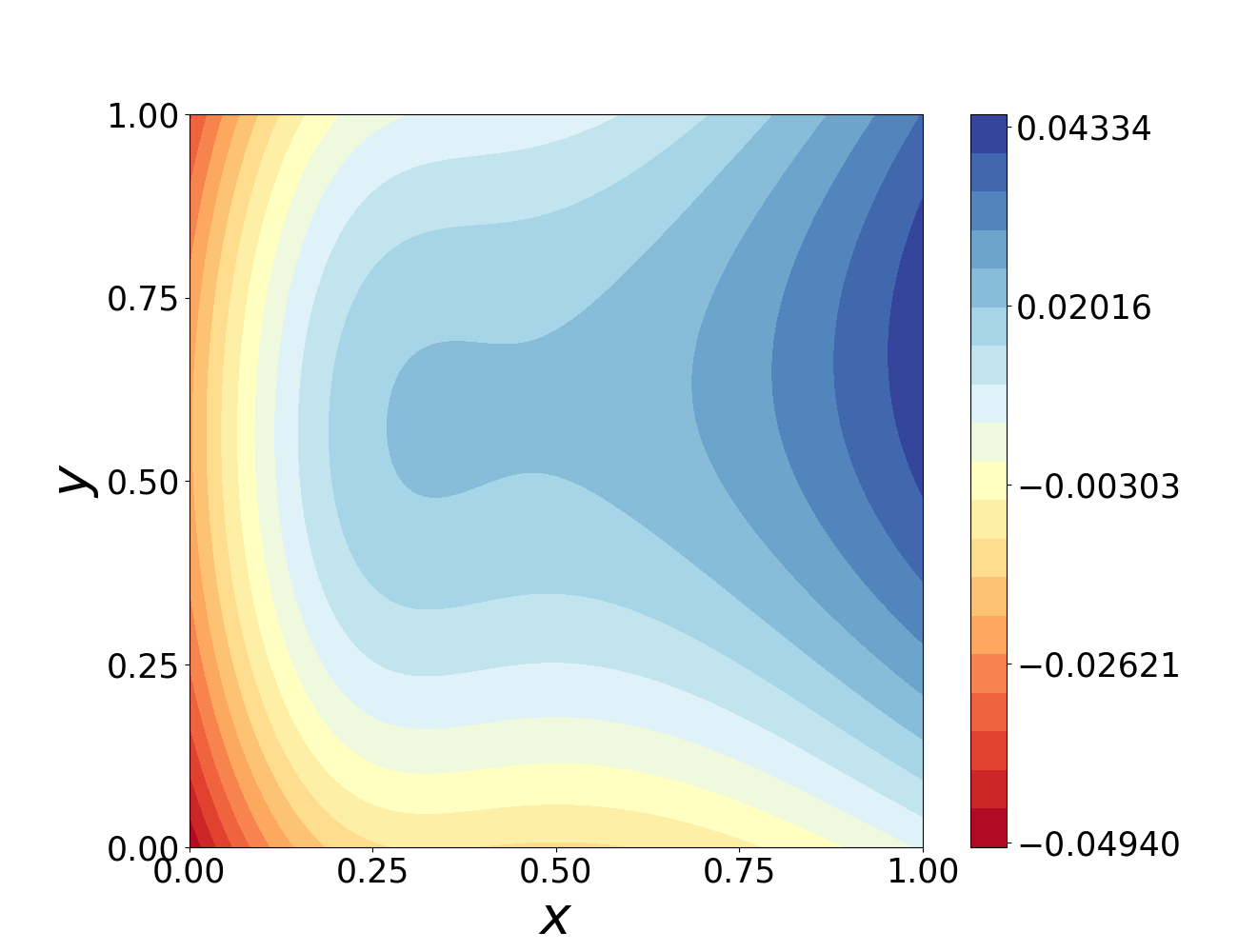}
    \\
    (a) Training history & (b) PINN solution
    \\
    \end{tabular} \\
    }
\caption{Darcy flow in Sec.~\ref{sec:channel}: (a) training loss of PINN; and (b) numerical solution obtained by PINN.}
\label{fig:channel_training}
\end{figure}

Solving such a problem using PINN is challenging due to the high contrast permeability. In fact, when using a fully-connected neural network (FCNN) with three hidden layers, each with 50 neurons and $\tanh(\cdot)$ as the activation function, training of the PINN can be slow and the resulting solution is inaccurate, see Fig.~\ref{fig:channel_training} for the training loss history and the numerical solution.

To resolve the sharp gradients induced by the high-contrast permeability field, we use the PoU-WTN method developed in Sec.~\ref{sec:pou_weaktransnet}.
The domain is divided into three nonoverlapping subdomains, with the decomposition aligned along the discontinuities in the permeability field. In each subdomain, we select 200 trial basis functions with the shape parameters being 1.  The number of test functions is take as 600.
The boundary weight $\beta$ and the interface weight $\lambda$ are both set to 1.
The relative $L_2$ error of the obtained $u_{\pou-\WTN}$ is $4.05\times 10^{-2}$, which successfully approximates the solution (see Fig.~\ref{fig:channel}). 

\begin{figure}[!htb]
\centerline{
    \begin{tabular}{cc} 
    \includegraphics[width=0.5\textwidth]{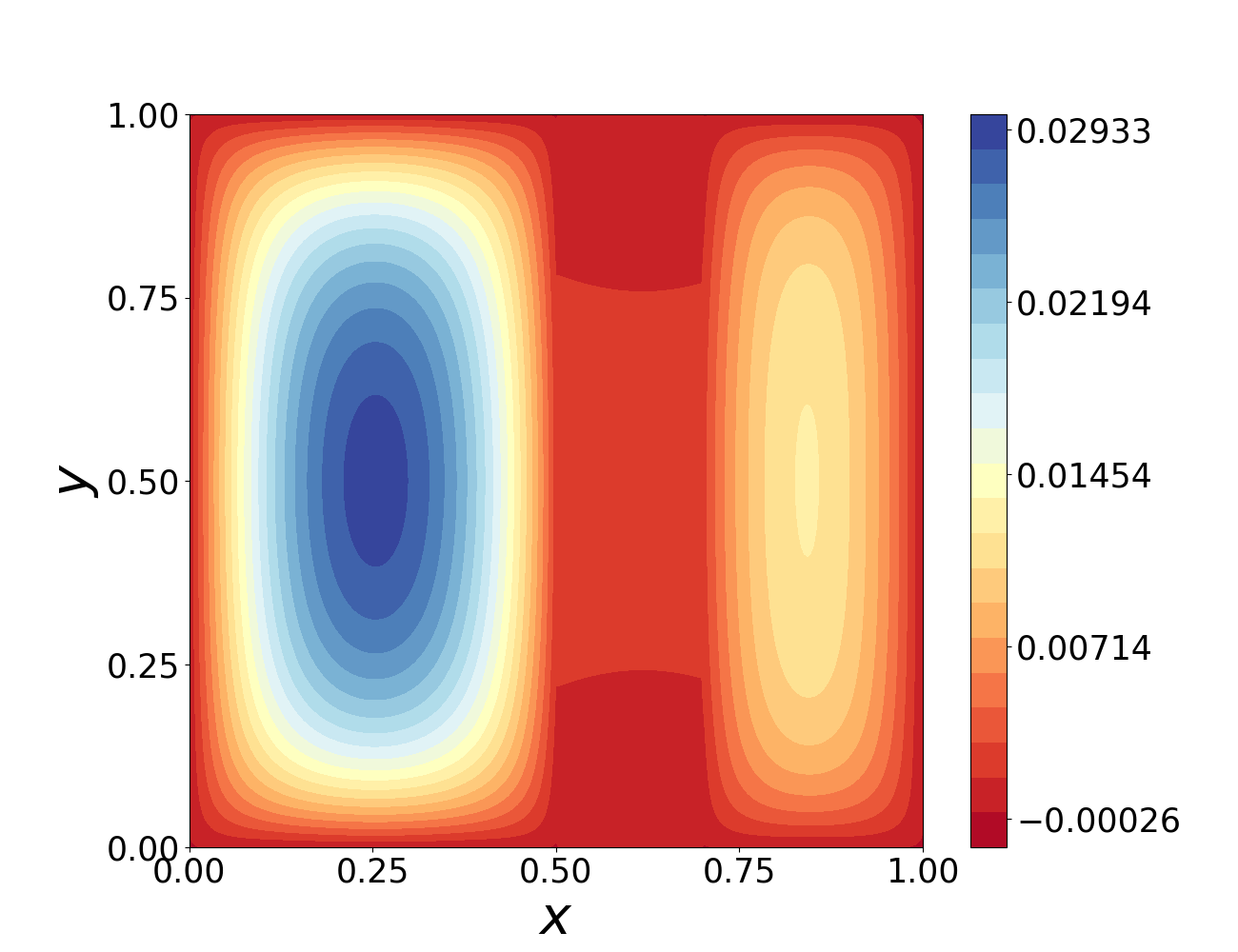} 
    & 
    \includegraphics[width=0.5\textwidth]{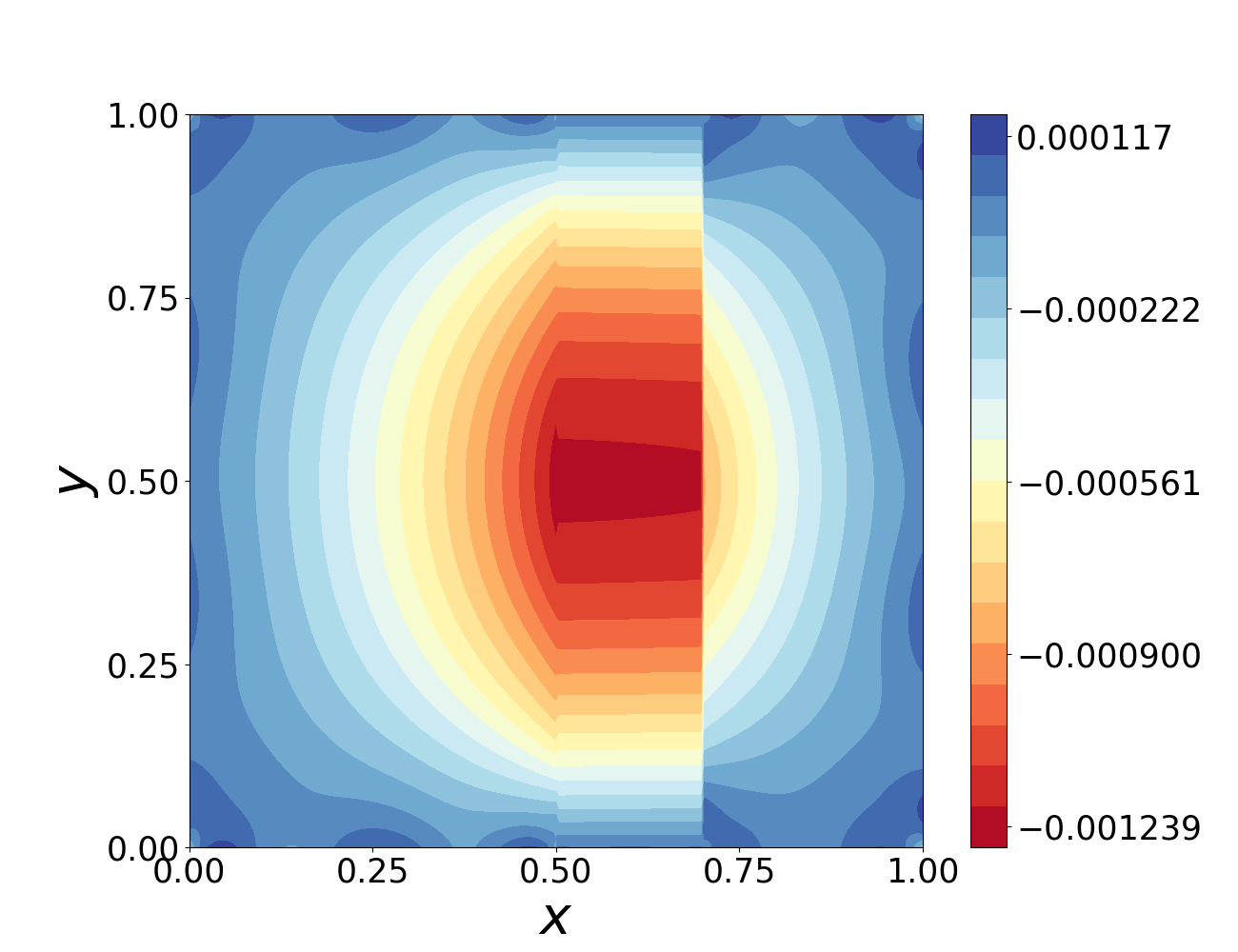}\\  
    (c) $u_{\pou-\WTN}$&(d) $u_{\pou-\WTN}-u^{\exact}$ \\
    \end{tabular} \\
    }
\caption{Darcy flow in Sec.~\ref{sec:channel}: (a) $u_{\pou-\WTN}$ approximation; and (b) pointwise error $u_{\pou-\WTN} - u^{\exact}$.}
\label{fig:channel2}
\end{figure}

\subsection{Poisson equation with sharp-gradient solutions}
\label{sec:vpinn}
\begin{figure}[!htb]
\centerline{
\begin{tabular}{cc}
    \includegraphics[width=0.6\linewidth]{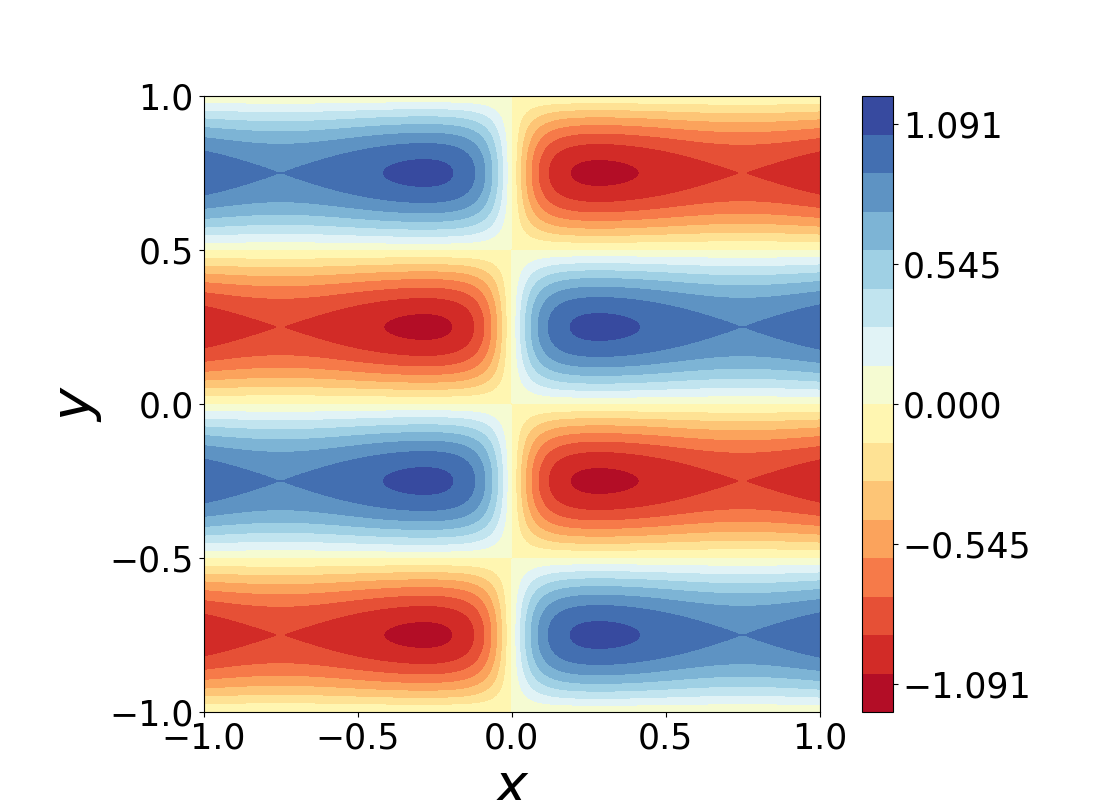}
    & 
    \includegraphics[width = 0.38\linewidth]{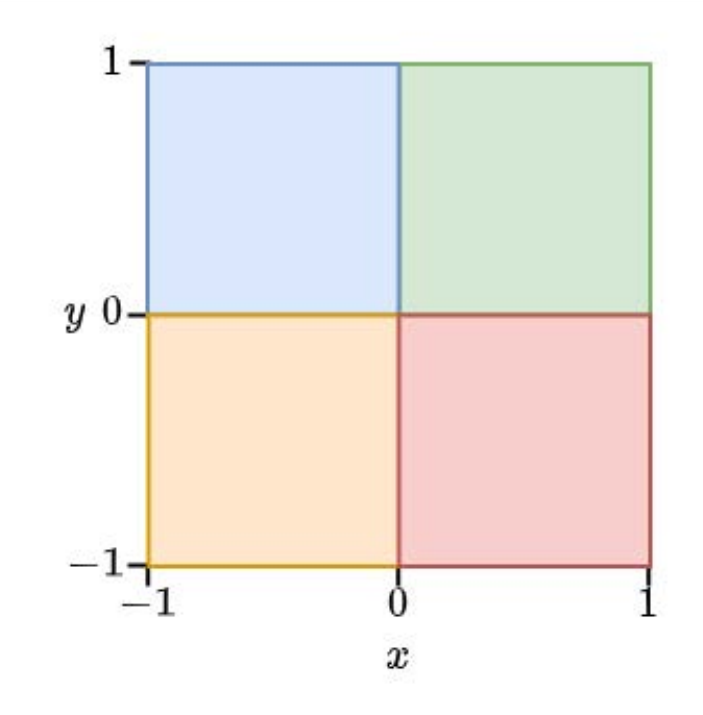}\\
    (a) & (b) 
    \end{tabular} \\}
    \caption{Poisson equation in Sec.\ref{sec:vpinn}: (a) Exact solution $u^\exact$; (b) The domain decomposition used for PoU.}
    \label{fig:vpinn-setting}
\end{figure}

We next consider the two-dimensional Poisson equation with Dirichlet boundary condition, that is, $\mathcal{L}[u] = -\Delta u$ and $\mathcal{B}[u]=u$, the source term $f$ and the boundary data $g$ are chosen such that 
the exact solution is prescribed by 
\[
u^{\exact}(x,y) = (0.1 \sin (2\pi x) + \tanh(10x)) \times \sin(2\pi y)\,.
\]
The exact solution is shown in Fig.~\ref{fig:vpinn-setting}(a). 

\begin{figure}[!htb]
\centerline{
    \begin{tabular}{cc}
    \includegraphics[width=0.5\textwidth]{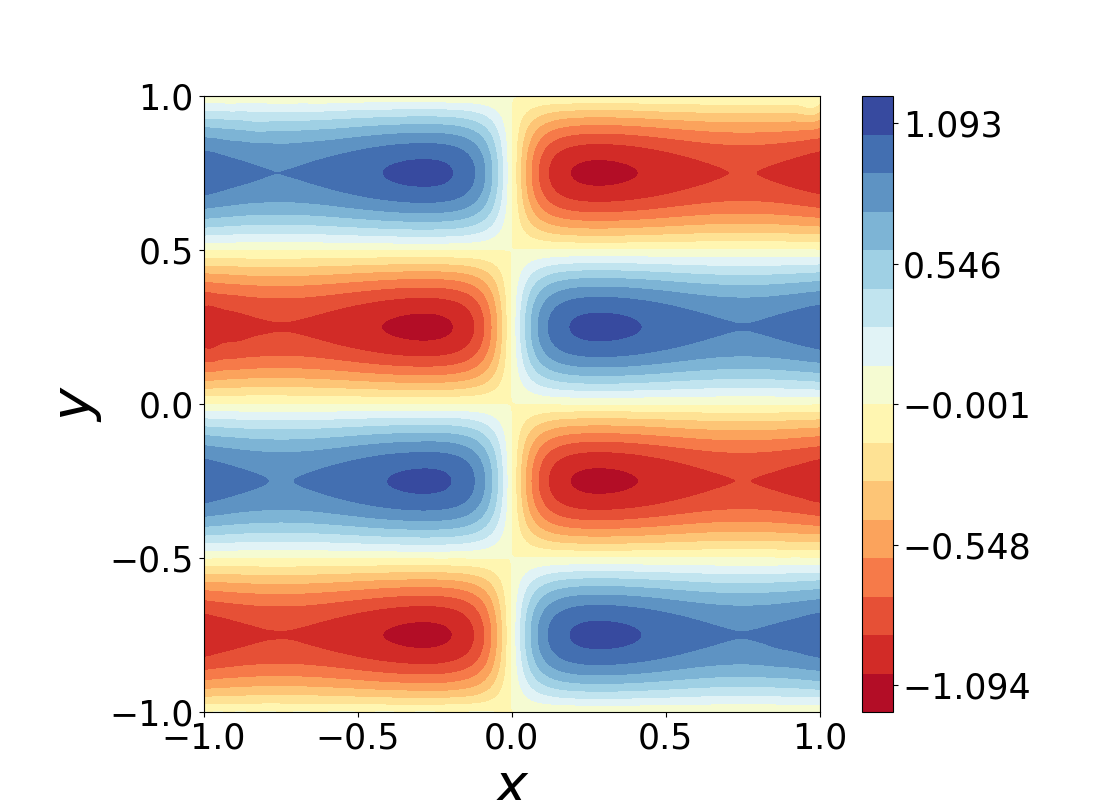} &
    \includegraphics[width=0.5\textwidth]{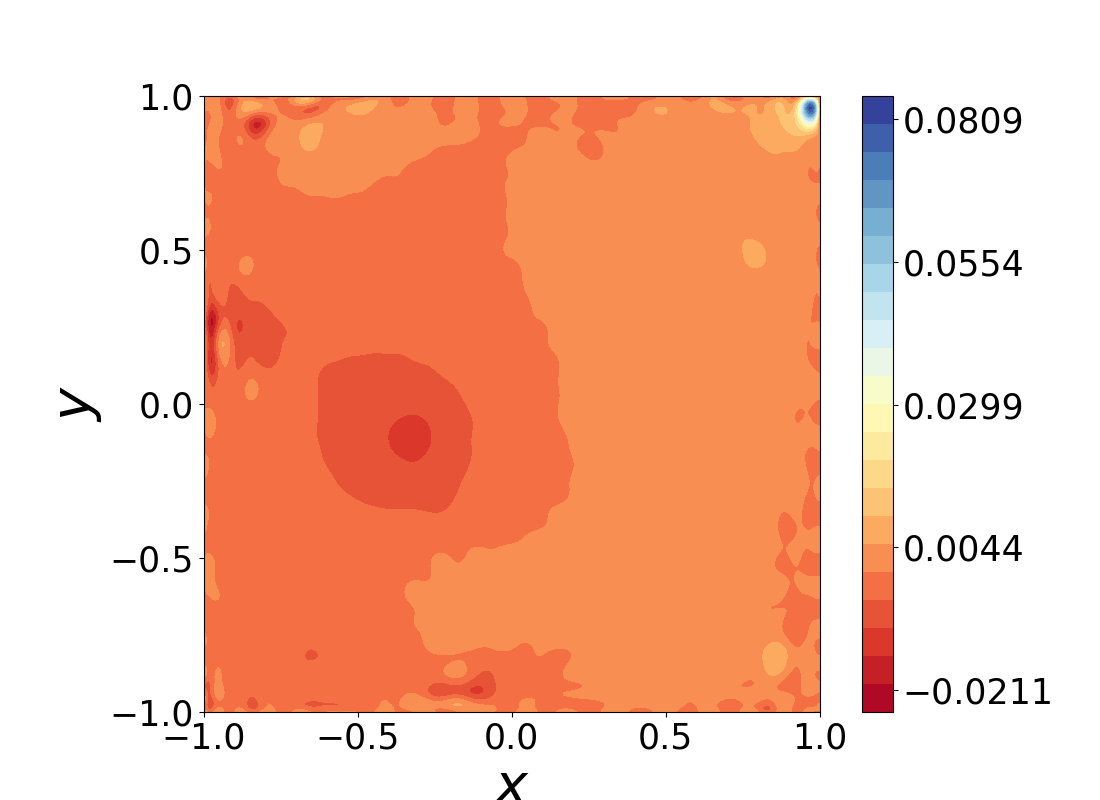}\\
    (a) $u_{\WTN}$ & (b) $u_{\WTN} - u^{\exact}$ \\ 
    \includegraphics[width=0.5\textwidth]{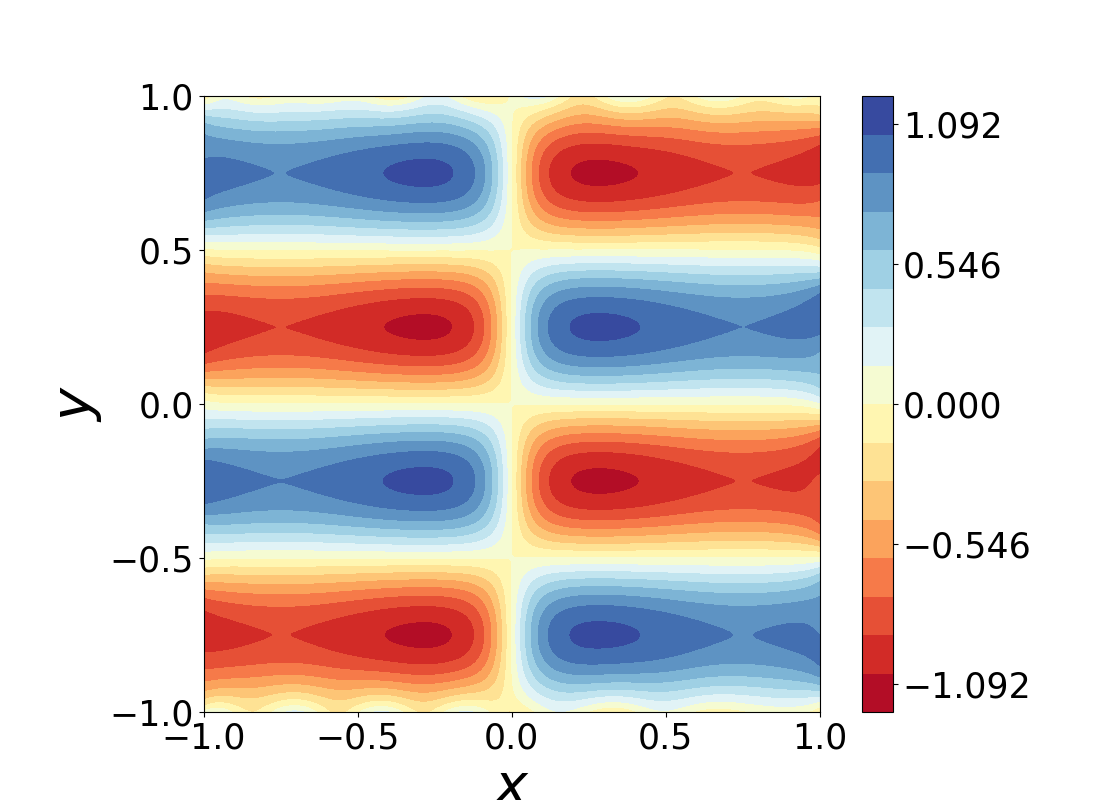} &
    \includegraphics[width=0.5\textwidth]{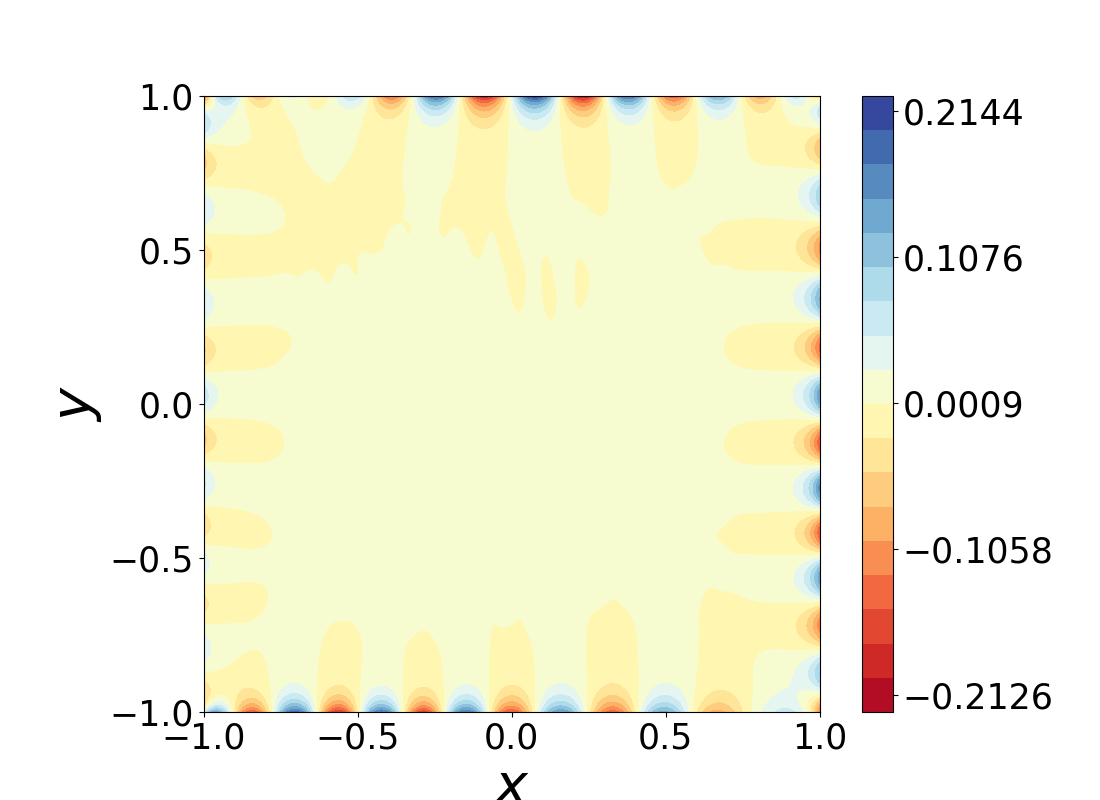} \\
    (c) $u_{\stf}$ & (d) $u_{\stf} - u^{\exact}$ \\  
     \includegraphics[width=0.5\textwidth]{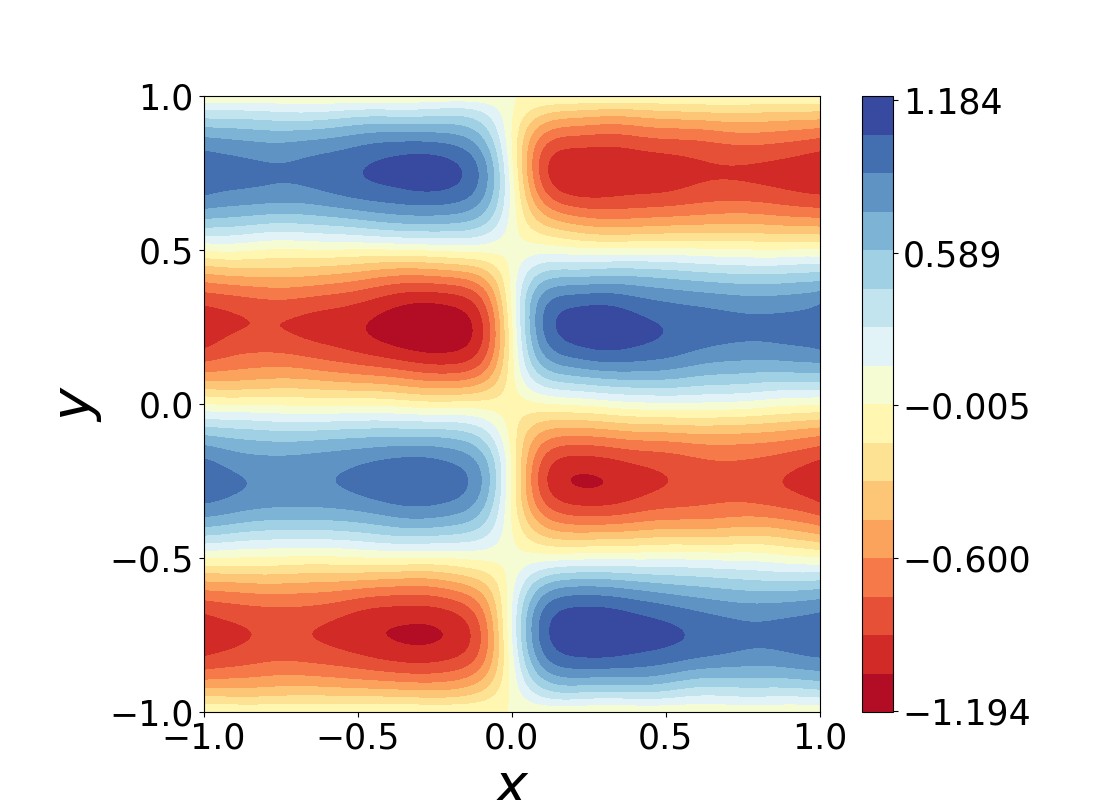} &
    \includegraphics[width=0.5\textwidth]{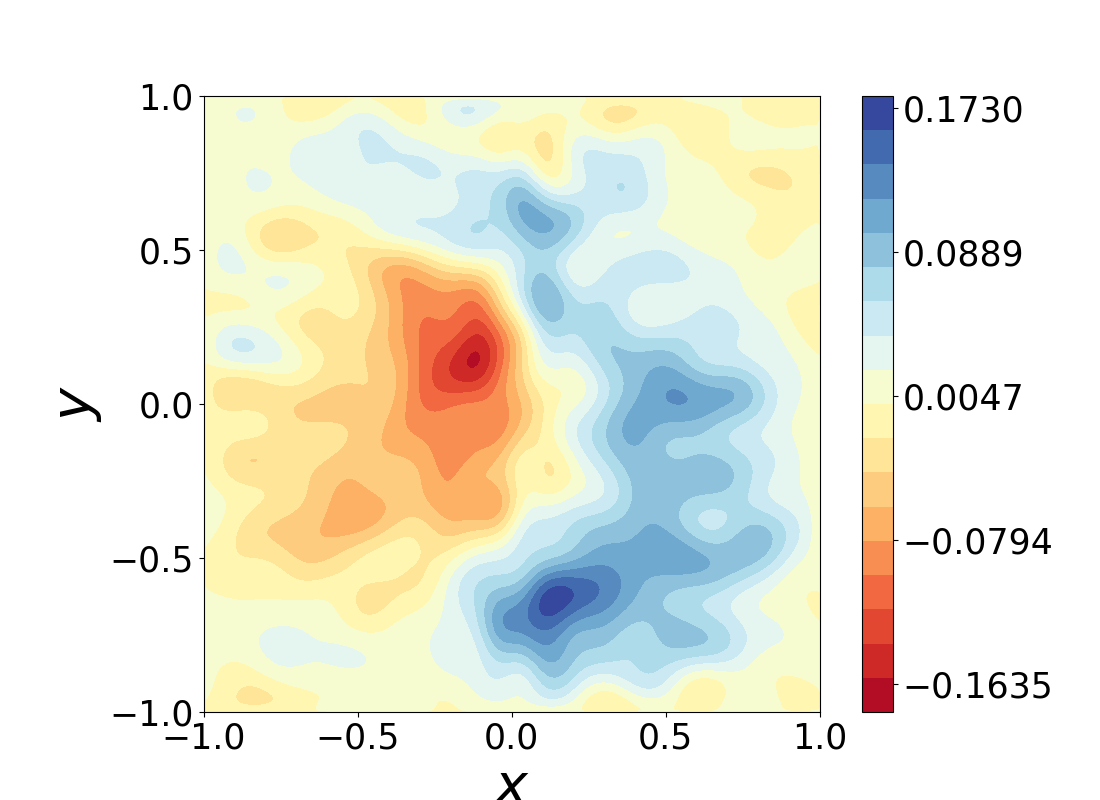} \\
    (e) $u_{\DRM}$ & (f) $u_{\DRM} - u^{\exact}$ \\   
    \end{tabular} \\
    }
\caption{Poisson equation in Sec.\ref{sec:vpinn}: (a) $u_{\WTN}$ obtained by Alg. \ref{alg:wtn}; (c) $u_{\stf}$ obtained by strong form loss; (e) $u_{\DRM}$ obtained by DRM loss; (b), (d), (f) show the pointwise errors $u_{\star} - u^{\exact}$ for $\star =\WTN, \stf$, and $\DRM$, respectively.}
\label{fig:results_vpinn}
\end{figure}

Even though the solution is smooth, it has sharp gradients near $x =0$, which makes it challenging to solve by neural network-based methods. 
Taking $M=1600$ trial basis functions with the shape parameters all being 5 and $N=1800$ test basis functions in WTN, we obtain $u_{\WTN}$. It is displayed in Fig.~\ref{fig:results_vpinn}(a) together with the associated pointwise error in Fig.~\ref{fig:results_vpinn}(b). 
We also compute the strong form solution and DRM solution using the same number of trial basis functions, the corresponding solutions and errors are shown in Fig.~\ref{fig:results_vpinn} (c-d) and (e-f), respectively. It is seen that $u_{\WTN}$ performs the best among these three methods.

\begin{figure}[!htb]
\centerline{
    \begin{tabular}{cc}
    \includegraphics[width=0.5\textwidth]{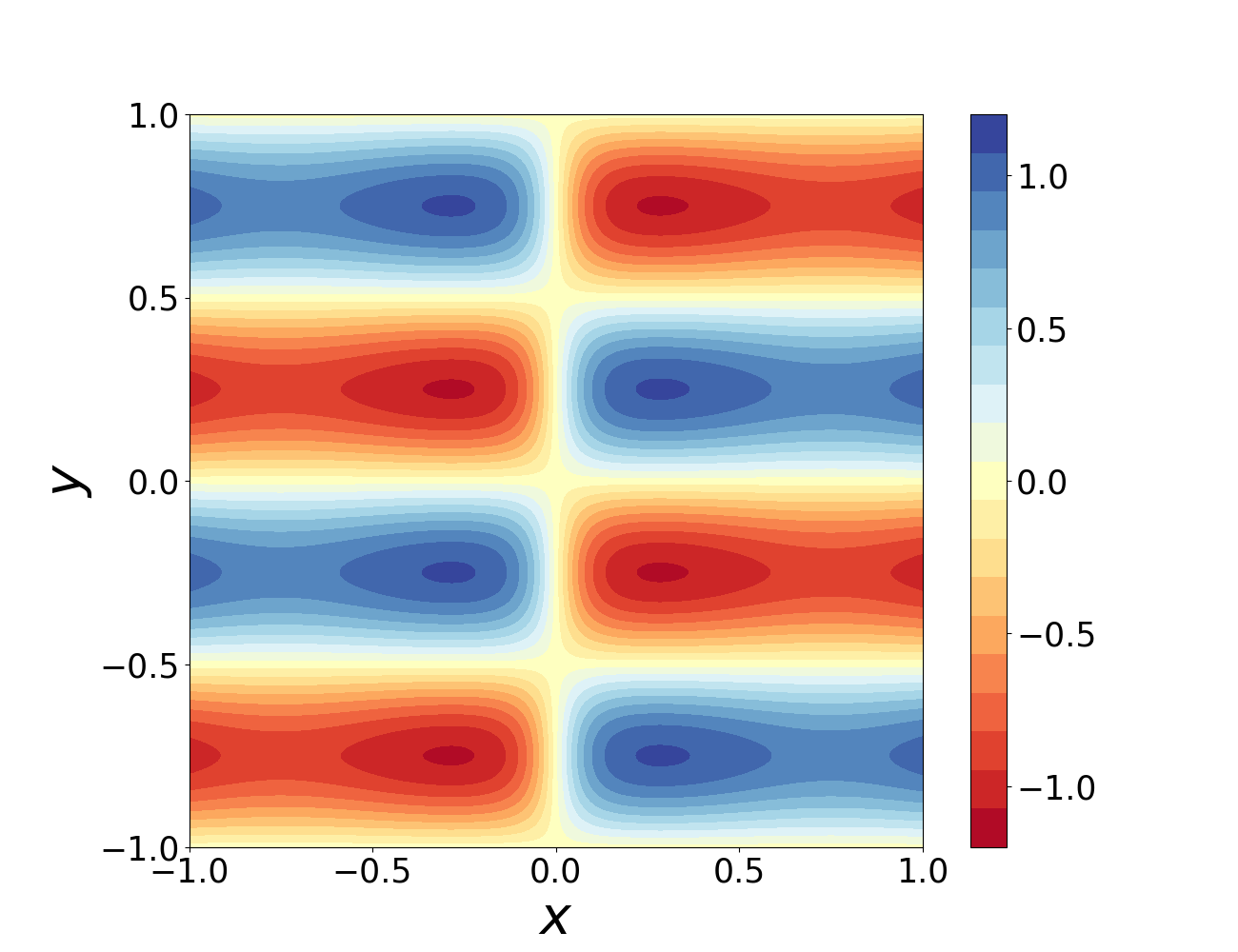} 
    &  \includegraphics[width=0.5\textwidth]{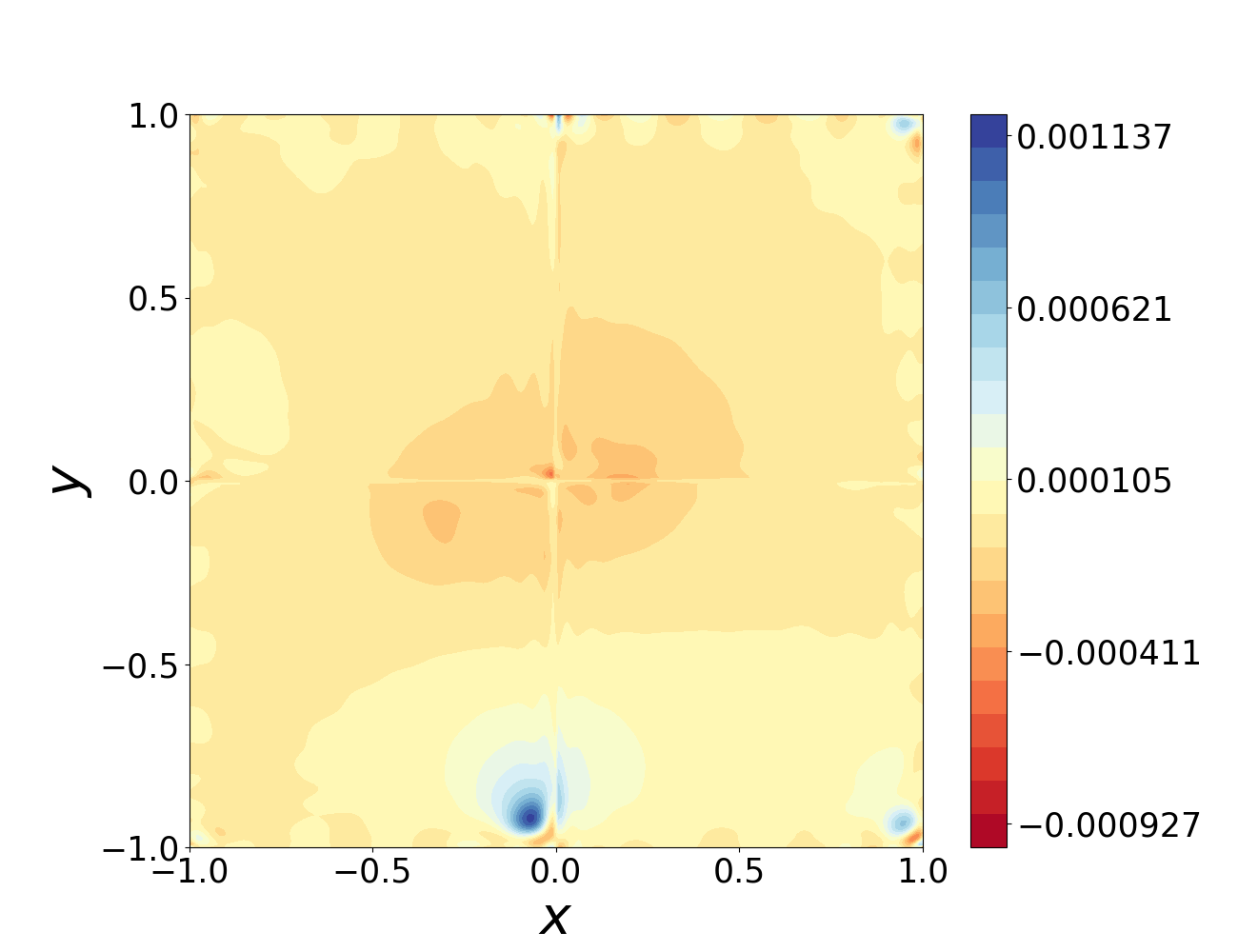} \\
    (a) $u_{\pou-\WTN}$ &  (b) $u_{\pou-\WTN}- u^{\exact}$  
    \end{tabular} \\
    }\caption{Poisson equation in Sec.~\ref{sec:vpinn}: (a) $u_{\pou-\WTN}$ obtained by Alg.~\ref{alg:pou-wtn}; and (b) pointwise error.}
\label{fig:results_vpinn_pou}
\end{figure}

 To overcome the difficulty induced by the sharp-gradient solution, we further utilize the PoU-WTN method (see Sec.~\ref{sec:pou_weaktransnet}). 
First, the domain is partitioned into 4 subdomains (see Fig.~\ref{fig:vpinn-setting} (b) for an illustration), the number of local basis in each subdomain is set to be $M^{(\ell)} = 400$ for $\ell =1,\ldots,4$, the shape parameter is set to be 5 in all neural basis functions, and the number of test functions is set as $N = 1800$. 
The boundary weight $\beta$ and the interface weight are both set to 1.
The numerical solution is shown in Fig.~\ref{fig:results_vpinn_pou}(a) along with the error displayed in Fig.~\ref{fig:results_vpinn_pou}(b). It is seen that it significantly reduces the approximation error, comparing with the WTN method. We list the relative $L_2$ errors of all these methods in Tab.~\ref{tab:re_vpinn}, which clearly demonstrate the effectiveness of the PoU-WTN method. 

\begin{table}[!htp]
	\caption{Poisson equation in Sec.~\ref{sec:vpinn}: Relative errors of $u_{\WTN}$, $u_{\stf}$, $u_\DRM$ and $u_{\pou-\WTN}$.}
    \vspace{0.1cm}
	\centering
    \adjustbox{max width=\textwidth}{
    \begin{tabular}{ccccc}
\hline
& $e_{\WTN}$ & $e_{\stf}$ & $e_{\DRM}$ & $e_{\pou-\WTN}$ \\
&$5.25\times 10^{-3} $ & $3.11\times 10^{-2}$ & $7.83\times 10^{-2}$ & $1.35\times 10^{-4}$ 
\\
\hline
\end{tabular}
}
	\label{tab:re_vpinn}
\end{table}

\subsection{L-shape problem with singularity} 
\label{sec:l_singularity}
\begin{figure}[!htb]
\centerline{
    \begin{tabular}{ccc}
    \includegraphics[width=0.4\textwidth]{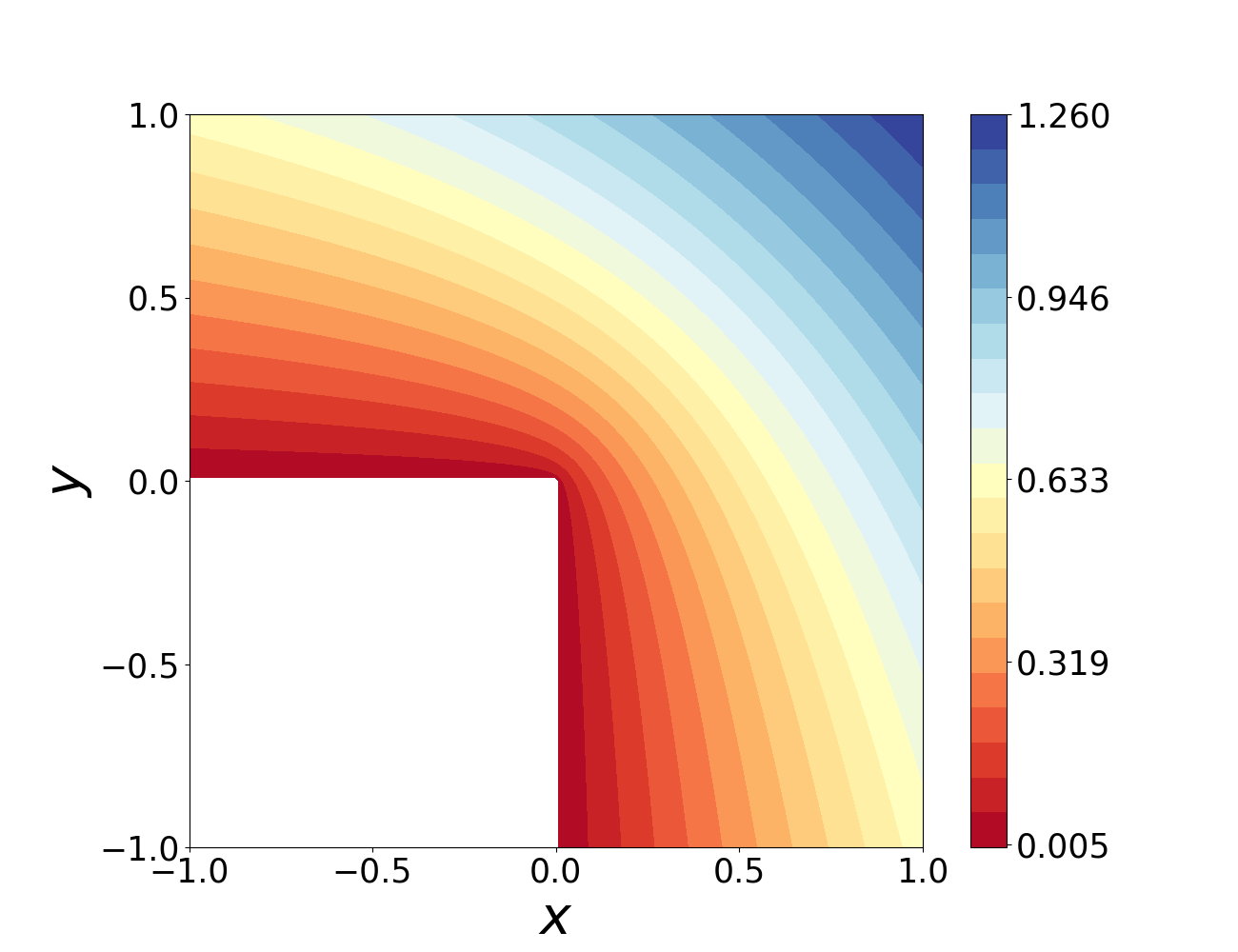}
    & 
    \includegraphics[width=0.4\textwidth]{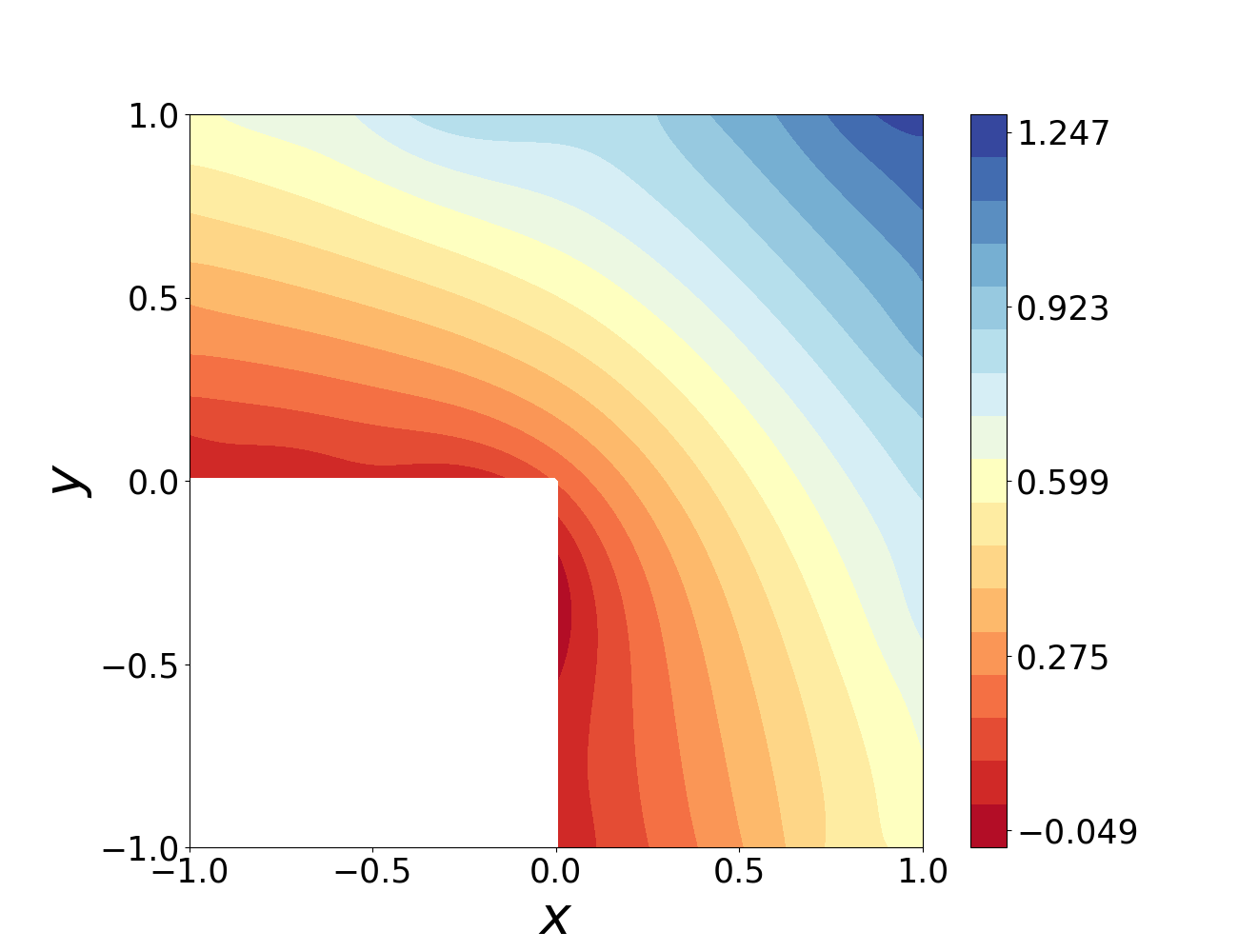} 
    & \includegraphics[width=0.4\textwidth]{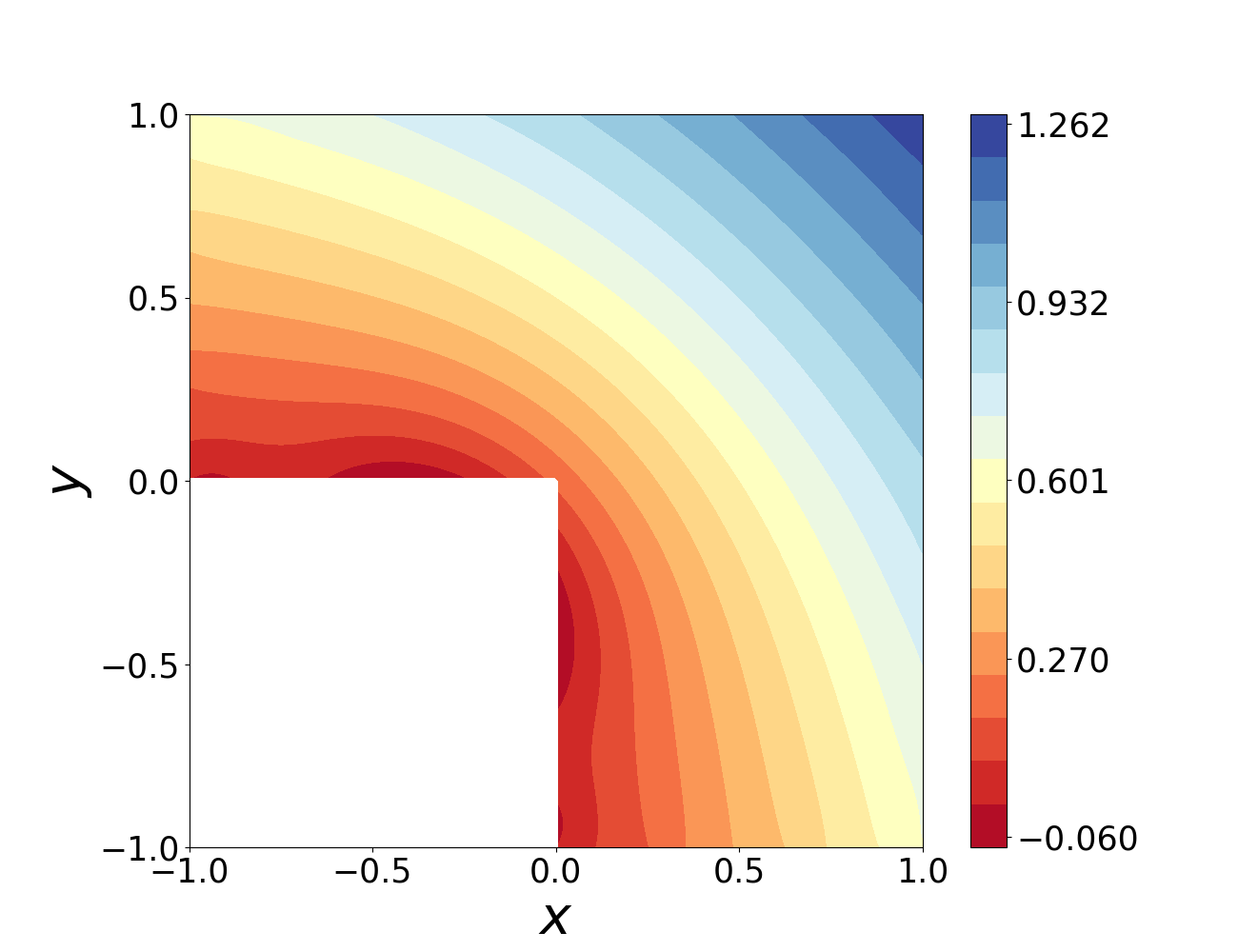}  
    \\    (a) $u^{\exact}$ & (b) $u_{\WTN}$ & (c) $u_{\stf}$
    \\
    & 
    \includegraphics[width=0.4\textwidth]{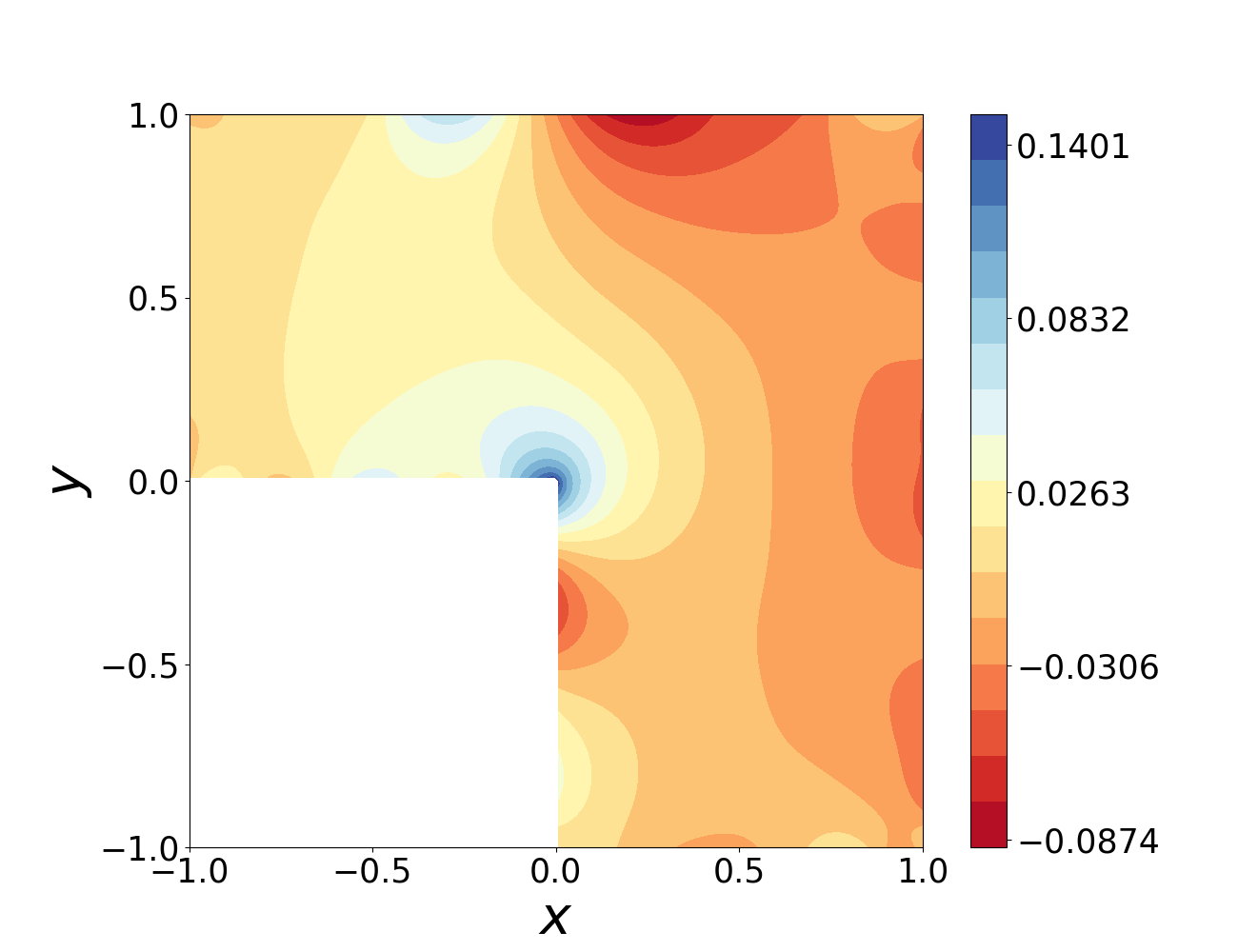} 
    & \includegraphics[width=0.4\textwidth]{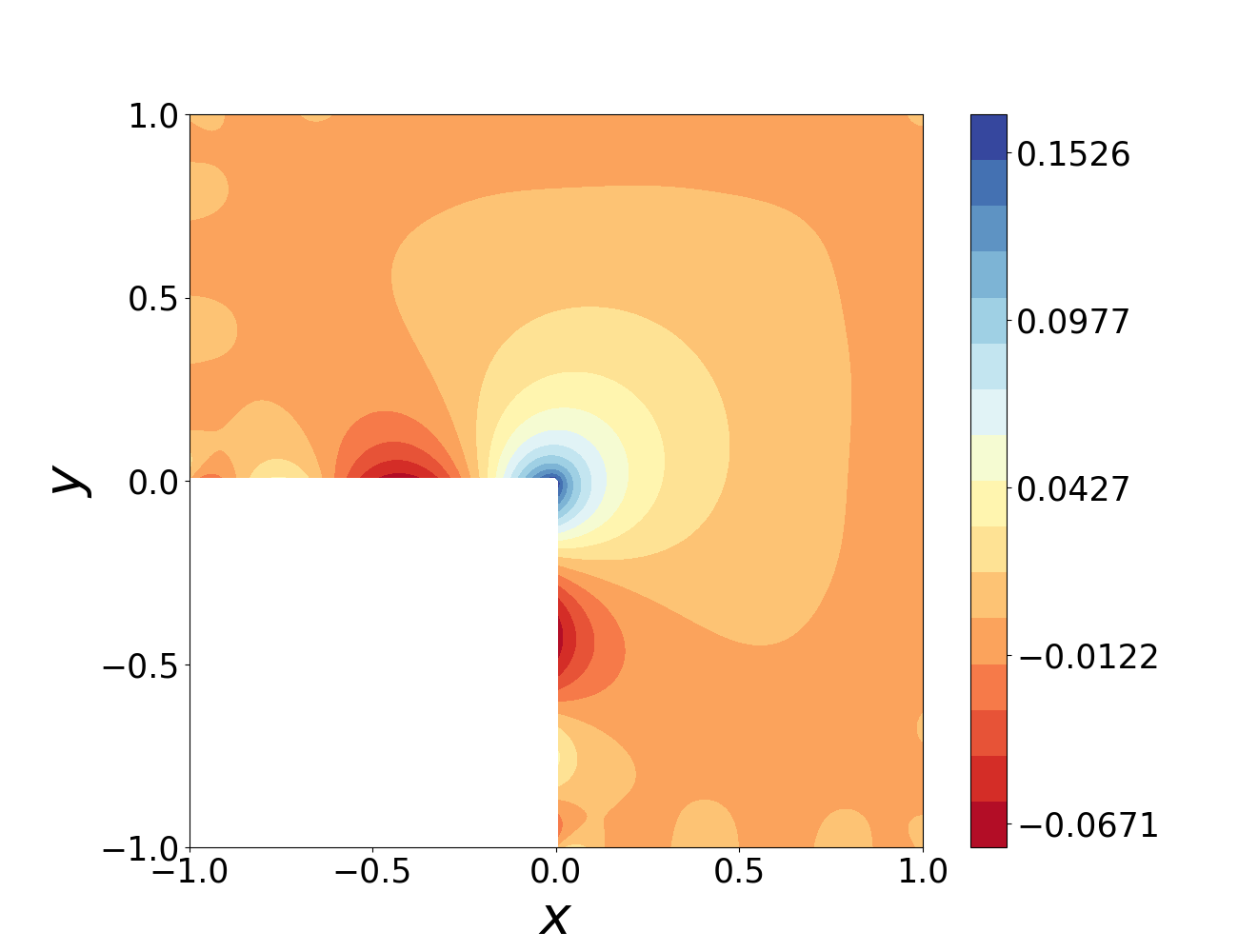}  \\
    & (d) $u_{\WTN} - u^{\exact}$ & (e) $u_{\stf} - u^{\exact}$ \\ 
    \end{tabular}  
    }
\caption{L-shape problem in Sec.~\ref{sec:l_singularity}: (a) The exact solution $u^{\exact}$; (b)$u_{\WTN}$ obtained by Alg. \ref{alg:wtn}; (c) $u_{\stf}$ obtained by strong form loss; (d) and (e) Pointwise error $u_{\star} - u^{\exact}$ for $\star =\WTN$ and $\stf$, respectively.}
\label{fig:results_L}
\end{figure}

Next, we consider again the Poisson's equation with Dirichlet boundary condition, but over an L-shape domain $\Omega = (-1,1)^2/(-1,0)^2$. 
The source term is $f = 0$, and the problem admits an analytical solution 
\[
u^{\exact}(r,\theta) = r^{2/3} \sin\left(\frac{2\theta + \pi}{3}\right)\,,
\]
where $r\ge 0$ represents the radial distance from the origin and $\theta$ denotes the angle measured counterclockwise from the positive vertical axis.  
The problem has a \textit{singularity} in the origin because the gradient of the solution becomes unbounded as $r\to 0$. The exact solution is shown in Fig.~\ref{fig:results_L}(a).

We first test SF and WTN. In both cases, we use $M=1200$ trial basis functions, where the shape parameter for all neurons are to set to 1, and $400\times 6$ boundary samples uniformly generated on the six edges. 
For the former case, $4\times 10^{4}$ interior samples are generated.
For the latter case, $N=1800$ test basis functions are selected. 
The numerical results and errors are shown in Fig.~\ref{fig:results_L} (b) and (d) for WTN, and (c) and (e) for SF, respectively.

\begin{figure}[htbp!]
\centerline{
    \begin{tabular}{cc}
    \includegraphics[width=0.27\textwidth]{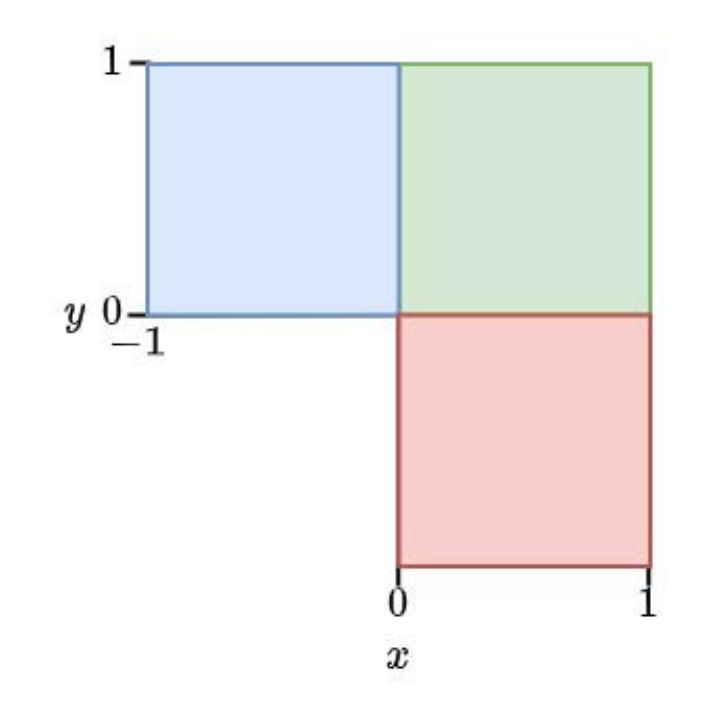}
    &
    \includegraphics[width=0.4\textwidth]{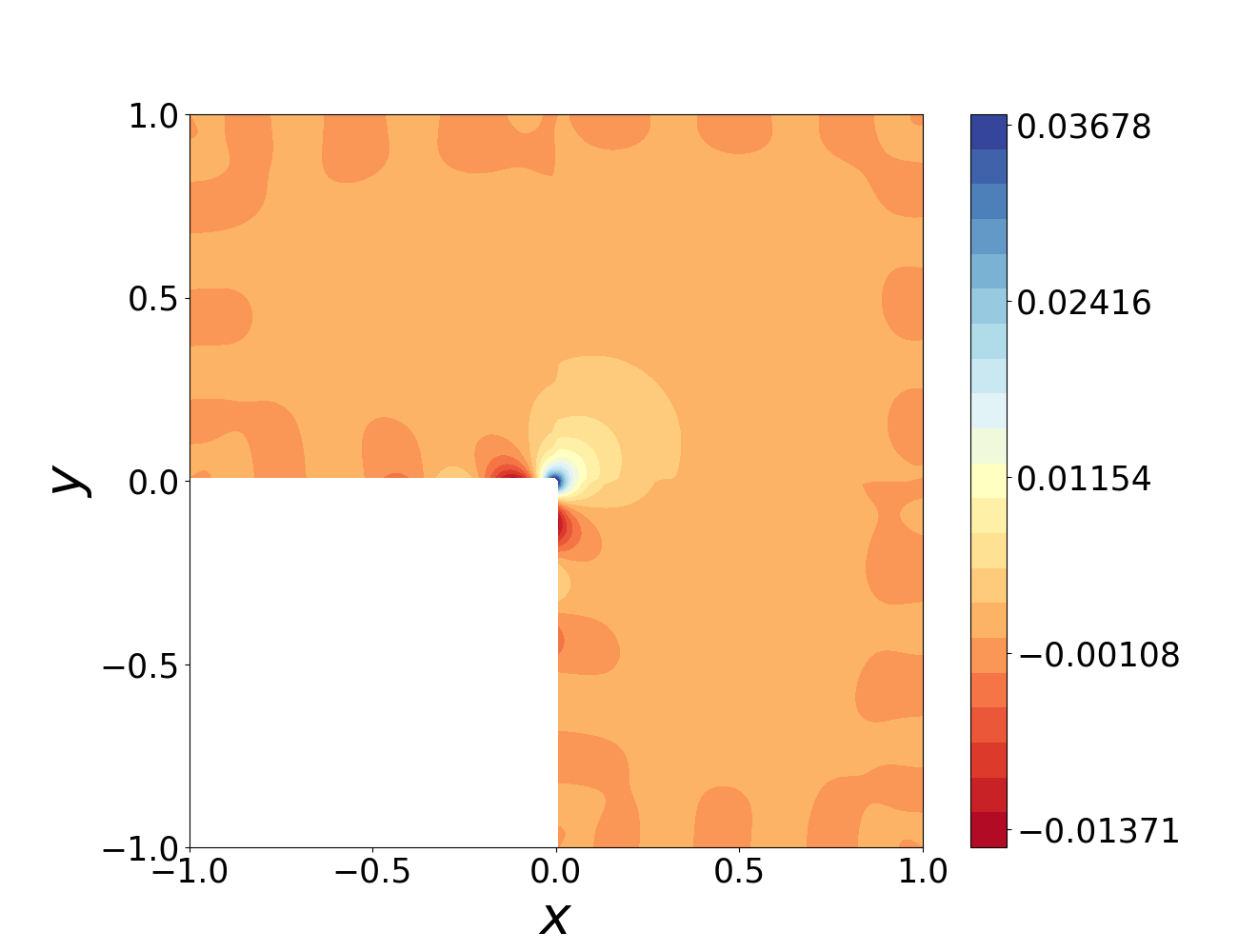} 
    \\
    (a) PoU strategy & (b) $u_{\pou-\WTN} - u^{\exact}$ \\
    \includegraphics[width=0.27\textwidth]{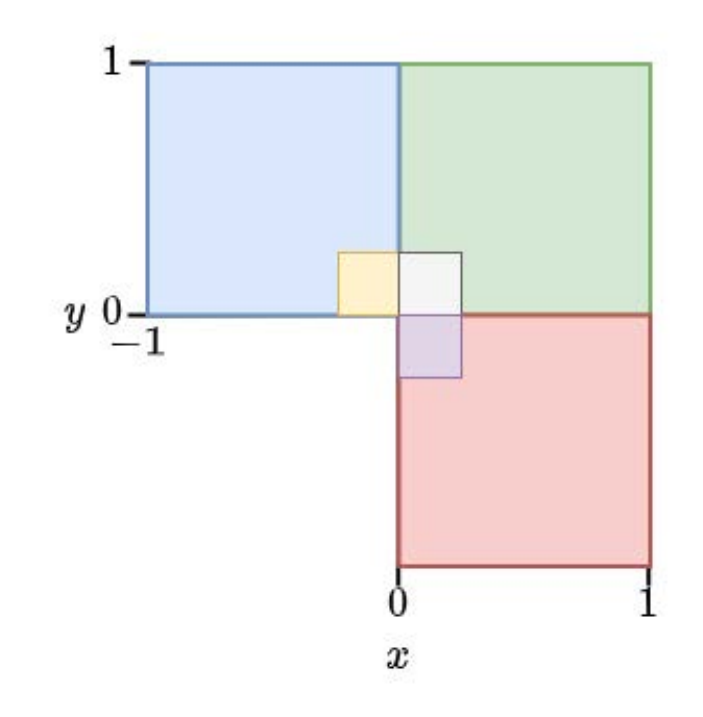}
    & 
    \includegraphics[width=0.4\textwidth]{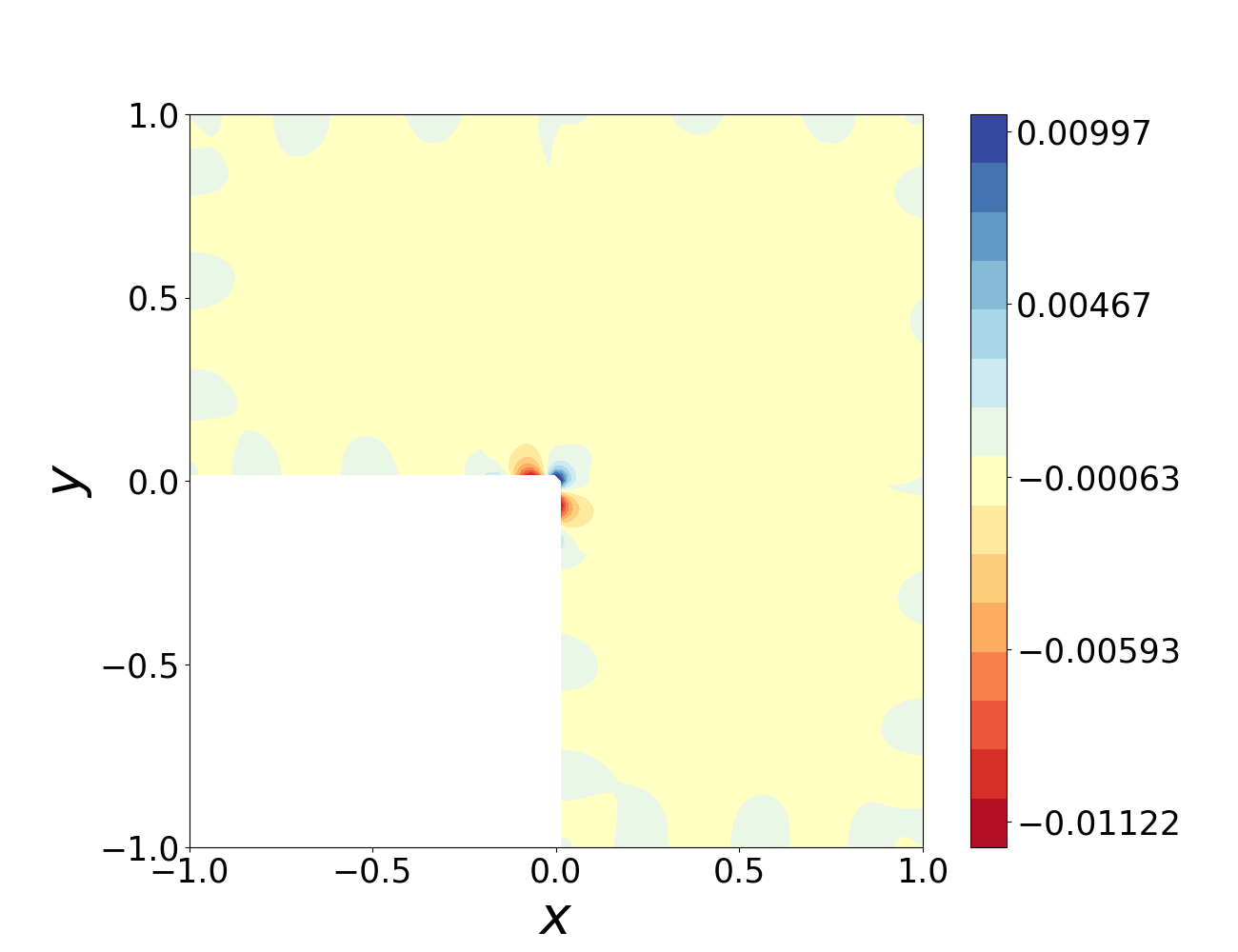} \\
    (c) PoU strategy & (d) $u_{\pou-\WTN-\R} - u^{\exact}$ \\
    \\
    \includegraphics[width=0.27\textwidth]{Figs/L_pou6_ill.pdf}
    & 
    \includegraphics[width=0.4\textwidth]{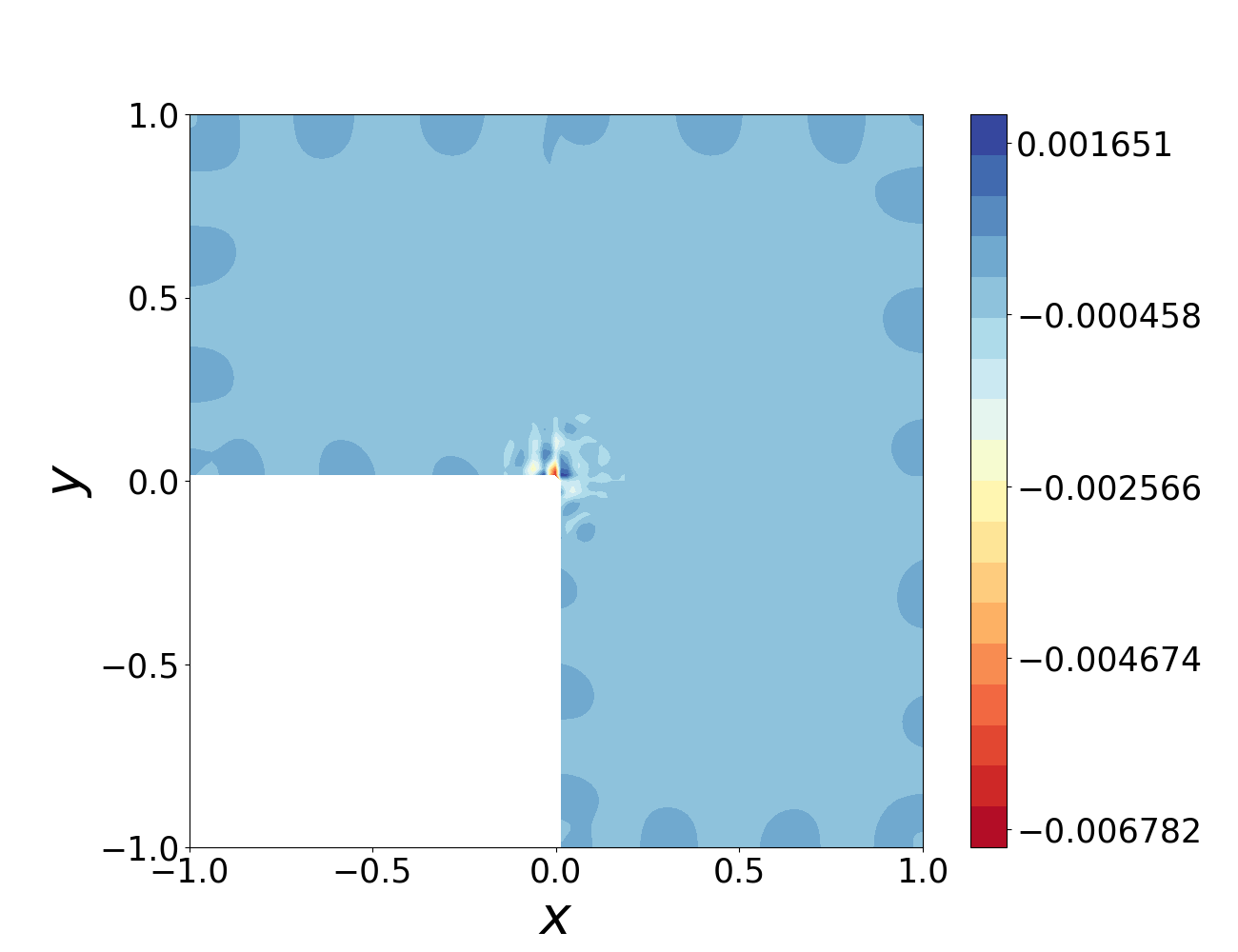}\\ 
    (e) PoU strategy & (f) $u_{\pou-\WTN-\M} - u^{\exact}$ 
    \end{tabular}
    }
    \caption{L-shape problem in Sec.~\ref{sec:l_singularity}: 
Each row shows a PoU strategy used in PoU-WTN and the associated numerical error. 
Top: The global domain is decomposed into three subdomains as shown in (a), with the shape parameter set to 1 for all neural basis functions in each subdomain. 
Middle: The global domain is decomposed into six subdomains as shown in (c), with the shape parameter set to 1 for all neural basis functions in each subdomain. 
Bottom: The PoU strategy is the same as (e), but in the three smaller subdomains, the shape parameters for neural basis functions are set to a mixture of values $(1,5,10)$.}\label{fig:results_L_pou}
\end{figure}

Observing that the error is predominantly localized near the origin,
inspired by $hp$-FEM, 
we employ the PoU-WTN and decompose the domain into three subdomains (see Fig.~\ref{fig:results_L_pou}(a)).
We set the interface weight $\lambda$ to be 1, 
and take $M^{(\ell)}= 400$ local trial basis functions with the shape parameter 1 in each subdomain. Setting the number of test basis functions to $N=1800$, we evaluate $u_{\pou-\WTN}$ and present the corresponding pointwise error in Fig.~\ref{fig:results_L_pou}(b), which becomes smaller than the error of $u_{\WTN}$, but the major error still concentrates in the small region around the origin.

Therefore, to further enhance the accuracy, we propose and investigate two new strategies. 
First, we separate out three smaller square subdomains of the size $0.2\times 0.2$ around the origin, as illustrated in Fig.~\ref{fig:results_L_pou}(c). 
In each of the six resulting subdomains, we use $M^{(\ell)} =  200$ local trial basis functions with the shared shape parameter 1, while maintaining $N=1800$ test basis functions to ensure dimensional consistency with the trial and test spaces used in Fig.~\ref{fig:results_L_pou}(a). 
The numerical solution, denoted by $u_{\pou-\WTN-\R}$, is computed, whose error is shown in Fig.~\ref{fig:results_L_pou}(d). From it, we observe that the numerical accuracy is greatly improved, indicating that the use of more local basis functions can improve the performance in this case. This is because the error primarily arises from the singularity at the origin, which motivates the next strategy.

Second, as the solution exhibits sharp gradients around the origin, it is natural to use local trial basis functions of large shape parameter in subdomains around the origin. Since there is no a prior knowledge of the optimal choice of the shape parameter, we propose a new \textit{mixing} strategy in which the shape parameters of local basis functions in the regions around the origin are set to several values. Hence, we use the same PoU functions and keep the number of trial basis functions $M^{(\ell)}$ unchanged, but in the three smaller subdomains, we divide the basis into three groups, each using a different shape parameter value from the set $\{1, 5, 10\}$.
This approach serves two key purposes. 
First, 
by including basis functions with larger shape parameter values, the resulting mixed basis can better capture steeper solution, thereby enhancing the expressive power of the local neural approximation space.
Second, as previously discussed, finding the optimal shape parameter $\gamma$ remains an open problem. Using a mixture of different shape parameters mitigates the risk associated with selecting a single, potentially suboptimal value.

We implement this strategy in the PoU-WTN framework, denoted by $u_{\pou-\WTN-\M}$ and display the numerical solution and error in Fig.~\ref{fig:results_L_pou}(e)-(f), respectively. 
Compared to Fig.~\ref{fig:results_L_pou}(e), we observe that such a treatment can better handle the singular solution and further improve the numerical accuracy. 
Tab.~\ref{tab:rel_err_L2} summarizes the relative $L_2$ errors of the solutions from $\stf$  and $\WTN$, comparing the case without PoU to three afore-discussed PoU strategies. We observe that using the PoU strategies improves both WTN and SF, with WTN achieving better accuracy. 

\begin{table}[!htp]
	\caption{L-shape problem in Sec.~\ref{sec:l_singularity}: Relative errors of solutions obtained with $\stf$ and $\WTN$, comparing the no-PoU case to three PoU strategies.}
	\vspace{0.1cm}
	\centering
    \adjustbox{max width=\textwidth}{
    \begin{tabular}{c|c|c}
    & $\WTN$ & $\stf$\\
    \hline
    No PoU & $4.30\times 10^{-2}$& $2.80\times 10^{-2}$ \\
    PoU strategy in Fig.\ref{fig:results_L_pou}(a) & $2.89\times 10^{-3}$ & $3.91\times 10^{-3}$ \\
    PoU strategy in Fig.\ref{fig:results_L_pou}(c) & $9.77\times 10^{-4}$ & $2.39\times 10^{-3}$ \\
    PoU strategy in Fig.\ref{fig:results_L_pou}(e)  & $3.06\times 10^{-4}$ & $1.54\times 10^{-3}$\\
 \hline
\end{tabular}
}
\label{tab:rel_err_L2}
\end{table}

\section{Conclusion remarks}\label{sec:conclusion} 

In this paper, we propose the WTN method for solving the elliptic PDEs, which uses neural basis functions represented by TransNet as trial functions and minimizes the weak residual of the PDE. This method follows a Petrov-Galerkin framework, in which the test space is spanned by RBFs. By leveraging the local support properties of RBFs, we can reduce the computational cost required to evaluate inner products. To address challenges in solving problems whose solutions are multiscale or have sharp gradients, we further propose the enhanced approaches: 
\begin{enumerate}
\item F-WTN, which incorporates Fourier features to WTN for finding solutions of multiple scales.
\item PoU-WTN, which uses localized neural basis functions to capture solutions of sharp gradients. 
\end{enumerate}
Through comprehensive numerical experiments, we demonstrate the effectiveness and accuracy of the proposed methods.
In particular, since a critical hyperparameter in TransNet is the shape parameter $\gamma$, which plays a key role in determining the approximability of the neural trial space, we propose using a mixture of shape parameters in local neural basis functions for PoU-WTN to better approximate the singular solution of the Poisson equation in an L-shape domain. 
In the future work, we will extend this method to multi-physics and time-dependent problems.

\section*{Acknowledgement}
Z.W. was partially supported by U.S. National Science Foundation under award numbers
DMS-2012469, DMS-2038080, DMS-2245097 and an ASPIRE grant from the Office of the Vice President for Research at the University of South Carolina.

\newpage
\appendix 
\section{Strong form loss for TransNet}
\label{sec:strong_loss}
For comparison, we briefly review the TransNet framework with strong form loss (or PINN loss~\cite{raissi2019physics}).
The loss function of the strong form, incorporating a penalty term for the boundary, is defined as 
\begin{equation}
    \mathfrak{L}_{\stf} =  \underbrace{\int_{\Omega} \|\mathcal{L}[u_{\TN}(\bx)] - f(\bx)\|^2 \dif \bx}_{\text{PDE residual}} + \underbrace{\beta_{\stf} \int_{\partial  \Omega}\|\mathcal{B}[u_{\TN}(\bx)] - g(\bx) \|^2\dif s}_{\text{boundary loss}}\,,
    \label{eq:stf_loss}
\end{equation}
where the boundary penalty term is weighted by $\beta_{\stf}$.  
There are other alternative approaches to handle boundary condition. For instance, 
the hPINNs method~\cite{lu2021physics} imposes hard constraint the penalty method and the augmented Lagrangian method. 

Given training sample sets $S_{\Omega} = \{\bx_{ \Omega}^{(m)}\}_{m=1}^{N_{\Omega}} \subset \Omega$ in the interior domain  and $S_{\partial \Omega} = \{\bx_{\partial \Omega}^{(m)}\}_{m=1}^{N_{\partial \Omega}}\subset \partial \Omega$ on the boundary , the loss \eqref{eq:stf_loss} can be approximated by the empirical risk using a Monte Carlo quadrature rule ~\cite{raissi2019physics}:
\begin{align}
    \mathfrak{L}_{\stf,\text{E}} &=  \frac{|\Omega|}{N_{\Omega}}\sum_{m=1}^{N_{\Omega}}\left[ \mathcal{L}[u_{\TN}(\bx^{(m)})] - f(\bx^{(m)})\right]^2 \nonumber\\
    &+\frac{\beta_{\stf} |\partial \Omega|}{N_{\partial \Omega}} \sum_{m=1}^{N_{\partial \Omega}}  \left[ \mathcal{B}[u_{\TN}(\bx^{(m)})] - g(\bx^{(m)}) \right]^2\,.
    \label{eq:stf}
\end{align}
For linear PDEs, it is equivalent to solving the following least-squares minimization:
\[
\balpha^{\star} = \argmin_{\balpha} \| \bm{L}_{\stf} \balpha - \bm{r}_{\stf}\|_2^2\,,
\]
where 
\[
\bm{L}_{\stf} = \begin{bmatrix}
   \bm{A}_{\stf} \\
    \widetilde{\beta} \bm{B}
\end{bmatrix}\,,\quad  \bm{r}_{\stf} = \begin{bmatrix}
    \bm{f}_{\stf} \\
    \widetilde{\beta} \bm{g}
\end{bmatrix}\,.
\]
Here, the matrix $\bm{A}_{\stf}\in \mathbb{R}^{N_{\Omega}\times (M+1)}$ with entry
 $(\bm{A}_{\stf})_{mj} = \mathcal{L}[\phi_j(\bx^{(m)})]$, 
 $\bm{f}_{\stf}\in \mathbb{R}^{N_{\Omega}}$ with entry $(\bm{f}_{\stf})_m = f(\bx^{(m)})$, $\bm{G}\in \mathbb{R}^{N_{\partial\Omega}\times (M+1)}$ with entry
$(\bm{B})_{mj} = \mathcal{B}[\phi_j(\bx^{(m)})]$, 
$\bm{g}\in \mathbb{R}^{N_{\partial \Omega}}$ with entry $\bm{g}_m = g(\bx^{(m)})$, and
$\widetilde{\beta} = \sqrt{\frac{\beta_{\stf}|\partial \Omega| N_{\Omega}}{|\Omega| N_{\partial \Omega}}}$ is the adjusted weight.   

\section{Ritz energy loss for TransNet}
\label{sec:drm_loss}
The Ritz energy is used as the loss function in the deep Ritz method (DRM)~\cite{yu2018deep}. When TransNet is used for trial, we have the loss function 
\begin{align}   
\mathfrak{L}_{\DRM}
& : = \underbrace{\mathcal{E}(u_{\TN})}_{\text{Ritz energy}} + \underbrace{\beta_{\DRM} \int_{\partial \Omega}\|\mathcal{B}[u_{\TN}(\bx)] - g(\bx) \|_2^2 \dif s }_{\text{boundary loss}}\,,\label{eq:drm_loss}
\end{align}
where $\mathcal{E}(u_{\TN})$ denotes the Ritz variational energy of $u_{\TN}$ and $\beta_{\DRM}$ is the weight of boundary loss -- a penalty to the boundary constraint.
The exact form of \eqref{eq:drm_loss} varies with the problem setting. 
Take the 2D Poisson equation for example, the DRM loss with boundary term has the following form:
\begin{align}    
\mathfrak{L}_{\DRM}
& : = \int_{\Omega}\left(\frac{1}{2}|\nabla u_{\TN}(\bx;\balpha)|^2\right)\dif \bx - \int_{\Omega}
 f(\bx)u_{\TN}(\bx;\balpha) \dif \bx  \nonumber \\
 & \quad + \beta_{\DRM} \int_{\partial \Omega} \|\mathcal{B}[u_{\TN}(\bx;\balpha)] - g(\bx) \|_2^2 \dif \bx\,.
 \label{eq:drm_loss_apdx}
\end{align}
Given the training sets 
$ S_{\Omega} = \{\bm{x}_\Omega^{(m)}\}_{m=1}^{N_{\Omega}} \subset \Omega$ and $S_{\partial \Omega} = \{\bx_{\partial \Omega}^{(m)}\}_{m=1}^{N_{\partial \Omega}}\subset \partial \Omega$, the first term in \eqref{eq:drm_loss_apdx} can be approximated by the Monte Carlo method as 
\begin{align*}
\mathfrak{L}_{\DRM}^{(1)} 
& =\int_{\Omega} \left(\frac{1}{2}|\nabla u_{\TN}(\bx;\balpha)|^2 \right)\dif \bx 
\approx \frac{|\Omega|}{N_{\Omega}}\sum_{m=1}^{N_{\Omega}} \frac{1}{2}|\nabla u_{\TN}(\bx_{\Omega}^{(m)};\balpha)|^2 \\
& = \frac{|\Omega|}{2N_{\Omega}}\sum_{m=1}^{N_{\Omega}} \left( \frac{\partial u_{\TN}(\bx;\balpha)}{\partial x}\bigg|_{\bx = \bx_{\Omega}^{(m)}}\right)^2  +  \left( \frac{\partial u_{\TN}(\bx;\balpha)}{\partial y}\bigg|_{\bx = \bx_{\Omega}^{(m)}}\right)^2\\
& = \frac{|\Omega|}{2N_{\Omega}}\left(\|\Phi_x \balpha\|_2^2 + \|\Phi_y \balpha\|_2^2 \right) = \frac{|\Omega|}{2N_{\Omega}} \left\| \begin{bmatrix}
    \Phi_x \\
    \Phi_y 
\end{bmatrix}\balpha
\right\|_2^2\,,
\end{align*}
where $\Phi_x$ is a $N_{\Omega}\times (M+1)$ matrix with entry $(\Phi_x)_{mj} = \frac{\partial \phi_j(\bx)}{\partial x}\big|_{\bx = \bx^{(m)}_{\Omega}}$,  $\Phi_y$ is a $N_{\Omega}\times (M+1)$ matrix with entry $(\Phi_y)_{mj} = \frac{\partial \phi_j(\bx)}{\partial y}\big|_{\bx = \bx^{(m)}_{\Omega}}$.
\paragraph{Case 1. $f(\bx) \equiv 0$}
In the case of the zero source function, the DRM loss \eqref{eq:drm_loss_apdx} can be rewritten as a least square problem:
\[
\mathcal{L}_{\DRM,MC} = \|\bm{L}_{\DRM}\balpha - \bm{r}_{\DRM}\|_2^2\,, 
\]
where 
\begin{equation}
\bm{L}_{\DRM} = \begin{bmatrix}
   \Phi_x \\
     \Phi_y \\
    \widetilde{\beta}\bm{B}
\end{bmatrix}\,, 
\quad 
\bm{r}_{\DRM} = \begin{bmatrix}
    \bm{0}\\
    \bm{0}\\
    \widetilde{\beta}\bm{g}
\end{bmatrix}\,.
\label{eq:drm_notation}
\end{equation}
where $\bm{B}$ and $\bm{g}$ follow the same definitions in \ref{sec:strong_loss}, and $\widetilde{\beta}\coloneqq \sqrt{\frac{2\beta_{\DRM} N_{\Omega} |\partial \Omega|}{N_{\partial \Omega}|\Omega|}|}$ is the adjusted weight.
\paragraph{Case 2. $f(\bx)$ is not a zero function}
In this case, the second term in \eqref{eq:drm_loss_apdx}:
\begin{align*}
\mathfrak{L}^{(2)}_{\DRM}
& \coloneqq \int_{\Omega}f(\bx)u_{\TN}(\bx;\balpha)\dif \bx\\
&\approx 
\frac{|\Omega|}{N_{\Omega}}\sum_{m=1}^{N_{\Omega}} f(\bx_{\Omega}^{(m)})\left( \sum_{j=0}^M \alpha_j\phi_j(\bx_{\Omega}^{(m)}) \right)=\frac{|\Omega|}{N_{\Omega}} \bm{f}^{\top}\Phi\balpha\,,
\end{align*}
where $\Phi$ is a $N_{\Omega}\times (M+1)$ matrix with entry $(\Phi)_{mj} = \phi_j(\bx_{\Omega}^{(m)})$, $\bm{f}$ is a $N_{\Omega}$-dim vector with entry $\bm{f}_m = f(\bx_{\Omega}^{(m)})$.
Then the DRM loss function is approximated by 
\[
\mathfrak{L}_{\DRM} = \|\bm{L}_{\DRM}\balpha- \bm{r}_{\DRM}\|_2^2 - \frac{|\Omega|}{N_{\Omega}} \bm{f}^{\top}\Phi\balpha 
\]
where the definitions of $\bm{L}_{\DRM}$ and $\bm{r}_{\DRM}$ follow \eqref{eq:drm_notation}.
Taking the gradient with respect to $\balpha$ as zero gives
\begin{align*}
\frac{\partial \mathfrak{L}_{\DRM}}{\partial \balpha} 
& = \frac{\partial}{\partial \balpha}  \left(\balpha^\top \bm{L}_{\DRM}^\top \bm{L}_{\DRM} \balpha - 2\bm{r}_{\DRM}^\top\bm{L}_{\DRM}\balpha - \frac{|\Omega|}{N_{\Omega}}(\bm{f}^T\Phi)\balpha\right) \\
& = 2  \bm{L}_{\DRM}^\top \bm{L}_{\DRM} \balpha - 2\bm{r}_{\DRM}^\top\bm{L}_{\DRM} - \frac{|\Omega|}{N_{\Omega}}(\bm{f}^T\Phi) = \bm{0}\,,  
\end{align*}
then the optimized $\balpha^\star$ is \footnote{When computing ($\bm{L}_{\DRM}^\top\bm{L}_{\DRM})^{-1}$, function \textsf{numpy.linalg.inv} from the \textsf{numpy} library is used. Also, a perturbation term $\epsilon \bm{I}$ with $\epsilon = 10^{-5}$ is introduced for numerical stability.}
\[
\balpha^\star = (\bm{L}_{\DRM}^\top\bm{L}_{\DRM})^{-1}\left(\bm{L}_{\DRM}^T\bm{r}_{\DRM} + \frac{|\Omega|}{2N_{\Omega}}(\Phi^T \bm{F})\right)\,.
\] 

When Darcy flow is considered, the differential operator is defined with a coefficient $\kappa(\bx)$, i.e.,  
\[
\mathcal{L}[u;\kappa] = - \nabla \cdot (\kappa \nabla u).  
\]
The first term in \eqref{eq:drm_loss_apdx}
changes to $\int_{\Omega} \left(\frac{1}{2} \kappa(\bx) |\nabla u_{\TN}(\bx;\balpha)|^2 \right) \dif \bx$, while the rest part in the loss function remains the same. Furthermore, 
the first term is approximated as  
\begin{align*}
\mathfrak{L}_{\DRM}^{(1)} 
& =\int_{\Omega} \left(\frac{1}{2} \kappa(\bx)|\nabla u_{\TN}(\bx;\balpha)|^2 \right) \dif \bx
\approx \frac{|\Omega|}{N_{\Omega}}\sum_{m=1}^{N_{\Omega}} \frac{1}{2} \kappa(\bx^{(m)}_{\Omega})|\nabla u_{\TN}(\bx^{(m)}_{\Omega};\balpha)|^2 \\
& = \frac{|\Omega|}{2N_{\Omega}}\sum_{m=1}^{N_{\Omega}} \kappa(\bx^{(m)}_{\Omega}) \left( \frac{\partial u_{\TN}(\bx;\balpha)}{\partial x}\bigg|_{\bx = \bx_{\Omega}^{(m)}}\right)^2  +  \kappa(\bx^{(m)}_{\Omega}) \left( \frac{\partial u_{\TN}(\bx;\balpha)}{\partial y}\bigg|_{\bx = \bx_{\Omega}^{(m)}}\right)^2\\
& = \frac{|\Omega|}{2N_{\Omega}}\left(\|\Phi_x(\kappa) \balpha\|_2^2 + \|\Phi_y(\kappa) \balpha\|_2^2 \right) = \frac{|\Omega|}{2N_{\Omega}} \left\| \begin{bmatrix}
    \Phi_x(\kappa) \\
    \Phi_y(\kappa) 
\end{bmatrix}\balpha
\right\|_2^2\,,
\end{align*}
where $\Phi_x(\kappa)$ is a $N_{\Omega}\times (M+1)$ matrix with entry $(\Phi_x(\kappa ))_{mj} = \sqrt{\kappa(\bx^{(m)}_{\Omega}})\frac{\partial \phi_j(\bx)}{\partial x}\big|_{\bx = \bx^{(m)}_{\Omega}}$,  and $\Phi_y(\kappa)$ is a $N_{\Omega}\times (M+1)$ matrix with entry $(\Phi_y(\kappa))_{mj} = \sqrt{\kappa(\bx^{(m)}_{\Omega}})\frac{\partial \phi_j(\bx)}{\partial y}\big|_{\bx = \bx^{(m)}_{\Omega}}$.


\bibliography{xu} 
\end{document}